\documentclass{article}

\usepackage[utf8]{inputenc}
\usepackage[T1]{fontenc}
\usepackage{mathpazo}

\RequirePackage[english]{babel}

\RequirePackage[margin=2.5cm]{geometry}
\RequirePackage{parskip}

\usepackage{tikz}
\usetikzlibrary{calc}
\usepackage{subcaption}

\makeatletter
\def\@fnsymbol#1{\ensuremath{\ifcase#1\or *\or \dagger\or \ddagger\or \mathsection\or \|\or **\or \dagger\dagger \or \ddagger\ddagger \else\@ctrerr\fi}}
\makeatother

\usepackage{xcolor}
\definecolor{cblue}{RGB}{0,70,140}
\definecolor{cgreen}{RGB}{100,140,0}
\definecolor{cred}{RGB}{190,10,50}

\usepackage[shortlabels]{enumitem}
\setlist{itemsep=0ex,topsep=0ex,parsep=0.4ex}

\usepackage{amsmath,amsfonts,amssymb,amsthm,mathtools,thmtools,thm-restate,bbm,centernot}

\usepackage[hyphens]{url} % Line breaking of urls
\usepackage[hidelinks,colorlinks,pagebackref]{hyperref} % Coloring of hyperlinks and backreferences in bibliography; must be loaded after url
\hypersetup{citecolor=cgreen,linkcolor=cblue,urlcolor=cblue}
\usepackage[capitalise,nameinlink,noabbrev,compress]{cleveref} % References in the document; must be loaded after hyperref and amsmath

\usepackage{todonotes} %Disable by adding [disable]
\definecolor{lightgreen}{rgb}{0.8,1.0,0.8}
\definecolor{lightblue}{rgb}{0.8,0.85,1}

\usepackage{csquotes} % Nice quotation marks

\renewcommand*{\backref}[1]{}
\renewcommand*{\backrefalt}[4]{
	\ifcase #1 Not cited.%
	\or $\uparrow$#2%
	\else $\uparrow$#2%
	\fi%
}

\theoremstyle{plain}
\newtheorem{theorem}{Theorem}[section]
\newtheorem{lemma}[theorem]{Lemma}
\newtheorem{claim}{Claim}

\newtheorem{observation}[theorem]{Observation}
\newtheorem{corollary}[theorem]{Corollary}
\newtheorem{conjecture}[theorem]{Conjecture}

\theoremstyle{definition}
\newtheorem{definition}[theorem]{Definition}
\newtheorem{question}[theorem]{Question}

\makeatletter
\renewenvironment{proof}[1][\proofname]
{\par\pushQED{\qed}
	\normalfont\topsep6\p@\@plus6\p@\relax\trivlist
	\item[\hskip\labelsep\bfseries#1\@addpunct{.}]
	\ignorespaces}
{\popQED\endtrivlist\@endpefalse}
\makeatother
\let\emptyset\varnothing

\newcommand{\BackedgeGraphUnlabeled}[2]{
\begin{tikzpicture}[scale=0.6]
  \tikzset{v/.style = {shape=circle, fill=black, draw, minimum size=5pt, inner sep=0pt}}
  \tikzset{edge/.style = {thick}}
  \foreach \i in {1,...,#2} { \node[v] (v\i) at (\i-1,0) {}; }
  \foreach \i/\j in {#1} { \draw[edge] (v\i) to[bend left=55] (v\j); }
\end{tikzpicture}
}

\newcommand{\BackedgeGraphLabeled}[3]{
\begin{tikzpicture}[scale=0.6]
  \tikzset{v/.style = {shape=circle, fill=black, draw, minimum size=5pt, inner sep=0pt}}
  \tikzset{edge/.style = {thick}}

  \foreach \i in {1,...,#2} {
    \node[v] (v\i) at (\i-1,0) {};
  }%
  \foreach \lab [count=\i from 1] in {#3} {
    \node[below=3pt] at (v\i.south) {\lab};
  }
  \foreach \i/\j in {#1} {
    \draw[edge] (v\i) to[bend left=55] (v\j);
  }
\end{tikzpicture}
}
\newcommand{\BackedgeGraphGeneral}[4]{
\begin{tikzpicture}[scale=0.6]
  \tikzset{v/.style = {shape=circle, fill=black, draw, minimum size=5pt, inner sep=0pt}}
  \tikzset{edge/.style = {thick}}
  \tikzset{lbl/.style = {below=5pt, font=\scriptsize}} 

  \foreach \d in {#4} {
    \global\expandafter\let\csname showdotat\d\endcsname\relax
  }

  \foreach \lab [count=\k] in {#3} {
    \global\expandafter\def\csname storedlab\k\expandafter\endcsname\expandafter{\lab}
  }

  \xdef\currentx{0}

  \foreach \i in {1,...,#2} {
    \node[v] (v\i) at (\currentx,0) {};
    
    \ifcsname storedlab\i\endcsname
        \node[lbl] at (v\i.south) {\csname storedlab\i\endcsname};
    \fi

    \ifcsname showdotat\i\endcsname
        \node at (\currentx + 1, 0) {$\cdots$};
        \pgfmathsetmacro{\nextx}{\currentx + 2}
    \else
        \pgfmathsetmacro{\nextx}{\currentx + 1}
    \fi
    
    \xdef\currentx{\nextx}
  }

  \foreach \i/\j in {#1} {
    \draw[edge] (v\i) to[bend left=55] (v\j);
  }

  \foreach \d in {#4} { \global\expandafter\let\csname showdotat\d\endcsname\undefined }
  \foreach \k in {1,...,#2} { \global\expandafter\let\csname storedlab\k\endcsname\undefined }
\end{tikzpicture}
}

\newcommand{\BackedgeGraphWithEllipsis}[4]{%
\begin{tikzpicture}[scale=0.7,
  v/.style   = {circle,fill=black,draw,minimum size=5pt,inner sep=0pt},
  edge/.style= {thick},
  lbl/.style = {below=5pt,font=\scriptsize}
]

  \foreach \lab [count=\k] in {#3} {
    \global\expandafter\def\csname storedlab\k\expandafter\endcsname\expandafter{\lab}
  }
  \def\gapindex{#4} % after which visible vertex the gap starts
  \def\gapwidth{2}  % horizontal length of the "omitted" block
  \def\dx{1}        % spacing between visible vertices

  \pgfmathsetmacro{\currentx}{0}
  \xdef\currentx{\currentx}

  \foreach \i in {1,...,#2} {%
    \node[v] (\i) at (\currentx,0) {};

    \ifcsname storedlab\i\endcsname
      \node[lbl] at (\i.south) {\csname storedlab\i\endcsname};
    \fi

    \ifnum\i=\gapindex
      \pgfmathsetmacro{\mid}{\currentx + 0.75*\gapwidth}
      \node at (\mid,0) {$\cdots$};

      \coordinate (EL) at ($(\currentx+0.65*\gapwidth,0.0)+(0,0.3)$);
      \coordinate (ER) at ($(\currentx+0.8*\gapwidth,0.0)+(0,0.3)$);

      \pgfmathsetmacro{\nextx}{\currentx + \gapwidth + \dx}
    \else
      \pgfmathsetmacro{\nextx}{\currentx + \dx}
    \fi

    \xdef\currentx{\nextx}
  }

  \foreach \a/\b in {#1} {%
    \draw[edge] (\a) to[bend left=55] (\b);
  }

  \foreach \k in {1,...,#2} {%
    \global\expandafter\let\csname storedlab\k\endcsname\undefined
  }

\end{tikzpicture}%
}

\title{Excluded Subtournaments for Matching Tournaments}
\author{Mathieu Rundstr\"om, Sophie Spirkl\thanks{We acknowledge the support of the Natural Sciences and Engineering Research Council of Canada (NSERC), [funding reference number RGPIN-2020-03912]. Cette recherche a \'et\'e financ\'ee par le Conseil de recherches en sciences naturelles et en g\'enie du Canada (CRSNG), [num\'ero de r\'ef\'erence RGPIN-2020-03912]. This project was funded in part by the Government of Ontario. This research was completed while Spirkl was an Alfred P. Sloan Fellow.  This research was undertaken, in part, thanks to funding from the Canada Research Chairs Program.}, Dinis Vitorino\thanks{Current address: Computer Science Institute, Charles University. The majority of the research involved in this paper was conducted while Vitorino was a Master's student at the David R. Cheriton School of Computer Science of the University of Waterloo},  Yun Xing\\University of Waterloo}

\date{\today}

\begin{document}

\maketitle

\begin{abstract}
    For an ordering $\sigma$ of the vertices of a tournament $T$, the backedge graph $B_\sigma(T)$ is the ordered graph whose edges are precisely the arcs of $T$ directed backwards with respect to $\sigma$.
    We study tournaments that admit an ordering whose backedge graph avoids prescribed ordered graphs, and provide excluded-subtournament characterizations of the resulting classes.
    Our main result is a characterization of the minimal tournaments that do not admit a matching backedge graph. We also give excluded-subtournament characterizations for some other classes of tournaments.
\end{abstract}

\section{Introduction}
\label{sec: intro}
A \emph{tournament} is an orientation of a complete graph. Given a tournament $T$ and a linear order $\sigma$ of its vertex set, one can define the \emph{backedge graph} of $T$ with respect to $\sigma$ (which we will also call the backedge graph of $T$ given by $\sigma$), denoted by $B_\sigma(T)$ ($B_{\sigma}$ when the tournament in question is unambiguous), as the undirected ordered graph with an edge between $v$ and $w$ if and only if $vw$ is an arc of $T$ and $w$ precedes $v$ in $\sigma$.
Backedge graphs are a convenient way to represent tournaments, as the backedge graph of a tournament with fewest edges is often substantially sparser than the original tournament.
Additionally, the set of backedge graphs of a given tournament exhibits useful closure properties: For example, if we consider isomorphic backedge graphs to be the same (forgetting about the orderings on the vertex set), the set of backedge graphs of a given tournament is closed under taking graph complements.
This sometimes results in backedge graphs acting as a bridge between the study of directed and undirected graphs.
The earliest application of backedge graphs that we are aware of is in the proof of Thm.~1.2 of \cite{erdos_hajnal_tournament_equivalence}, i.e., the proof of equivalence of the graph and tournament versions of the Erdős-Hajnal conjecture.
Similar ideas also seem to play a role in the proof of Rao's conjecture \cite{CHUDNOVSKY201444}.
More recent papers where backedge graphs are studied more explicitly include \cite{kimStructureD122025,aubian2024computing,aboulker2023cliquenumbertournaments,crew2026characterizinglargecliquenumber,CHUDNOVSKY2024103786}.
In particular, \cite{kimStructureD122025} gives a structure theorem for all tournaments excluding a certain five-vertex tournament, and the structure theorem is formulated in terms of backedge graphs. 

In this paper, we are concerned with the problem of determining the minimal tournaments (w.r.t.\ subtournaments) with no backedge graph in some hereditary classes of ordered graphs $\mathcal{C}$ (where hereditary means closed under induced ordered subgraphs).
Such characterizations are related to the decision problem of determining whether a tournament has a backedge graph isomorphic to a graph in a hereditary class $\mathcal{C}$ of ordered graphs.
Aboulker et al. \cite{aboulker_finding_2026} denote said problem\footnote{Actually, \cite{aboulker_finding_2026} considers classes $\mathcal{C}$ of \emph{unordered} graphs. This is easily seen to be a special of considering ordered graph classes: Instead of the class $\mathcal{C}$ of unordered graphs, we could instead consider the class of all ordered graphs whose underlying unordered graphs are in $\mathcal{C}$.} by $\mathcal{C}$-FAS (after \emph{feedback arc set}).
If the $\mathcal{C}$-FAS problem is NP-complete for some class $\mathcal{C}$, then $\mathcal{C}$ is unlikely to admit a finite forbidden-subtournament characterization, since one could then test membership in polynomial time by brute force over the fixed obstruction set.
The contrapositive of this statement informs us that exhibiting a finite list of forbidden tournaments gives a (trivial) polynomial time algorithm for solving the $\mathcal{C}$-FAS problem.
Similarly, whenever it can be checked whether a tournament belongs to the class of minimal forbidden sub-tournaments for $\mathcal{C}$ in polynomial time, the $\mathcal{C}$-FAS problem is in the intersection of NP and CoNP, and therefore unlikely to be NP-hard.

Aboulker et al. \cite{aboulker_finding_2026} show that $\mathcal{C}$-FAS is NP-complete for the class of forests.
Aubian shows that $\mathcal{C}$-FAS is NP-complete for the class $\mathcal{C}$ of graphs with clique number at most $k$, for all $k\geq 3$, and conjectures the same is true for $k=2$.
In forthcoming work with David Evangelista, Sahab Hajebi, Hidde Koerts, and Eileen Robinson, the first, second, and third authors prove that $\mathcal{C}$-FAS is NP-complete for $\mathcal{C}$ as above and $k=2$. Additionally, Thm. 1.6 of \cite{hardness2Colorability} states that $\mathcal{C}$-FAS is NP-complete for $\mathcal{C}$ the class of bipartite graphs.

The $\mathcal{C}$-FAS problem can trivially be solved in polynomial time for $\mathcal{C}$ the class of graphs with no edges, since it is equivalent to determining whether a tournament contains a directed triangle (i.e., whether it is transitive).
The second author, in collaboration with David Evangelista (see \cite{evangelista2026backedge}), has also shown that one can check in polynomial time whether a tournament admits a backedge graph which is a disjoint union of cliques. 
Davot et al. \cite{davot_degreewidth_2023} define the \emph{degreewidth} of a tournament as the minimum maximum degree over all backedge graphs.
They show that one can determine if a tournament has degreewidth at most one in time $O(n^3)$.

\begin{figure}[ht]
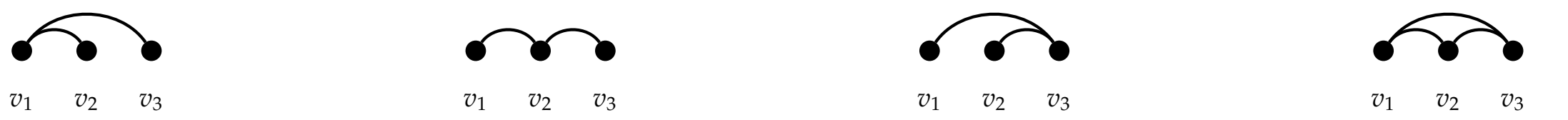

    \centering
    
    \begin{subfigure}{0.24\linewidth}\centering
    \BackedgeGraphGeneral{1/2,1/3}{3}{$v_1$, $v_2$, $v_3$}{}
    \caption{$B_1$}
    \end{subfigure}\hfill
    \begin{subfigure}{0.24\linewidth}\centering
    \BackedgeGraphGeneral{1/2,2/3}{3}{$v_1$, $v_2$, $v_3$}{}
    \caption{$B_2$}
    \end{subfigure}\hfill
    \begin{subfigure}{0.24\linewidth}\centering
    \BackedgeGraphGeneral{1/3,2/3}{3}{$v_1$, $v_2$, $v_3$}{}
    \caption{$B_3$}
    \end{subfigure}\hfill
    \begin{subfigure}{0.24\linewidth}\centering
    \BackedgeGraphGeneral{1/2,1/3,2/3}{3}{$v_1$, $v_2$, $v_3$}{}
    \caption{$B_4$}
    \end{subfigure}\hfill

    \caption{The four minimal (ordered) backedge graphs with minimum degree at least two}\label{fig:forb_backedge_graphs_matching}
\end{figure}

Following the terminology of \cite{liuVariousTheoremsTournaments2012}, we say a tournament is a \emph{matching tournament} if it has a backedge graph of maximum degree one, i.e., we call tournaments with degreewidth one matching tournaments.
This class of tournaments has also been referred to as the class of \emph{sparse tournaments}, cf. \cite{davot_degreewidth_2023,bessyArcdisjointCyclePacking2018}, but we will generally avoid this terminology.
A tournament is matching if and only if it has a backedge graph avoiding the four ordered graphs below (as ordered induced subgraphs); these are precisely the minimal ordered graphs containing a vertex of degree two.
More specifically, consider $v_1 \prec v_2 \prec v_3$.
Let $B^1, B^2, B^3$, $B^4$ be backedge graphs with vertex set $\{v_1,v_2,v_3\}$ and the following edge sets $E(B^1) = \{v_1v_2, v_1v_3\}$, $E(B^2) = \{v_1v_2, v_2v_3\}$, $E(B^3) = \{v_1v_3, v_2v_3\}$, and $E(B^4)=\{v_1v_2,v_1v_3,v_2v_3\}$, cf. \Cref{fig:forb_backedge_graphs_matching}.
Then a tournament $T$ is a matching tournament if and only if $T$ has a backedge graph $B_\sigma$ that does not contain $B^i$ for $i\in [4]$ as an (induced) ordered subgraph.
This phrasing of the matching property is more cumbersome, but not pointless: it lends itself to a natural generalization of forbidding ordered graphs in the backedge graphs of tournaments.
Instead of forbidding all four $B^i$, we can ask about forbidding just $B^1,B^3,B^4$; this gives the class of \emph{paving} tournaments from \cite{kimStructureD122025}. 

Augmenting a result in \cite{liuVariousTheoremsTournaments2012}, where it was shown that there are infinitely many minimal forbidden subtournaments for membership in the class of matching tournaments, we give a complete characterization of the minimal forbidden subtournaments for membership in this class in \Cref{section:main}.
It is clear from our characterization that one can check whether a tournament contains a minimal forbidden subtournament for this class in polynomial time, and this result can therefore be seen as resolving some of the tension between the aforementioned results in \cite{liuVariousTheoremsTournaments2012} (there are infinitely many forbidden subtournaments for the class of matching tournaments) and \cite{davot_degreewidth_2023} (one can check whether a tournament is a matching tournament in polynomial time).
We additionally give characterizations of the forbidden subtournaments for the class of paving tournaments in \Cref{section:paving}, and for the class of matching tournaments having a backedge graph where the edges do not cross in \Cref{section:NCM}.

\section{Preliminaries}
\label{sec: prelims}

Given a tournament $T$ and two disjoint subsets $X$ and $Y$ of $V(T)$, we write $X\Rightarrow Y$ if all arcs between $X$ and $Y$ have a tail in $X$ and head in $Y$. If $X$ or $Y$ are singletons, we may write the element instead of the set in this notation, e.g., $v\Rightarrow Y$. Given an ordering $\sigma$ of the vertices of $T$, we write $w\prec_{\sigma} v$ to mean that $\sigma(w) < \sigma(v)$. 

The \emph{mirror} of a backedge graph $B$ of a tournament $T$ with respect to the ordering $\sigma$ is the backedge graph $B^*$ of the tournament $T^*$, obtained by inverting each arc of $T$, with respect to $\sigma^*$, the inverse of $\sigma$ (i.e., the ordering with $v\prec_{\sigma^*} w$ if and only if $w\prec_{\sigma} v$). Note that $B$ is obtained from $B^*$ by reversing the order of the vertices. 

An ordered graph $H$ is contained in an ordered graph $G$ if $G$ has an order-preserving induced copy of $H$.  Following \cite{berger2013tournaments}, for tournaments $A,B,C$, let $\Delta(A,B,C)$ denote the tournament obtained from a directed triangle by substituting $A,B,C$ for its three vertices, with $A\Rightarrow B\Rightarrow C\Rightarrow A$.  We write $\Delta(a,b,c)$ when $A,B,C$ are transitive tournaments of orders $a,b,c$. 

An \emph{induced monotone path} is an ordered graph with with ordered vertex set $v_1,\dots,v_n$ and edges $v_{i}v_{i+1}$ for all $i\in[1,n-1]$.
A \emph{paving tournament}, introduced in \cite{kimStructureD122025} to characterize $\Delta(1,2,2)$--free tournaments, is a tournament that has a backedge graph whose connected components are monotone paths.
% where each vertex has at most one neighbour before in the ordering, earlier at most one neighbour later in the ordering. The connected components in such a backedge graph are called \emph{monotone paths}.
In \Cref{section:paving}, we characterize the minimal non-paving tournaments by exhibiting a finite list of minimal obstructions.

\begin{figure}[ht]
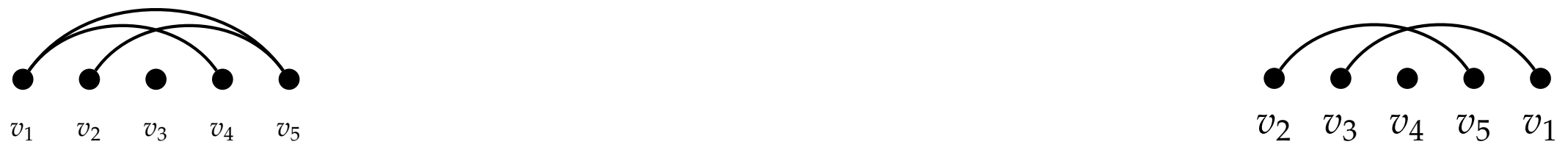

    \label{fig: backedge graphs of weird finite things minimal nonmatching}
    \centering
    \begin{subfigure}{0.32\linewidth}\centering
    \BackedgeGraphGeneral{1/4,1/5,2/5}{5}{$v_1$, $v_2$, $v_3$, $v_4$, $v_5$}{}
    \caption{$T_5$}\label{subfig:T5}
    \end{subfigure}\hfill
    \begin{subfigure}{0.32\linewidth}\centering
    \BackedgeGraphLabeled{1/4,2/5}{5}{$v_2$,$v_3$,$v_4$,$v_5$,$v_1$}
    \caption{$U_5$}\label{subfig:U5}
    \end{subfigure}
    %
    % \begin{subfigure}{0.32\linewidth}\centering
    % \BackedgeGraphUnlabeled{1/4,1/5,2/6}{6}
    % \caption{$U_5^{(5)}$}\label{subfig:1-in-degen-U55}
    % \end{subfigure}\hfill
    %
    %
    % \begin{subfigure}{0.48\linewidth}\centering
    % \BackedgeGraphGeneral{1/5, 2/6, 3/7, 3/4, 4/5}{7}{$w_1$, $w_2$, $v_1$, $v_2$, $v_3$, $w_3$, $w_4$}{}
    % \caption{$M_7$}\label{subfig:M7}
    % \end{subfigure}
    % \begin{subfigure}{0.48\linewidth}\centering
    % % \BackedgeGraphUnlabeled{1/4,1/6,2/5,3/6}{6}
    % \BackedgeGraphUnlabeled{1/3,1/6,2/5,4/6,3/4}{6}
    % \caption{$P_7^-$}\label{subfig:1-in-degen-P7^-3}
    % \end{subfigure}
    \caption{Backedge graphs of $T_5$ and $U_5$.}
\end{figure}

\begin{figure}[ht]
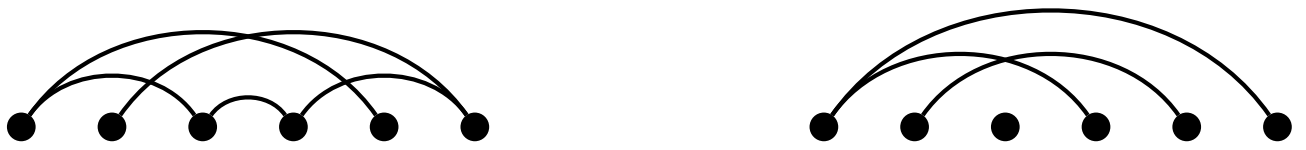
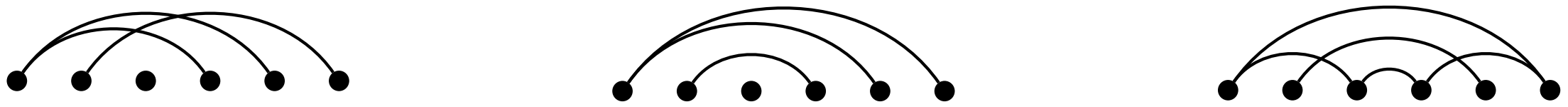

    \centering
    
    \begin{subfigure}{0.32\linewidth}\centering
    \BackedgeGraphUnlabeled{1/3,1/5,2/6,4/6,3/4}{6}
    \caption{$\Delta(2,2,2)$}\label{subfig:1-in-degen-D222}
    \end{subfigure}\begin{subfigure}{0.32\linewidth}\centering
    \BackedgeGraphUnlabeled{1/4,1/6,2/5}{6}
    \caption{$Q_6$}\label{subfig:1-in-degen-Q6}
    \end{subfigure}
    
    \begin{subfigure}{0.32\linewidth}\centering
    \BackedgeGraphUnlabeled{1/4,1/5,2/6}{6}
    \caption{$U_5^{(5)}$}\label{subfig:1-in-degen-U55}
    \end{subfigure}\hfill
    \begin{subfigure}{0.32\linewidth}\centering
    \BackedgeGraphUnlabeled{1/5,1/6,2/4}{6}
    \caption{$\Delta(1,C_3,2)$}\label{subfig:1-in-degen-D1C32}
    \end{subfigure}\hfill\begin{subfigure}{0.32\linewidth}\centering
    % \BackedgeGraphUnlabeled{1/4,1/6,2/5,3/6}{6}
    \BackedgeGraphUnlabeled{1/3,1/6,2/5,4/6,3/4}{6}
    \caption{$P_7^-$}\label{subfig:1-in-degen-P7^-}
    \end{subfigure}
    
    \caption{Backedge graphs of the five minimal non-1-in-degenerate tournaments on six vertices.}\label{fig:1-in-degen}
\end{figure}

Our proofs of statements about the sets of minimal nonpaving and minimal nonmatching tournaments rely on a case-distinction depending on the \emph{degeneracy} of the tournaments under considerations.
For a non-negative integer $k$, a tournament $T$ is \emph{$k$-in-degenerate} if every subtournament of $T$ has a vertex of in-degree at most $k$; a tournament being \emph{$k$-out-degenerate} is defined analogously.
We denote by $d_T^-(v)$ and $d_T^+(v)$ the in-degree and out-degree of a vertex $v$ in $T$.
Previously, in unpublished work joint with David Evangelista and Xinyue Fan (see also \cite{xing2026excluded}), the third and fourth author determined the set of minimal non-1-in-degenerate tournaments, constituting a list of six small tournaments.
We begin by bounding the size of a minimal non-$k$-in-degenerate tournament (a similar result appeared in \cite{xing2026excluded}).
\begin{lemma}\label{lemma:degen-size-bound}
    Let $k\geq0$. If $T$ be a minimal non-$k$-in-degenerate tournament, or minimal non-$k$-out-degenerate, then 
    \[
        2k+3\leq |V(T)|\leq \lfloor \frac{1}{8}(2k+5)^2\rfloor.
    \]
\end{lemma}
\begin{proof}
    We may assume after reversing all edges $T$ is minimal non-$k$-in-degenerate. By minimality, every vertex of $T$ has in-degree at least $k+1$. The average in-degree of a tournament is 
    \[
        \frac{\binom{|V(T)|}{2}}{|V(T)|} = \frac{|V(T)|-1}{2},
    \]
    and the average is at least the minimum. So $k+1\leq \frac{|V(T)|-1}{2}$, or equivalently, $2k+3\leq |V(T)|$. \\
    
    Let $S\subseteq V(T)$ be the set of vertices in $T$ of in-degree exactly $k+1$. By the minimality of $T$, every $v\in V(T)$ has an out-neighbour $u\in S$. The total in-degree of $S$ in $T$ is $(k+1)|S|$. The number of arcs from $V(T)\setminus S$ to $S$ is $(k+1)|S|-\binom{|S|}{2}$, and as every vertex has at least one out-neighbour in $S$, we have
        \begin{equation}
            |V(T)|-|S| \leq (k+1)|S|-\binom{|S|}{2}.
        \end{equation}
    Thus 
    \begin{equation}\label{eq:1}
        |V(T)|\leq (k+2)|S|-\binom{|S|}{2}.
    \end{equation}
    The right side of \eqref{eq:1} is maximized at $|S| = k+\frac{5}{2}$, giving the upper bound of 
    \[
    |V(T)|\leq \lfloor \frac{1}{8}(2k+5)^2\rfloor.
    \]
\end{proof}
We now wish to describe the set of minimal non-1-in-degenerate tournaments.
To do so, we need to define some special tournaments.
$T_5$ is the tournament on $v_1, \dots, v_5$ where $v_i\Rightarrow v_j$ if $j-i\equiv 1, 2 \ (\textnormal{mod} \ 5)$, see \Cref{subfig:T5}.
It can easily be verified that $T_5$ is the unique tournament on five vertices where all vertices have in-degree two.
We denote by $U_5$ the tournament obtained from $T_5$ by replacing the edge $v_1v_2$ by the edge $v_2v_1$, see \Cref{subfig:U5}.
We denote by $U_5^{(i)}$ for $i\in[5]$ the tournament obtained from $U_5$ by substituting a tournament of size two for $v_i$, see e.g., \Cref{subfig:1-in-degen-U55} for a backedge graph of $U_5^{(5)}$. 
$P_7^-$ is the unique tournament defined by the backedge graph in \Cref{subfig:1-in-degen-P7^-}.

\begin{theorem}\label{theorem:1-in-degen}
    $T$ is a minimal non-1-in-degenerate tournament if and only if $T$ is isomorphic to one of $T_5$, $\Delta(2,2,2)$, $\Delta(1, C_3, 2)$, $P_7^-$, $Q_6$, $U_5^{(5)}$ (see \Cref{fig:1-in-degen}).
\end{theorem}
\begin{proof}
    One can verify the six tournaments are all minimal non-1-in-degenerate by examining the backedge graphs in \Cref{fig:1-in-degen}. Now, suppose that $T$ is minimal non-1-in-degenerate. Then, every vertex of $T$ has in-degree at least two. By \cref{lemma:degen-size-bound}, $5\leq |V(T)|\leq \lfloor \frac{1}{8}(2+5)^2\rfloor  = 6$. Let $S$ be the vertices of $T$ of in-degree exactly 2. By equation \eqref{eq:1} with $k=1$, we have
    \[
    |V(T)|\leq 3|S|-\binom{|S|}{2}.
    \]
    The right side of the above is positive only if $|S|\leq 5$. If $|S|=5$, then $|V(T)|=5$, and $T = T_5$, as it is the only five-vertex tournament satisfying the minimum in-degree condition. Moreover, if $|S|\leq 2$, then $|V(T)|<5$, a contradiction. The only remaining possibilities are $|S|\in \{3,4\}$. 
    
    %  Then, every vertex of $T$ has in-degree at least two. Let $S\subseteq V(T)$ be the set of vertices in $T$ of in-degree exactly 2. By the minimality of $T$, every $v\in V(T)$ has an out-neighbour $u\in S$. 

    % \begin{claim}\label{claim:fewer-6}
    %     $|V(T)|\leq 6$.
    % \end{claim}
    % \Dinis{I think what we show here (and use later) is that $|S|$ is between $3$ and $5$, $|V(T)|=6$ when $|S|\in[3,4]$, and $|V(T)|=5$ when $|S|=5$, no?}
    % \begin{subproof}
    %     The total in-degree of $S$ in $T$ is $2|S|$. The number of arcs from $V(T)\setminus S$ to $S$ is $2|S|-\binom{|S|}{2}$, and as every vertex has at least one out-neighbour in $S$, we have
    %     \begin{equation}\label{eq:1}
    %         |V(T)|-|S| \leq 2|S|-\binom{|S|}{2}.
    %     \end{equation}
    %     The right side of \cref{eq:1} is negative for $|S|\geq 6$, so $|S|\leq 5$. If $|S|=5$, then $|S| = |V(T)|$. By the minimum degree condition on $T$, we know $|V(T)|\geq 5$, which does not hold if $|S|\leq 2$. Thus, $|S|\geq 3$. If $|S|=3$, \cref{eq:1} implies $|V(T)|\leq 6$, and if $|S|=4$ we similarly have $|V(T)|\leq 6$, as desired. 
    % \end{subproof}

    % \vspace{2em}
    % If $|V(T)| = 5$ then $T = T_5$, as it is the only five-vertex tournament satisfying the minimum in-degree condition. It remains to show that if $T$ has six vertices, then it is one of the tournaments in \Cref{fig:1-in-degen}. As observed previously, the only possibilities are $|S|\in \{3,4\}$. 

    \textbf{Case $|S|=3$}.
    By the handshaking lemma, since $|V(T)|=6$, $T$ has three vertices of in-degree exactly 3, i.e., $T$ has $(3,3,3,2,2,2)$ as its in-degree sequence. Since every vertex in $T$ has an out-neighbour in $S$, $T[S]$ and $T[V(T)\setminus S]$ both form directed triangles. Furthermore, each vertex not in $S$ has an out-neighbour in $S$ (since $T$ is minimal non-1-in-degenerate). From the degree constraints, it follows that each vertex not in $S$ has a  unique and distinct out-neighbour in $S$, in other words, the edges from $V(T)\setminus S$ to $S$ form a perfect matching. After accounting for cyclic symmetry, there are two such perfect matchings, yielding the tournaments $P_7^-$ (\Cref{subfig:1-in-degen-P7^-}) and $\Delta(2,2,2)$ (\Cref{subfig:1-in-degen-D222}). 

    \textbf{Case $|S|=4$}.
    By the handshaking lemma, since $|V(T)|=6$, $T$ has one vertex of in-degree exactly 3, and one vertex of in-degree exactly 4, i.e., $T$ has $(4,3,2,2,2,2)$ as its in-degree sequence. Let $p$ and $q$ be the vertices of in-degree 4 and 3, respectively. Since $p$ has an out-neighbour, say $x$, in $S$, we have $q\Rightarrow p$. It follows that each of $q$ and $p$ has a unique out-neighbour in $S$. Let $y$ denote the out-neighbour of $q$ in $S$. If $x=y$, then $x\Rightarrow S\setminus x$, implying that $T[S\setminus x]$ is a directed triangle. This corresponds to $\Delta(1,C_3,2)$ (\Cref{subfig:1-in-degen-D1C32}). If $x\neq y$, we observe that $T[S]$ has in-degree sequence $(2,2,1,1)$. There is only one tournament with such an in-degree sequence. Now, we distinguish the cases $x\Rightarrow y$ and $y\Rightarrow x$. These yield the tournaments $Q_6$ (\Cref{subfig:1-in-degen-Q6}) and $U_5^{(5)}$ (\Cref{subfig:1-in-degen-U55}), respectively.
\end{proof}
\begin{corollary}\label{corollary:out-degen}
    $T$ is a minimal non-1-out-degenerate tournament if and only if $T$ is isomorphic to one of $T_5$, $\Delta(2,2,2)$, $\Delta(1, 2, C_3),$ $P_7^-, Q_6^*,$ and $U_5^{(3)}$.
\end{corollary}
\begin{proof}
    Since $T$ is is minimal non-1-out-degenerate if and only if $T^*$ is minimal non-1-in-degenerate, it follows from \Cref{theorem:1-in-degen} that $T$ is minimal non-1-out-degenerate if and only if $T$ is isomorphic to one of $T_5^*, \Delta(2,2,2)^*, \Delta(1, C_3, 2)^*, P_7^{-*}, Q_6^*, U_5^{(5)*}$. It is clear by definition that $T_5^* = T_5$, $\Delta(2,2,2)^* = \Delta(2,2,2)$, $\Delta(1, C_3, 2)^* = \Delta(1, 2, C_3)$, from \Cref{fig:1-in-degen} that $P_7^{-*} = P_7^-$, and from \Cref{subfig:U5} that $U_5^{(5)*} = U_5^{(3)}$. 
\end{proof}

An answer to the following question would be a natural follow-up to the results above.
\begin{question}
    Let $k\geq 2$. What are the excluded subtournaments for the classes of $k$-in-degenerate and $k$-out-degenerate tournaments?
\end{question}
Lemma \ref{lemma:degen-size-bound} shows that for each fixed $k$, the set of minimal excluded subtournaments is finite. 

With this, we can state our characterization of minimal nonpaving tournaments.
\begin{restatable}{theorem}{paving}\label{theorem:min-non-paving}
    A tournament $T$ is paving if and only if $T$ does not contain $T_5$ or a tournament with a backedge graph in \Cref{fig:min-non-paving} as a subtournament.
\end{restatable}
\begin{figure}[h]
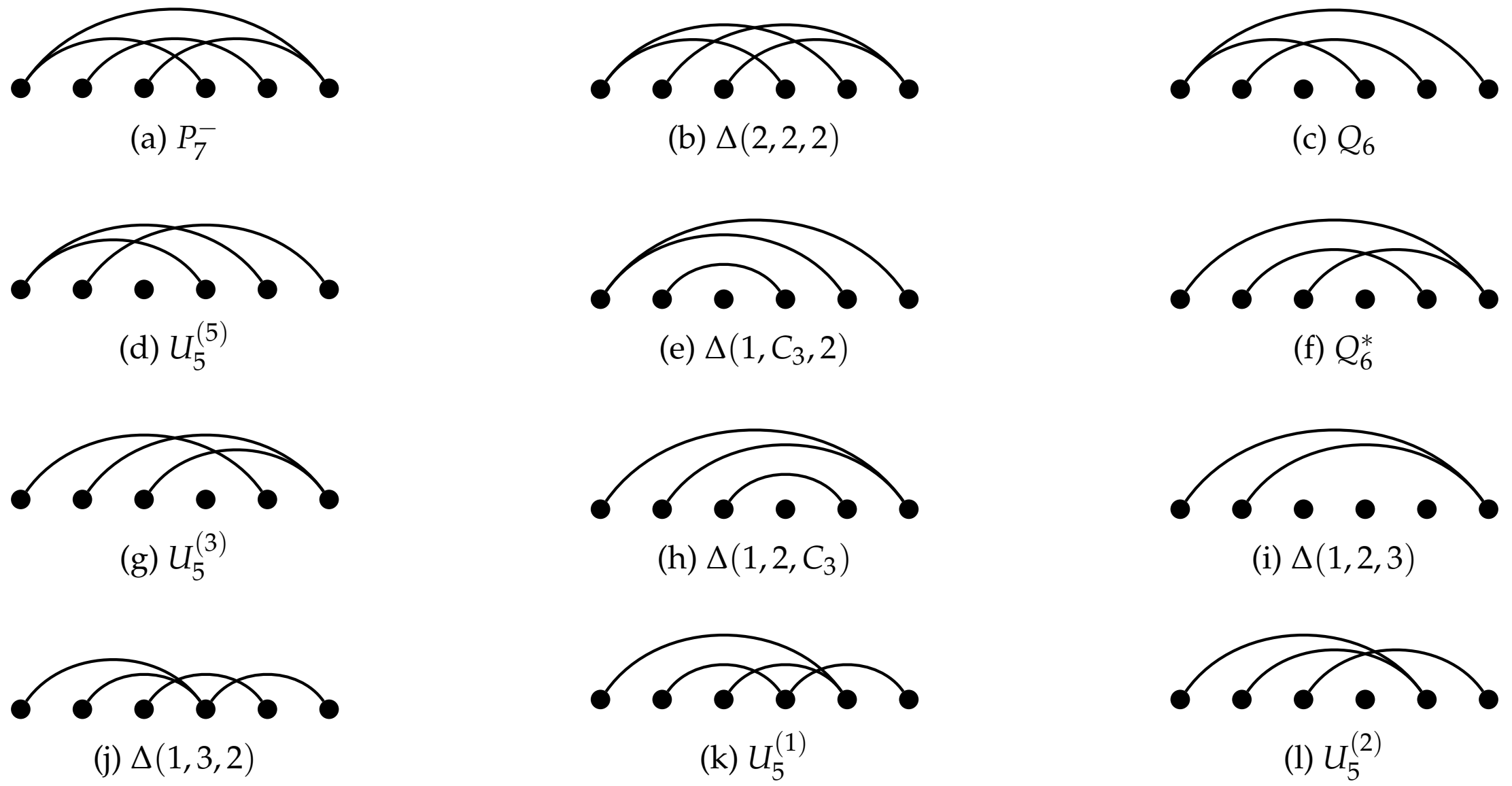

    \centering
    %1,2,4
    \begin{subfigure}{0.32\linewidth}\centering
    \BackedgeGraphUnlabeled{1/4,1/6,2/5,3/6}{6}\caption{$P_7^-$}
    \end{subfigure}\hfill
    \begin{subfigure}{0.32\linewidth}\centering
    \BackedgeGraphUnlabeled{1/4,1/5,2/6,3/6}{6}\caption{$\Delta(2,2,2)$}
    \end{subfigure}\hfill
    \begin{subfigure}{0.32\linewidth}\centering
    \BackedgeGraphUnlabeled{1/4,1/6,2/5}{6}\caption{$Q_6$}
    \end{subfigure}
    
    \vspace{3mm}

    %5,6,7
    \begin{subfigure}{0.32\linewidth}\centering
    \BackedgeGraphUnlabeled{1/4,1/5,2/6}{6}\caption{$U_5^{(5)}$}
    \end{subfigure}\hfill
    \begin{subfigure}{0.32\linewidth}\centering
    \BackedgeGraphUnlabeled{1/5,1/6,2/4}{6}\caption{$\Delta(1,C_3,2)$}\label{subfig:min-non-paving-Q6*}
    \end{subfigure}\hfill
    \begin{subfigure}{0.32\linewidth}\centering
    \BackedgeGraphUnlabeled{1/6,2/5,3/6}{6}\caption{$Q_6^*$}
    \end{subfigure}
    
    \vspace{3mm}

    %9,10,11
    \begin{subfigure}{0.32\linewidth}\centering
    \BackedgeGraphUnlabeled{1/5,2/6,3/6}{6}\caption{$U_5^{(3)}$}\label{subfig:min-non-paving=U53}
    \end{subfigure}\hfill
    \begin{subfigure}{0.32\linewidth}\centering
    \BackedgeGraphUnlabeled{1/6,2/6,3/5}{6}\caption{$\Delta(1,2,C_3)$}\label{subfig:min-non-paving-D(12C3)}
    \end{subfigure}\hfill
    \begin{subfigure}{0.32\linewidth}\centering
    \BackedgeGraphUnlabeled{1/6,2/6}{6}\caption{$\Delta(1,2,3)$}
    \end{subfigure}
    
    \vspace{3mm}
    
    %12,13,14
    \begin{subfigure}{0.32\linewidth}\centering
    \BackedgeGraphUnlabeled{1/4,2/4,3/5,4/6}{6}\caption{$\Delta(1,3,2)$}
    \end{subfigure}\hfill
    \begin{subfigure}{0.32\linewidth}\centering
    \BackedgeGraphUnlabeled{1/5,2/4,3/5,4/6}{6}\caption{$U_5^{(1)}$}
    \end{subfigure}\hfill
    \begin{subfigure}{0.32\linewidth}\centering
    \BackedgeGraphUnlabeled{1/5,2/5,3/6}{6}\caption{$U_5^{(2)}$}
    \end{subfigure}
    
    \caption{The minimal non-paving tournaments on six vertices.}
    \label{fig:min-non-paving}
\end{figure}
The tournaments in \Cref{fig:min-non-paving} comprise the five minimal non-1-in-degenerate tournaments on six vertices we obtain from \Cref{theorem:1-in-degen}, the additional three minimal non-1-out-degenerate tournaments from \Cref{corollary:out-degen}, and four more tournaments. The additional four are $\Delta(1,2,3), \Delta(1,3,2)$, $U_5^{(1)}$, and $U_5^{(2)}$.

\begin{figure}[ht]
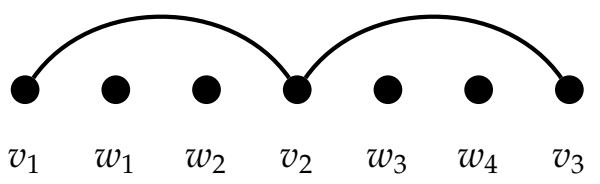
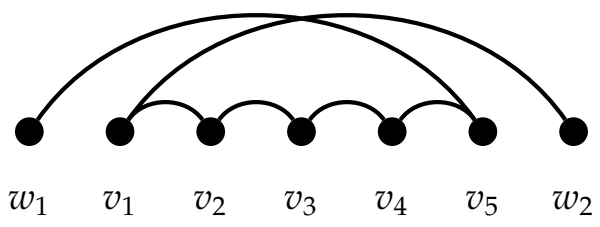
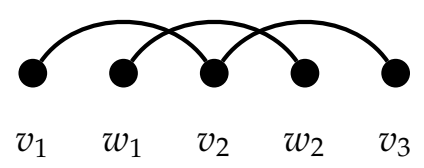
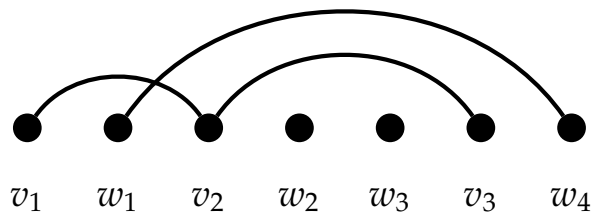
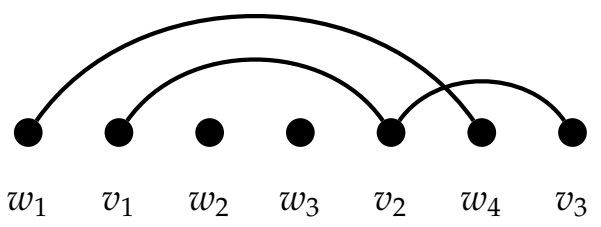

    \centering
    
    \begin{subfigure}{0.48\linewidth}\centering
    \BackedgeGraphGeneral{1/4,4/7}{7}{$v_1$, $w_1$, $w_2$, $v_2$, $w_3$, $w_4$, $v_3$}{} % type 2, k=1
    \caption{$M_{2,2}^{(3)}$}
    \end{subfigure}\hfill
    \begin{subfigure}{0.48\linewidth}\centering
    \BackedgeGraphGeneral{1/6,2/7,2/3,3/4,4/5,5/6}{7}{$w_1$, $v_1$, $v_2$, $v_3$, $v_4$, $v_5$, $w_2$}{} % type 3, k=2
    \caption{$M_{1,1,is}^{(5)}$}
    \end{subfigure}

    \vspace{4mm}

    \begin{subfigure}{0.48\linewidth}\centering
    \BackedgeGraphGeneral{1/3,2/4,3/5}{5}{$v_1$, $w_1$, $v_2$, $w_2$, $v_3$}{} % type 5, k=1
    \caption{$M_{1,1,adj}^{(3)}$}
    \end{subfigure}\hfill
    \begin{subfigure}{0.48\linewidth}\centering
    \BackedgeGraphGeneral{1/3,2/7,3/6}{7}{$v_1$, $w_1$, $v_2$, $w_2$, $w_3$, $v_3$, $w_4$}{} % type 1, k=1
    \caption{$M_{1,2}^{(3)}$}
    \end{subfigure}\hfill
    
    \vspace{4mm}
    
    \begin{subfigure}{0.48\linewidth}\centering
    \BackedgeGraphGeneral{1/6,2/5,5/7}{7}{$w_1$, $v_1$, $w_2$, $w_3$, $v_2$, $w_4$, $v_3$}{} % type 4, k=1
    \caption{$M_{2,1}^{(3)}$}
    \end{subfigure}
    \caption{Backedge graphs of the smallest tournaments in each of the five infinite classes of minimal nonmatching tournaments}\label{fig:base-tournaments}
\end{figure}

To state our second main result (\Cref{theorem:main}), we require the following classes of tournaments, see Figure \ref{fig:base-tournaments}. 
\begin{definition}
    \label{def: infinite classes minimal nonmatching}
    Let $k$ be a positive integer.
    \begin{itemize}
        \item[a)] Let $B$ be the backedge graph obtained from an induced monotone path on $2k+1$ vertices $v_1,\dots,v_{2k+1}$ by adding two isolated vertices $w_1,w_2$ between $v_1$ and $v_2$ of the path, and two isolated vertices $w_3,w_4$ between $v_{2k}$ and $v_{2k+1}$.
            Formally, $V=\{v_1,\dots,v_{2k+1},w_1,\dots,w_4\}$, $E=\{v_iv_{i+1}\mid 1\leq l<2k+1\}$, and the ordering inducing $B$ is $\sigma=v_1,w_1,w_2,v_2,\dots,v_{2k},w_3,w_4,v_{2k+1}$.
            We denote the tournament defined by this backedge graph and ordering by $M_{2,2}^{(2k+1)}$.
        \item[b)] Let $B$ be the backedge graph obtained from an induced monotone path on $2k+3$ vertices $v_1,\dots,v_{2k+3}$ by adding a vertex $w_1$ before $v_1$, a vertex $w_2$ after $v_{2k+3}$, and the edges $w_1v_{2k+3}$ and $v_1w_2$.
            Formally, $V=\{v_1,\dots,v_{2k+3},w_1,w_2\}$, $E=\{v_iv_{i+1}\mid 1\leq l<2k+3\}\cup \{w_1v_{2k+3},v_1w_2\}$, and the ordering inducing $B$ is $\sigma=w_1,v_1,\dots,v_{2k+3},w_2$.
            We denote the tournament defined by this backedge graph and ordering by $M_{1,1,is}^{(2k+3)}$.

        \item[c)] Let $B$ be the backedge graph obtained from an induced monotone path on $2k+1$ vertices $v_1,\dots,v_{2k+1}$ by adding one vertex $w_1$ between the $v_1$ and $v_2$, a vertex $w_2$ between $v_{2k}$ and $v_{2k+1}$, and the edge $w_1w_2$.
            Formally, $V=\{v_1,\dots,v_{2k+1},w_1,w_2\}$, $E=\{v_iv_{i+1}\mid 1\leq l<2k+1\}\cup \{w_1w_2\}$, and the ordering inducing $B$ is $\sigma=v_1,w_1,v_2,\dots,v_{2k},w_2,v_{2k+1}$.
            We denote the tournament defined by this backedge graph and ordering by $M_{1,1,adj}^{(2k+1)}$. (This is the infinite family of minimal non-matching tournaments described in \cite{liuVariousTheoremsTournaments2012}.)

        % \Dinis{The notation $M_{1,1,adj}^{(2k+1)}$ is pretty terrible, but I couldn't immediately come up with any substantial improvements.
        % Better alternatives are, naturally, welcome.}

        % \Dinis{I think this is the infinite family $Q_n$ of minimally nonmatching tournaments in \cite{liuVariousTheoremsTournaments2012} in case anyone was wondering/thinks that this should be mentioned somewhere.}

        \item[d)] Let $B$ be the backedge graph obtained from an induced monotone path on $2k+1$ vertices $v_1,\dots,v_{2k+1}$ by adding one vertex $w_1$ between $v_1$ and $v_2$, two vertices $w_2,w_3$ between $v_{2k}$ and $v_{2k+1}$, one vertex $w_4$ after $v_{2k+1}$, and the edge $w_1w_4$.
            Formally, $V=\{v_1,\dots,v_{2k+1},w_1,\dots,w_4\}$, $E=\{v_iv_{i+1}\mid 1\leq l<2k+1\}\cup \{w_1w_4\}$, and the ordering inducing $B$ is $\sigma=v_1,w_1,v_2,\dots,v_{2k},w_2,w_3,v_{2k+1},w_4$.
            We denote the tournament defined by this backedge graph and ordering by $M_{1,2}^{(2k+1)}$.

        \item[e)] Let $B$ be the backedge graph obtained from an induced monotone path on $2k+1$ vertices $v_1,\dots,v_{2k+1}$ by adding one vertex $w_1$ before $v_1$, two vertices $w_2,w_3$ between $v_1$ and $v_2$, one vertex $w_4$ between $v_{2k+1}$ and $v_{2k+1}$, together with the edge $w_1w_4$.
            More precisely, $B$ is the mirror of the backedge graph used to define $M_{1,2}^{(2k+1)}$.
            We denote the tournament defined by this backedge graph and ordering by $M_{2,1}^{(2k+1)}$.
    \end{itemize}

    We respectively refer to each of these five families as types a) to e).
    When the parameter $k$ is fixed, we write $T_a,T_b,T_c,T_d,T_e$ for representatives of each of these five types.
\end{definition}
We additionally require the tournament $M_7$ defined by the backedge graph in \Cref{fig:M7}.
\begin{figure}[ht]
\centering
\BackedgeGraphGeneral{1/5, 2/6, 3/7, 3/4, 4/5}{7}{$w_1$, $w_2$, $v_1$, $v_2$, $v_3$, $w_3$, $w_4$}{}
\caption{$M_7$}\label{fig:M7}
\end{figure}
With this, we can finally describe the set of minimal nonmatching tournaments.
\begin{restatable}{theorem}{matching}
    \label{theorem:main}
    $T$ is a minimal nonmatching tournament if and only if $T$ is $T_5$, $M_7$, $P_7^-$, or one of $M_{2,2}^{(2k+1)}$, $M_{1,1,is}^{(2k+3)}$, $M_{1,1,adj}^{(2k+1)}$, $M_{1,2}^{(2k+1)}$, $M_{2,1}^{(2k+1)}$ where $k$ is a positive integer.
\end{restatable}

After studying (non)matching tournaments, we thought that it might be just as natural (if not more so) to consider tournaments $T$ with backedge graphs $B_\sigma$ where the edges form a matching, and there are no two edges of the form $ik, jl$, where $i\prec_\sigma j\prec_\sigma k \prec_\sigma l$.
We call such tournaments noncrossing matching (NCM for short) tournaments.
Our final contribution is the following characterization of the set of NCM tournaments.
\begin{restatable}{theorem}{ncm}
    \label{thm:NCM}
    A tournament $T$ is NCM if and only if $T$ contains no $\Delta(1,2,2), T_5, U_5, W_5$. 
\end{restatable}

\section{Paving tournaments}\label{section:paving}
In this section we prove a characterization of paving tournaments by first reducing to minimal non-$1$-in-degenerate and non-$1$-out-degenerate tournaments.
Recall that we aim to prove the following theorem.
\paving

\begin{figure}[ht]
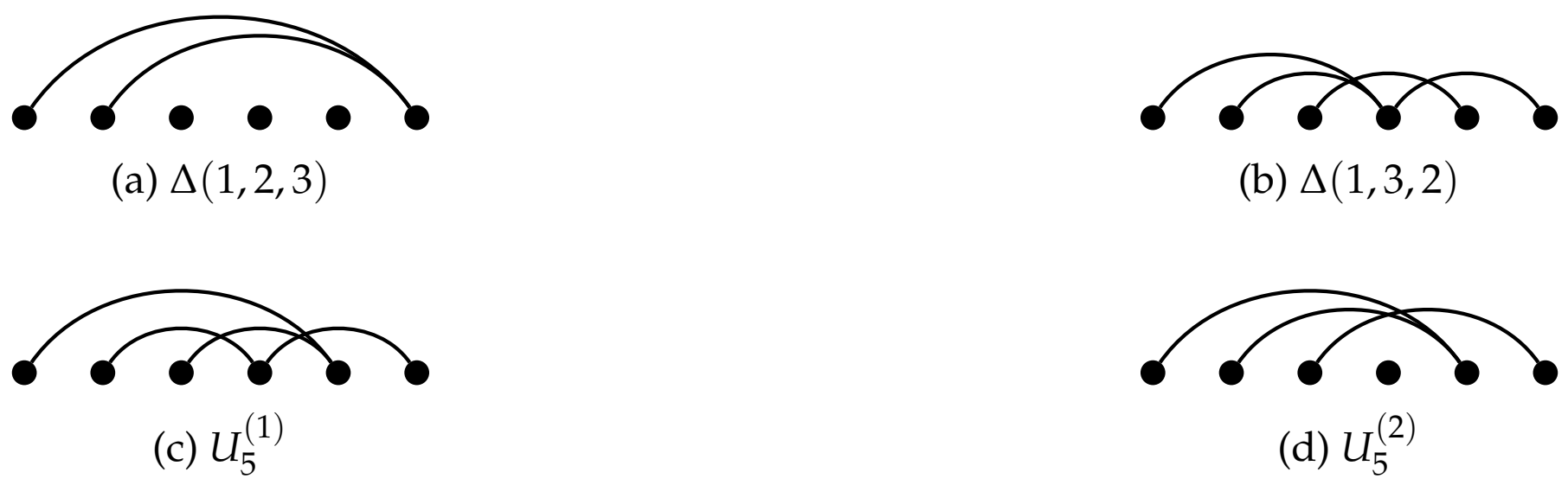

    \centering
    \begin{subfigure}{0.48\linewidth}\centering
    \BackedgeGraphUnlabeled{1/6,2/6}{6}\caption{$\Delta(1,2,3)$}\label{subfig:min-non-paving-1-in-D(123)}
    \end{subfigure}\hfill
    \begin{subfigure}{0.48\linewidth}\centering
    % updated order: (1,3,5,2,4,6)
    \BackedgeGraphUnlabeled{1/4,2/4,3/5,4/6}{6}\caption{$\Delta(1,3,2)$}\label{subfig:min-non-paving-1-in-D(132)}
    \end{subfigure}
    
    \vspace{4mm}
    
    \begin{subfigure}{0.48\linewidth}\centering
    % updated order: (1,2,4,6,3,5)
    \BackedgeGraphUnlabeled{1/5,2/4,3/5,4/6}{6}\caption{$U_5^{(1)}$}\label{subfig:min-non-paving-1-in-U51}
    \end{subfigure}\hfill
    \begin{subfigure}{0.48\linewidth}\centering
    \BackedgeGraphUnlabeled{1/5,2/5,3/6}{6}\caption{$U_5^{(2)}$}\label{subfig:min-non-paving-1-in-U52}
    \end{subfigure}
    
    \caption{Backedge graphs of the 1-in-degenerate minimal non-paving tournaments on six vertices.}
    \label{fig:min-non-paving-1-in}
\end{figure}

There are four tournaments on four vertices up to isomorphism. It is easy to check that they all admit a paving ordering.
\begin{observation}\label{observation:4-vtx-paving}
    Every four-vertex tournament is paving.
\end{observation}
Moreover, it is known that there is precisely one five-vertex tournament that is nonpaving.
\begin{lemma}[\cite{kimStructureD122025}]\label{lemma:C5-unique-non-paving}
    $T_5$ is the only nonpaving tournament on five vertices.
\end{lemma}
% \begin{proof}
%     Let $T$ be a tournament on five vertices distinct from $T_5$. Then $T$ has a vertex $v$ with $d^-(v)\leq 1$, so $N^+(v)$ contains at least three vertices, say $u_1, u_2, u_3$. 
% \end{proof}
In addition, the sets of minimal non-1-in- and minimal non-1-out-degenerate tournaments are related to that of minimal nonpaving tournaments.
\begin{lemma}\label{lemma:non-1-degen-implies-non-paving}
    If $T$ is a minimal non-1-in-degenerate or minimal non-1-out-degenerate, then $T$ is minimal nonpaving. 
\end{lemma}
\begin{proof}
    By symmetry in the definition of paving with respect to in- and out-neighbours, it suffices to verify that if $T$ is a minimal non-1-in-degenerate tournament, then $T$ is minimal nonpaving. As $T_5$ has minimum in-degree 2, $T_5$ is non-paving, and thus minimal-non-paving by \Cref{observation:4-vtx-paving}. The remaining tournaments have six vertices, and do not contain $T_5$ as a subtournament by \Cref{theorem:1-in-degen}. By \Cref{lemma:C5-unique-non-paving}, all their five-vertex subtournaments are paving. Therefore, all minimal non-1-in-degenerate tournaments are also minimal nonpaving.
\end{proof}
Further, all of the following tournaments are minimal nonpaving.
\begin{lemma}\label{lemma:paving-sufficient}
    The tournaments $\Delta(1,2,3), \Delta(1,3,2), U_5^{(1)}, U_5^{(2)}$ are minimal nonpaving.
\end{lemma}
\begin{proof}
    Since $U_5^* \cong U_5$ as exhibited by the permutation $\pi = (v_1\ v_2)(v_3\ v_5)$, we have $U_5^{(1)*}\cong U_5^{(2)}$. Additionally, $\Delta(1,2,3)^*\cong \Delta(1,3,2)$. Since the paving property is preserved by mirroring, it thus suffices to prove that $\Delta(1,2,3)$ and $U_5^{(1)}$ are minimal nonpaving. We first show these tournaments are nonpaving. 

    Let $T\cong \Delta(1,2,3)$ have vertex set $a, b_1, b_2, c_1,c_2,c_3$ where $a\Rightarrow \{b_1,b_2\}\Rightarrow \{c_1,c_2,c_3\}\Rightarrow a$. Suppose $T$ has a paving ordering $\sigma$. Since $b_1$ is the only vertex with in-degree at most one, and $c_3$ is the only vertex with out-degree at most one, we have that $b_1$ is first and $c_3$ is last in $\sigma$. The in-neighbour of $b_1$ is $a$, and $b_1\prec_\sigma a$, so we have $a\prec_\sigma b_2$, since $a$ has at most one left neighbour in $B_\sigma$. But since $\{c_1, c_2,c_3\}\Rightarrow a$, and $a\prec_\sigma c_3$, we also have $c_1, c_2\prec_\sigma a$. Thus
    \[
    c_1, c_2\prec_\sigma a\prec_\sigma b_2,
    \]
    so $b_2$ has two left-neighbours in $B_\sigma$, a contradiction.

    Let $T'\cong U_5^{(1)}$. Recall that $U_5$ is obtained from $T_5$ by reversing the arc $v_1v_2$; $U_5^{(1)}$ is obtained from $U_5$ by substituting a tournament on two vertices for $v_1$. Say the two vertices are $u_1, w_1$, where $u_1\Rightarrow w_1$ in $T'$. Observe that $v_2$ is the unique vertex in $T'$ with in-degree one, and $w_1$ is the unique vertex with out-degree one. Suppose $T'$ has a paving ordering $\sigma'$. Then $v_2$ is first and $w_1$ is last in $\sigma'$. Note that $v_3$ is the only out-neighbour of $w_1$. Since $v_3$ has at most one right-neighbour in $B_{\sigma'}$, $u_1\prec_{\sigma'}v_3\prec_{\sigma'}w_1$. Since $v_5\Rightarrow \{u_1,w_1\}$ and $v_5\Rightarrow v_2$, and $v_2$ is first in $\sigma'$, we have $v_5\prec_{\sigma'}u_1$. Thus
    \[
        v_2\prec_{\sigma'}v_5\prec_{\sigma'}u_1\prec_{\sigma'}v_3\prec_{\sigma'}w_1.
    \]
    Hence, in $B_{\sigma'}$ has a left-neighbour $v_2$ and a right-neighbour $v_3$. It follows that $v_4\prec_{\sigma'}v_5$, but then $v_3$ has two left-neighbours, a contradiction. 
    
    It remains to check minimality, which equates to verifying $\Delta(1,2,3)$ and $U_5^{(1)}$ do not contain $T_5$, by \Cref{observation:4-vtx-paving} and \Cref{lemma:C5-unique-non-paving}. This follows from the fact that each of $\Delta(1,2,3)$ and $U_5^{(1)}$ is $1$-in-degenerate (see Figure \ref{fig:min-non-paving-1-in}), but $T_5$ is not (by Theorem \ref{theorem:1-in-degen}). 
\end{proof}

Finally, the tournaments above are the only 1-in-degenerate and 1-out-degenerate minimal non-paving tournaments.
\begin{lemma}\label{lemma:paving-necessary}
    Let $T$ be a 1-in-degenerate and 1-out-degenerate minimal non-paving tournament. Then $T$ is one of $\Delta(1,2,3)$, $\Delta(1,3,2)$, $U_5^{(1)}$, and $U_5^{(2)}$ (see \Cref{fig:min-non-paving-1-in}).
\end{lemma}
\begin{proof}
    Let $T$ be a 1-in-degenerate and 1-out-degenerate, minimal non-paving tournament, and let $v$ be a vertex in $T$ with $d^-(v)\leq 1$. If $d^-(v)=0$, we obtain a paving ordering of $T$ from one of $T\setminus v$. Hence, $d^-(v)=1$. Let $\sigma = v_1,\dots, v_n$ be a paving ordering of $T\setminus v$, and $B_\sigma = B_\sigma(T\setminus v)$. After potentially swapping the order of consecutive vertices in $\sigma$, we may assume there are no edges between consecutive vertices in $B_\sigma$.
  
    Since $T$ is nonpaving, $v, v_1, \dots, v_n$ is not a paving ordering. Let $v_i$ be the in-neighbour of $v$. Then $v_i$ has exactly one left-neighbour $v_j$, and $j\leq i-2$. 
    We consider cases depending on the position of $j$ relative to $i$.
    \begin{itemize}[align=left,leftmargin=*]
        \item [Case 1:] $j\leq i-4$.
            Then $v_i\Rightarrow \{v, v_j\}\Rightarrow \{v_{i-3}, v_{i-2}, v_{i-1}\}\Rightarrow v_i$, so $T$ either contains $\Delta(1,2,3)$ (\Cref{subfig:min-non-paving-1-in-D(123)}), or $\Delta(1,2,C_3)$ (\Cref{subfig:min-non-paving-D(12C3)}), depending on whether $T[\{v_{i-3}, v_{i-2}, v_{i-1}\}]$ is transitive or a cyclic triangle. The former being one of our desired outcomes, the latter contradicting 1-out-degeneracy, by \Cref{corollary:out-degen}. 
            % \BackedgeGraphUnlabeled{1/6,2/6,3/5}{6}
            % \BackedgeGraphUnlabeled{1/6,2/6}{6}
            % depending on whether $T[\{v_{i-3}, v_{i-2}, v_{i-1}\}]$ is transitive or a cyclic triangle.
            
        \item [Case 2:] $j=i-3$.
            If $v_j$ has a left-neighbour or $v_i$ has a right-neighbour, we may rearrange the order to obtain a situation handled by Case 1.
            That is, $v, v_{i-3}, v_{i-2}, v_{i-1}, v_i$ and the left-neighbour of $v_j$ or the right neighbour of $v_i$ together induce $\Delta(1, 2, 3)$ (\Cref{subfig:min-non-paving-1-in-D(123)}) or $\Delta(1, 2, C_3)$ (\Cref{subfig:min-non-paving-D(12C3)}). We may thus assume $v_j$ has no left neighbour and $v_i$ has no right-neighbour. If $v_{i-2}$ has a right-neighbour, said neighbour is to the right of $v_i$, since $v_i$ is not adjacent to $v_{i-2}$ and $v_{i-2}$ is not adjacent to $v_{i-1}$. Thus, we obtain $U_5^{(2)}$ (\Cref{subfig:min-non-paving-1-in-U52}), one of the four desired outcomes. So we may assume that $v_{i-2}$ has no right-neighbour. If $v_{i-1}$ has a neighbour to the left, said neighbour is to the left of $v_{j}$, so we obtain $Q_6^*$ (\Cref{subfig:min-non-paving-Q6*}), which is not 1-out-degenerate (\Cref{corollary:out-degen}). Thus, $v_{i-1}$ has no left-neighbour. We may thus reorder $v_{i-3},\dots, v_i$ into $v_{i-2},v_iv_{i-3},v_{i-1}$ to obtain a paving ordering of $T$. This reordering adds the edges $v_iv_{i-1}, v_{i-2}v_{i-3}$ and removes the edge $v_{i-3}v_i$, and, by the previous observations, this is a paving ordering of $T$, as shown in Figure \ref{fig:paving3}, which is impossible. 
            \begin{figure}[ht]
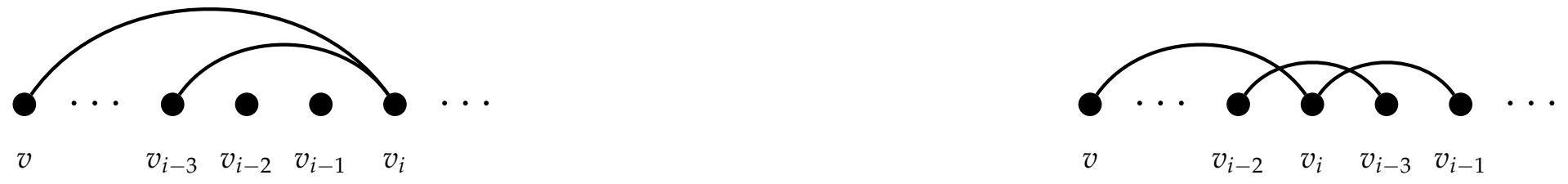

                \centering
                \begin{subfigure}{0.48\linewidth}\centering
                \BackedgeGraphGeneral{1/5, 2/5}{5}{$v$, $v_{i-3}$, $v_{i-2}$, $v_{i-1}$, $v_i$}{1, 5}
                \end{subfigure}\hfill
                \begin{subfigure}{0.48\linewidth}\centering
                \BackedgeGraphGeneral{1/3, 2/4, 3/5}{5}{$v$, $v_{i-2}$, $v_{i}$, $v_{i-3}$, $v_{i-1}$}{1, 5}
                \end{subfigure}
                \caption{Reordering to obtain a paving ordering in Case 2.}
                \label{fig:paving3}
            \end{figure}
            
        \item[Case 3:] $j=i-2$. We may assume that either $v_{i-2}$ has no left-neighbour, or $v_i$ has no right-neighbour, as otherwise $T$ again contains $\Delta(1,2,3)$ (\Cref{subfig:min-non-paving-1-in-D(123)}), noting that said left- and right-neighbours together with $v_{i-1}$ form an independent set of size three in $B_\sigma$.
    
            If $v_{i-2}$ has no left-neighbour and $v_i$ has a right-neighbour $v_l$, we may reorder $v_{i-2}, v_{i-1}, v_i$ into $v_{i-1}, v_{i}, v_{i-2}$. This produces the edge $v_{i-2}v_{i-1}$ and results in a paving ordering unless $v_{i-1}$ already has a right-neighbour $v_k$ in $B_\sigma$. If $v_{i-1}$ has a left-neighbour, we obtain $U_5^{(1)}$ (\Cref{subfig:min-non-paving-1-in-U51}) using the original ordering $v_{i-2},v_{i-1},v_{i}$ together with $v$ and the two neighbours of $v_{i-1}$. 
            Thus, we may assume that $v_{i-1}$ has no left-neighbour in $B_\sigma$, but has a right-neighbour $v_k$. If $k>l$ we obtain $U_5^{(2)}$ (\Cref{subfig:min-non-paving-1-in-U52}), by placing $v_l$ directly to the left of $v_i$ in $\sigma$ as shown in Figure \ref{fig:paving2}, one of our desired outcomes.
            \begin{figure}[ht]
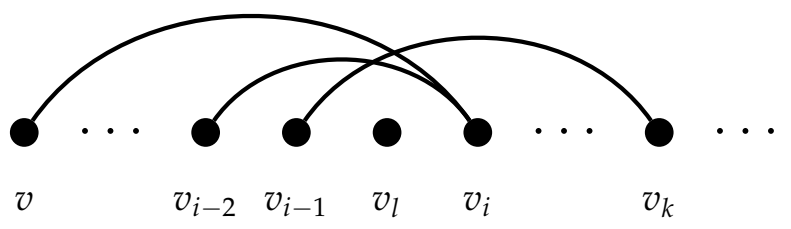

                \centering
                \BackedgeGraphGeneral{1/5,2/5,3/6}{6}{$v$, $v_{i-2}$, $v_{i-1}$, $v_{l}$, $v_i$, $v_k$}{1, 5, 6}
                \caption{A copy of $U_5^{(2)}$ obtained when $v_i$ has a right-neighbour $v_l$; and $v_{i-1}$ has a right-neighbour $v_k$ with $k > l$.}
                \label{fig:paving2}
            \end{figure}
            
            Now, $k\neq l$, as otherwise $v_k$ has two left-neighbours, $v_i$ and $v_{i-1}$, contradicting the fact that  $B_{\sigma}\setminus v$ is paving.
            If $k<l$, we obtain $\Delta(1,3,2)$ (\Cref{subfig:min-non-paving-1-in-D(132)}), one of the four desired outcomes. 
        
            In the case where $v_i$ has no right-neighbour, we may reorder $v_{i-2}, v_{i-1}, v_i$ into $v_i, v_{i-2}, v_{i-1}$. This produces the edge $v_iv_{i-1}$ and results in a paving ordering unless $v_{i-1}$ already had a left-neighbour $v_k$ in $B_\sigma$. In addition, reordering $v_{i-2}, v_{i-1}, v_i$ into $v_{i-1}, v_{i}, v_{i-2}$ results in a paving ordering unless $v_{i-1}$ already had a right neighbour $v_{l}$ in $B_\sigma$ (which gives $U_5^{(1)}$, as seen in Figure \ref{fig:paving1}), or $v_{i-2}$ already had a left-neighbour $v_m$ in $B_{\sigma}$. 
        
            \begin{figure}[ht]
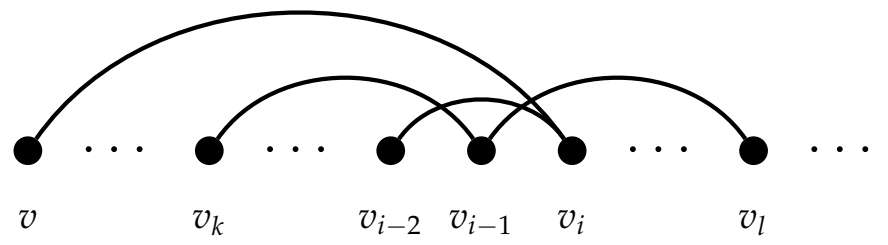

                \centering
                \BackedgeGraphGeneral{1/5,2/4,3/5,4/6}{6}{$v$, $v_k$,$v_{i-2}$, $v_{i-1}$, $v_i$, $v_l$}{1,2, 5, 6}
                \caption{A copy of $U_5^{(1)}$ obtained when $v_{i-1}$ has a left-neighbour $v_k$ and a right-neighbour $v_l$.}
                \label{fig:paving1}
            \end{figure}
        
            If $k < m$, then we obtain $Q_6^*$ (by moving $v_m$ between $v_{i-2}$ and $v_{i-1}$; see \Cref{subfig:min-non-paving-Q6*}). If $k > m$, and we obtain $U_5^{(3)}$ (\Cref{subfig:min-non-paving=U53}) as shown in Figure \ref{fig:pavingfinal}. 
            Once more, $k=m$ is not an option, since the ordering of $B_{\sigma}-v$ is a paving backedge graph.
            \begin{figure}[ht]
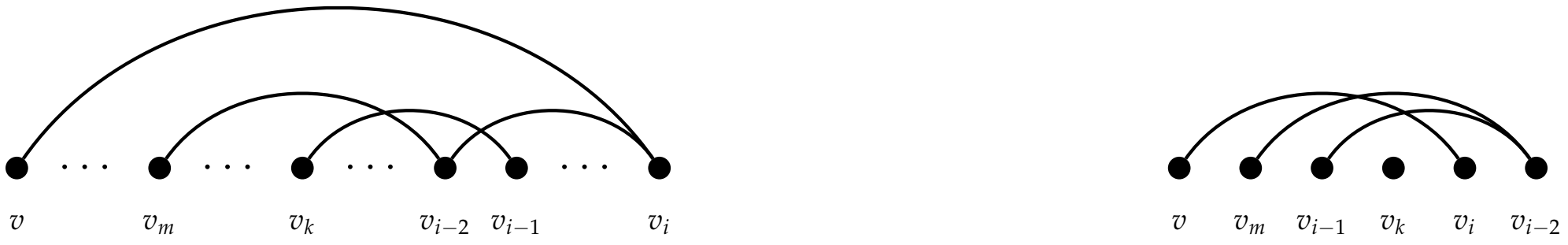

                \centering
        
                \begin{subfigure}{0.48\linewidth}\centering
                \BackedgeGraphGeneral{1/6,3/5,4/6,2/4}{6}{$v$, $v_m$, $v_k$,$v_{i-2}$, $v_{i-1}$, $v_i$}{1,2, 5, 3}
                \end{subfigure}\hfill
                \begin{subfigure}{0.48\linewidth}\centering
                \BackedgeGraphGeneral{1/5, 2/6, 3/6}{6}{$v$, $v_m$, $v_{i-1}$,$v_k$, $v_{i}$, $v_{i-2}$}{}
                \end{subfigure}     
                \caption{Obtaining $U_5^{(3)}$ in Lemma \ref{lemma:paving-necessary}.}
                \label{fig:pavingfinal}
            \end{figure}
    \end{itemize}
    This concludes the proof.
\end{proof}

With this, proving \Cref{theorem:min-non-paving} is straightforward.
\begin{proof}[Proof of \Cref{theorem:min-non-paving}]
    By \Cref{lemma:non-1-degen-implies-non-paving}, \Cref{theorem:1-in-degen}, and \Cref{lemma:paving-sufficient}, every tournament with a backedge graph in \Cref{fig:min-non-paving} is minimal non-paving. If $T$ is non-1-in-degenerate or non-1-out-degenerate, then $T$ contains a tournament with a backedge graph in \Cref{fig:min-non-paving} by \Cref{theorem:1-in-degen} and \Cref{corollary:out-degen}. If $T$ is non-paving, 1-in-degenerate and 1-out-degenerate, then $T$ contains a tournament with a backedge graph in \Cref{fig:min-non-paving} by \Cref{lemma:paving-necessary}. This concludes the proof. 
\end{proof}

\section{Matching tournaments}\label{section:main}
In this section, we focus on matching tournaments and prove that the aforementioned characterization of minimal nonmatching tournaments is correct.
For convenience, the theorem that we aim to prove is restated below.
\matching
The proof of \Cref{theorem:main} has four parts.  First we show that the listed tournaments are nonmatching.  Second, we show that their proper subtournaments are matching.  Third, we record that the listed obstructions are pairwise incomparable.  Finally, we prove unavoidability: every minimal nonmatching tournament appears in the list.
The last step, which is by far the hardest, is somewhat similar to the proof of unavoidability of the determined obstacles for being paving.
As there, we distinguish cases depending on whether the tournament under consideration is $1$-in-degenerate.
When it is not, we use \Cref{theorem:1-in-degen} to conclude the proof.
When it is, we latch on to a vertex $s$ of in-degree one and conduct a (somewhat cumbersome) case-analysis dependent on the subgraph induced by the second-neighbourhood of $s$ in a backedge graph that would become matching upon the deletion of $s$.

\begin{figure}
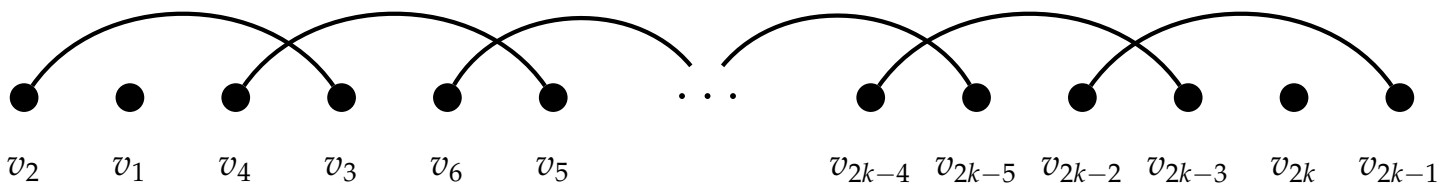

    \centering
    \begin{subfigure}{0.48\linewidth}\centering
    \BackedgeGraphGeneral{}{1}{$v_1$}{}
    \caption{$R_1$}\label{subfig:R1}
    \end{subfigure}\hfill

    \vspace{5mm}
    
    \begin{subfigure}{0.48\linewidth}\centering
    \BackedgeGraphGeneral{1/2}{2}{$v_1$, $v_2$}{}
    \caption{$B_{\sigma_1}$ of $R_2$}\label{subfig:R2 1}
    \end{subfigure}\hfill
    \begin{subfigure}{0.48\linewidth}\centering
    \BackedgeGraphGeneral{}{2}{$v_2$, $v_1$}{}
    \caption{$B_{\sigma_2}$ of $R_2$}\label{subfig:R2 2}
    \end{subfigure}\hfill

    \vspace{5mm}

    \begin{subfigure}{0.24\linewidth}\centering
    \BackedgeGraphGeneral{1/2, 2/3}{3}{$v_1$, $v_2$, $v_3$}{}
    \caption{$B_{\sigma}$ of $R_3$}\label{subfig:R3 nonmatching}
    \end{subfigure}\hfill
    \begin{subfigure}{0.24\linewidth}\centering
    \BackedgeGraphGeneral{1/3}{3}{$v_1$, $v_3$, $v_2$}{}
    \caption{$B_{\sigma_1}$ of $R_3$}\label{subfig:R3 matching 1}
    \end{subfigure}\hfill
    \begin{subfigure}{0.24\linewidth}\centering
    \BackedgeGraphGeneral{1/3}{3}{$v_2$, $v_1$, $v_3$}{}
    \caption{$B_{\sigma_2}$ of $R_3$}\label{subfig:R3 matching 2}
    \end{subfigure}\hfill
    \begin{subfigure}{0.24\linewidth}\centering
    \BackedgeGraphGeneral{1/3}{3}{$v_3$, $v_2$, $v_1$}{}
    \caption{$B_{\sigma_3}$ of $R_3$}\label{subfig:R3 matching 3}
    \end{subfigure}\hfill

    \vspace{5mm}

    \begin{subfigure}{0.24\linewidth}\centering
    \BackedgeGraphGeneral{1/2, 2/3, 3/4}{4}{$v_1$, $v_2$, $v_3$, $v_4$}{}
    \caption{$B_{\sigma}$ of $R_4$}\label{subfig:R4 nonmatching}
    \end{subfigure}\hfill
    \begin{subfigure}{0.24\linewidth}\centering
    \BackedgeGraphGeneral{1/3, 2/4}{4}{$v_1$, $v_3$, $v_2$, $v_4$}{}
    \caption{$B_{\sigma_1}$ of $R_4$}\label{subfig:R4 matching 1}
    \end{subfigure}\hfill
    \begin{subfigure}{0.24\linewidth}\centering
    \BackedgeGraphGeneral{1/4}{4}{$v_2$, $v_1$, $v_4$, $v_3$}{}
    \caption{$B_{\sigma_2}$ of $R_4$}\label{subfig:R4 matching 2}
    \end{subfigure}\hfill
    \begin{subfigure}{0.24\linewidth}\centering
    \BackedgeGraphGeneral{1/4, 2/3}{4}{$v_2$, $v_4$, $v_1$, $v_3$}{}
    \caption{$B_{\sigma_3}$ of $R_4$}\label{subfig:R4 matching 3}
    \end{subfigure}\hfill

    \vspace{5mm}

    \begin{subfigure}{0.9\linewidth}\centering
        \BackedgeGraphWithEllipsis{1/2, 2/3, 3/4, 4/5, 5/6, 6/7, 7/EL, ER/8, 8/9, 9/10, 10/11, 11/12, 12/13, 13/14}{14}{$v_1$, $v_2$, $v_3$, $v_4$, $v_5$, $v_6$, $v_7$, $v_{2k-5}$, $v_{2k-4}$, $v_{2k-3}$, $v_{2k-2}$, $v_{2k-1}$, $v_{2k}$ ,$v_{2k+1}$}{7}
        \caption{$B_{\sigma}$ of $R_{2k+1}$}\label{subfig:R2k+1 nonmatching}
    \end{subfigure}
    \vspace{2mm}
    \begin{subfigure}{0.9\linewidth}\centering
        \BackedgeGraphWithEllipsis{1/3, 2/5, 4/7, 6/EL, ER/9, 8/11, 10/13}{13}{$v_1$, $v_3$, $v_2$, $v_5$, $v_4$, $v_7$, $v_6$, $v_{2k-3}$, $v_{2k-4}$, $v_{2k-1}$, $v_{2k-2}$, $v_{2k+1}$, $v_{2k}$}{7}
        \caption{$B_{\sigma_1}$ of $R_{2k+1}$}\label{subfig:R2k+1 matching 1}
    \end{subfigure}
    \vspace{2mm}
    \begin{subfigure}{0.9\linewidth}\centering
        \BackedgeGraphWithEllipsis{1/4, 3/6, 5/EL, ER/8, 7/10, 9/12, 11/13}{13}{$v_2$, $v_1$, $v_4$, $v_3$, $v_6$, $v_5$, $v_{2k-4}$, $v_{2k-5}$, $v_{2k-2}$, $v_{2k-3}$, $v_{2k}$, $v_{2k-1}$, $v_{2k+1}$}{6}
        \caption{$B_{\sigma_2}$ of $R_{2k+1}$}\label{subfig:R2k+1 matching 2}
    \end{subfigure}

    \vspace{5mm}

    \begin{subfigure}{0.9\linewidth}\centering
        \BackedgeGraphWithEllipsis{1/2, 2/3, 3/4, 4/5, 5/6, 6/7, 7/EL, ER/8, 8/9, 9/10, 10/11, 11/12, 12/13, 13/14}{14}{$v_1$, $v_2$, $v_3$, $v_4$, $v_5$, $v_6$, $v_7$, $v_{2k-6}$, $v_{2k-5}$, $v_{2k-4}$, $v_{2k-3}$, $v_{2k-2}$, $v_{2k-1}$ ,$v_{2k}$}{7}
        \caption{$B_{\sigma}$ of $R_{2k}$}\label{subfig:R2k nonmatching}
    \end{subfigure}
    \vspace{2mm}
    \begin{subfigure}{0.9\linewidth}\centering
        \BackedgeGraphWithEllipsis{1/3, 2/5, 4/7, 6/EL, ER/9, 8/11, 10/13, 12/14}{14}{$v_1$, $v_3$, $v_2$, $v_5$, $v_4$, $v_7$, $v_6$, $v_{2k-5}$, $v_{2k-6}$, $v_{2k-3}$, $v_{2k-4}$, $v_{2k-1}$, $v_{2k-2}$, $v_{2k}$}{7}
        \caption{$B_{\sigma_1}$ of $R_{2k}$}\label{subfig:R2k matching 1}
    \end{subfigure}
    \vspace{2mm}
    \begin{subfigure}{0.9\linewidth}\centering
        \BackedgeGraphWithEllipsis{1/4, 3/6, 5/EL, ER/8, 7/10, 9/12}{12}{$v_2$, $v_1$, $v_4$, $v_3$, $v_6$, $v_5$, $v_{2k-4}$, $v_{2k-5}$, $v_{2k-2}$, $v_{2k-3}$, $v_{2k}$, $v_{2k-1}$}{6}
        \caption{$B_{\sigma_2}$ of $R_{2k}$}\label{subfig:R2k matching 2}
    \end{subfigure}
    
    \caption{Backedge graphs of $R_n$ for different values of $n$}
    \label{fig: backedge graphs of paths as paths and matchings}
\end{figure}

We begin by proving that each of the graphs in the set described in \Cref{theorem:main} is minimally nonmatching.
It is easy, though somewhat tedious, to check that this is true of $T_5, M_7$, $P_7^-$, and the instances of a) to e) where the involved induced paths have length below five.
To prove that it is also true of the remaining elements of the set, we will use the precise characterization of the matching orderings of tournaments with induced monotone paths as backedge graphs from a combination of Corollary 5.6 and some comments made on page 17 of \cite{liuVariousTheoremsTournaments2012}.
\begin{lemma}
    \label{lemma: matching orderings of induced paths}
    Let $n$ be a positive integer and $R_n$ be the tournament with an induced monotone path of on $n$ vertices as a backedge graph.
    Additionally suppose that the vertex set of $R_n$ is $v_1,\dots,v_n$ and the ordering of $V(R_n)$ that induces its backedge graph that is an induced monotone path is $\sigma=v_1,\dots,v_n$.
    \begin{itemize}
        \item If $n=1$, then $R_n$ has exactly one matching ordering, namely $\sigma$.
        \item If $n=2$, then $R_n$ has exactly two matching orderings, namely $\sigma_1=v_1,v_2$ and $\sigma_2=v_2,v_1$.
        \item If $n=3$, then $R_n$ has exactly three matching orderings, namely $\sigma_1=v_1,v_3,v_2$, $\sigma_2=v_2,v_1,v_3$, and $\sigma_3=v_3,v_2,v_1$.
        \item If $n=4$, then $R_n$ has exactly three matching orderings, namely $\sigma_1=v_1,v_3,v_2,v_4$, $\sigma_2=v_2,v_1,v_4,v_3$, and $\sigma_3=v_2,v_4,v_1,v_3$.
        \item If $n\geq 5$, then $R_n$ has exactly two matching orderings.
        
        When $n=2k+1$ is odd, these are
        \begin{equation}
            \label{eq: matching ordering induced paths odd length at least five v2k has degree zero}
            \sigma_1=v_1,v_3,v_2,v_5,v_4,\dots,v_{2k+1},v_{2k}
        \end{equation}
        and
        \begin{equation}
            \label{eq: matching ordering induced paths odd length at least five v1 has degree zero}
            \sigma_2=v_2,v_1,v_4,v_3,\dots,v_{2k},v_{2k-1},v_{2k+1}.
        \end{equation}
        
        When $n=2k$ is even, the matching orderings of $R_n$ are
        \begin{equation}
            \label{eq: matching ordering induced paths even length at least five everyone has degree one}
            \sigma_1=v_1,v_3,v_2,v_5,v_4,\dots,v_{2k-1},v_{2k-2},v_{2k}
        \end{equation}
        and
        \begin{equation}
            \label{eq: matching ordering induced paths even length at least five two degree zero vertices}
            \sigma_2=v_2,v_1,v_4,v_3,\dots,v_{2k},v_{2k-1}.
        \end{equation}
    \end{itemize}
    We give each of the backedge graphs involved in this lemma in \Cref{fig: backedge graphs of paths as paths and matchings}.
\end{lemma}

We are particularly interested in $R_n$ for $n$ odd and at least five where the following observation is useful.

\begin{observation}
    \label{obs: precise matching orderings induced paths odd length at least five}

    Let $n=2k+1$ be an odd integer that is greater than four and $R_n$ be the tournament with an induced monotone path on $n$ vertices as a backedge graph.
    Additionally suppose that the vertex set of $R_n$ is $v_1,\dots,v_n$ and the ordering of $V(R_n)$ that induces its backedge graph that is an induced monotone path is $\sigma=v_1,\dots,v_n$.
    Let $\sigma_1$ and $\sigma_2$ be the two matching orderings of $R_n$ from \Cref{lemma: matching orderings of induced paths}.

    Every vertex except for $v_{2k+1}$ has degree one in $B_{\sigma_1}(R_n)$ and $v_{2k+1}$ has degree zero in this backedge graph of $R_n$.
    Similarly, every vertex except for $v_{1}$ has degree one in $B_{\sigma_2}(R_n)$ and $v_{1}$ has degree zero in this backedge graph of $R_n$.
\end{observation}

\subsection{Nonmatching-ness}
We now examine the tournaments with backedge graphs of types a) to e) where the involved paths have length at least five and argue that they are nonmatching.

Suppose such a tournament $T$ on vertices $v_1,\dots,v_{2k+1}$ and $w_1,w_2$ (in the case of $M_{1,1,is}^{(2k+1)}$ and $M_{1,1,adj}^{(2k+1)}$), or $w_1,\dots,w_4$ ($M_{2,2}^{(2k+1)}$ and $M_{1,2}^{(2k+1)}$, and $M_{2,1}^{(2k+1)}$) has a matching backedge graph $B_{\sigma}$ and let $\sigma=u_1,\dots,u_{2k+3}$, resp. $u_1,\dots,u_{2k+5}$, be an ordering inducing such a backedge graph.

Since the set of tournaments we are examining is closed under reversing all arcs, the mirror of a matching backedge graph of a tournament $T$ is a matching backedge graph of $T^*$, and by \Cref{lemma: matching orderings of induced paths} $v_{1}$ precedes $v_{2}$ in $\sigma$ or $v_{2k+1}$ follows $v_{2k}$ in $\sigma$, we may assume without loss of generality that $v_1$ is before $v_2$ in $\sigma$.
In particular, the ordering of $v_1,\dots,v_{2k+1}$ that is induced by $\sigma$ must be $\sigma_1$ from \Cref{eq: matching ordering induced paths odd length at least five v2k has degree zero}.
Thus, $v_1$ has degree one in the backedge graph induced by $v_1,\dots,v_{2k+1}$ (cf. \Cref{obs: precise matching orderings induced paths odd length at least five}).
% We consider two cases, depending on whether $v_1$ is before or after $v_2$ in $\sigma$.
% If $v_1$ is before $v_2$ in this ordering, then the ordering of $v_1,\dots,v_{2k+1}$ that is induced by $\sigma$ must be that in \Cref{eq: matching ordering induced paths odd length at least five v2k has degree zero}.
% In particular, $v_1$ has degree one in the backedge graph induced by $v_1,\dots,v_{2k+1}$ (cf. \Cref{obs: precise matching orderings induced paths odd length at least five}).
We now consider five different cases depending on the type of $T$:
\begin{itemize}
    \item[a)]
    If $T$ is $M_{2,2}^{(2k+1)}$, then at most one of $w_3$ and $w_4$ is after $v_{2k+1}$.
    In particular, one of $w_3,w_4$ is before $v_{2k+1}$ and therefore before $v_{2k}$.
    This is impossible, since it would result in $v_{2k}$ having degree two in $B_{\sigma}$.
    \item[b)]
    If $T$ is $M_{1,1,is}^{(2k+1)}$, then $w_1$ and $w_2$ must both be before $v_1$ in $\sigma$.
    But then they are also before $v_{2k+1}$ in $\sigma$ which implies the existence of edges between $v_{2k+1}$ and both $w_1$ and $w_2$ in $B_{\sigma}$.
    This is impossible, since it would result in $v_{2k+1}$ having degree two in $B_{\sigma}$.
    \item[c)]
    If $T$ is $M_{1,1,adj}^{(2k+1)}$, then $w_2$ must be after $v_{2k}$ in $\sigma$.
    It is therefore also after $v_{2k+1}$ and $w_2v_{2k+1}$ is an edge in $B_{\sigma}$.
    It follows that $w_1$ is after $w_2$ in $\sigma$, since they would otherwise be adjacent in $B_{\sigma}$.
    But then, $w_1$ is also after $v_{2k+1}$ in $\sigma$
    This is impossible, since it would result in $v_{2k+1}$ having degree two in $B_{\sigma}$.
    \item[d)]
    If $T$ is $M_{1,2}^{(2k+1)}$, then at most one of $w_3$ and $w_4$ is after $v_{2k+1}$ in $\sigma$.
    Much as in case a), this results in a contradiction to the fact that $v_{2k}$ has degree one in $B_{\sigma}$.
    \item[e)]
    If $T$ is $M_{2,1}^{(2k+1)}$, then $w_1$ is before $v_1$ in $\sigma$, and $v_{2k+1}, v_{2k}, w_4$ appear in this order in $\sigma$. But now $w_4$ has degree at least two in $\sigma$, a contradiction. 

    %Sophie: don't think you can use symmetry here since symmetry was broken earlier. 
    %then $T$ has a backedge graph isomorphic to the mirror of a backedge graph of $M_{1,2}^{(2k+1)}$.
    %Since $M_{1,2}^{(2k+1)}$ is a nonmatching tournament, the same must be true of $T$.
\end{itemize}

\subsection{Minimality}
\begin{figure}[ht!]
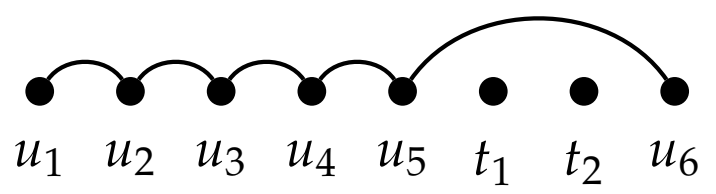
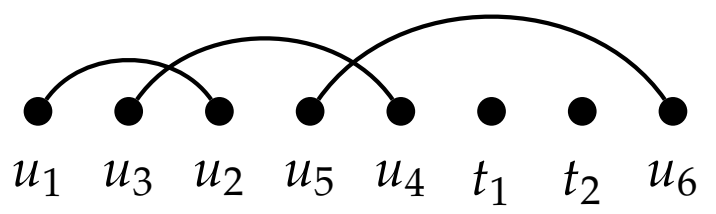

    \centering
    % \begin{subfigure}{0.9\linewidth}\centering
    %     \BackedgeGraphWithEllipsis{1/2, 2/3, 3/4, 4/5, 5/6, 6/7, 7/EL, ER/8, 8/9, 9/10, 10/11, 11/12, 12/13, 13/16}{16}{$u_1$, $u_2$, $u_3$, $u_4$, $u_5$, $u_6$, $u_7$, $u_{2l-6}$, $u_{2l-5}$, $u_{2l-4}$, $u_{2l-3}$, $u_{2l-2}$, $u_{2l-1}$, $t_1$, $t_2$, $u_{2l}$}{7}
    %     \caption{$B_{\sigma}$ of $T$ in \Cref{lemma: even path with two isolated vertices on one end}}
    % \end{subfigure}
    % \begin{subfigure}{0.9\linewidth}\centering
    %     \BackedgeGraphWithEllipsis{1/3, 2/5, 4/7, 6/EL, ER/9, 8/11, 10/13, 12/16}{16}{$u_1$, $u_3$, $u_2$, $u_5$, $u_4$, $u_7$, $u_6$, $u_{2l-5}$, $u_{2l-6}$, $u_{2l-3}$, $u_{2l-4}$, $u_{2l-1}$, $u_{2l-2}$, $t_1$, $t_2$, $u_{2l}$}{7}
    %     \caption{$B_{\sigma'}$ of $T$ in \Cref{lemma: even path with two isolated vertices on one end}}
    % \end{subfigure}
    \begin{subfigure}{0.49\linewidth}\centering
        \BackedgeGraphLabeled{1/2, 2/3, 3/4, 4/5, 5/8}{8}{$u_1$, $u_2$, $u_3$, $u_4$, $u_5$, $t_1$, $t_2$, $u_{6}$}
        \caption{$B_{\sigma}$ of $T$ in \Cref{lemma: even path with two isolated vertices on one end}}
    \end{subfigure}
    \begin{subfigure}{0.49\linewidth}\centering
        \BackedgeGraphLabeled{1/3, 2/5, 4/8}{8}{$u_1$, $u_3$, $u_2$, $u_5$, $u_4$, $t_1$, $t_2$, $u_{6}$}
        \caption{$B_{\sigma'}$ of $T$ in \Cref{lemma: even path with two isolated vertices on one end}}
    \end{subfigure}
    \caption{Backedge graphs of the tournaments in \Cref{lemma: even path with two isolated vertices on one end}}
    \label{fig: even path with two isolated vertices on one end}
\end{figure}

Secondly, we argue that every proper subtournament of any of the tournaments $M_{2,2}^{(2k+1)}$, $M_{1,1,is}^{(2k+1)}$, $M_{1,1,adj}^{(2k+1)}$, $M_{1,2}^{(2k+1)}$, $M_{2,1}^{(2k+1)}$ has a matching backedge graph.
The following sequence of results related to tournaments with backedge graphs largely composed of induced monotone paths is of some use in this respect.
\begin{lemma}
    \label{lemma: even path with two isolated vertices on one end}
    Let $T$ be a tournament on the vertex set $u_1,\dots,u_{2l},t_1,t_2$ with a backedge graph $B_{\sigma}$ induced by the ordering $\sigma=u_1,\dots,u_{2l-1},t_1,t_2,u_{2l}$ and $E(G)\coloneqq\{u_{i}u_{i+1}\mid 1\leq i<2l\}$.
    Then the ordering
    \[\sigma'=u_1,u_3,u_2,\dots,u_{2l-3},u_{2l-4},u_{2l-1},u_{2l-2},t_1,t_2,u_{2l}\]
    induces a matching backedge graph of $T$.

    We give an example with $l=3$ in \Cref{fig: even path with two isolated vertices on one end}.
\end{lemma}

\begin{proof}
    Note that the ordering of $u_1,\dots,u_{2l}$ induced by $\sigma'$ is a matching ordering of the subtournament of $T$ induced by this vertex set (due to \Cref{lemma: matching orderings of induced paths}).
    Moreover, the positions of $t_1,t_2$ relative to vertices in $V(T)\setminus \{t_1,t_2\}$ is the same in $\sigma$ as in $\sigma'$ (and these vertices therefore still have degree at most one in $B_{\sigma'}$).
    Thus, $\sigma'$ is a matching ordering of $T$.
\end{proof}

\begin{figure}[ht]
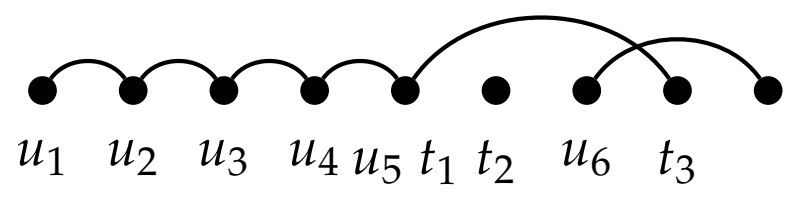
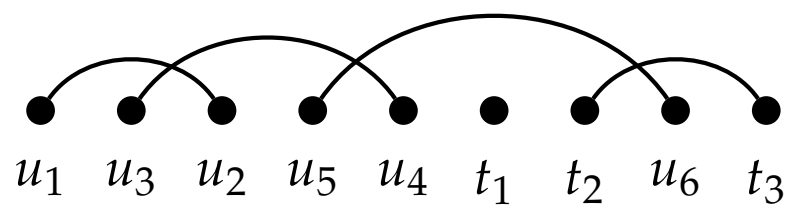

    \centering
    % \begin{subfigure}{0.9\linewidth}\centering
    %     \BackedgeGraphWithEllipsis{1/2, 2/3, 3/4, 4/5, 5/6, 6/7, 7/EL, ER/8, 8/9, 9/10, 10/11, 11/12, 12/13, 13/16, 15/17}{17}{$u_1$, $u_2$, $u_3$, $u_4$, $u_5$, $u_6$, $u_7$, $u_{2l-6}$, $u_{2l-5}$, $u_{2l-4}$, $u_{2l-3}$, $u_{2l-2}$, $u_{2l-1}$, $t_1$, $t_2$, $u_{2l}$, $t_3$}{7}
    %     \caption{$B_{\sigma}$ of $T$ in \Cref{lemma: even path with two isolated vertices on one end plus an incoming edge}}
    % \end{subfigure}
    % \begin{subfigure}{0.9\linewidth}\centering
    %     \BackedgeGraphWithEllipsis{1/3, 2/5, 4/7, 6/EL, ER/9, 8/11, 10/13, 12/16, 15/17}{17}{$u_1$, $u_3$, $u_2$, $u_5$, $u_4$, $u_7$, $u_6$, $u_{2l-5}$, $u_{2l-6}$, $u_{2l-3}$, $u_{2l-4}$, $u_{2l-1}$, $u_{2l-2}$, $t_1$, $t_2$, $u_{2l}$, $t_3$}{7}
    %     \caption{$B_{\sigma'}$ of $T$ in \Cref{lemma: even path with two isolated vertices on one end plus an incoming edge}}
    % \end{subfigure}
    \begin{subfigure}{0.49\linewidth}\centering
        \BackedgeGraphLabeled{1/2, 2/3, 3/4, 4/5, 5/8, 7/9}{9}{$u_1$, $u_2$, $u_3$, $u_4$, $u_5$ $t_1$, $t_2$, $u_{6}$, $t_3$}
        \caption{$B_{\sigma}$ of $T$ in \Cref{lemma: even path with two isolated vertices on one end plus an incoming edge}}
    \end{subfigure}
    \begin{subfigure}{0.49\linewidth}\centering
        \BackedgeGraphLabeled{1/3, 2/5, 4/8, 7/9}{9}{$u_1$, $u_3$, $u_2$, $u_5$, $u_4$, $t_1$, $t_2$, $u_{6}$, $t_3$}
        \caption{$B_{\sigma'}$ of $T$ in \Cref{lemma: even path with two isolated vertices on one end plus an incoming edge}}
    \end{subfigure}
    \caption{Backedge graphs of the tournaments in \Cref{lemma: even path with two isolated vertices on one end plus an incoming edge}}
    \label{fig: even path with two isolated vertices on one end plus an incoming edge}
\end{figure}

\begin{lemma}
    \label{lemma: even path with two isolated vertices on one end plus an incoming edge}
    Let $T$ be a tournament on the vertex set $u_1,\dots,u_{2l},t_1,t_2$ with a backedge graph $B_{\sigma}$ induced by the ordering $\sigma=u_1,\dots,u_{2l-1},t_1,t_2,u_{2l},t_3$ and $E(G)\coloneqq\{u_{i}u_{i+1}\mid 1\leq i<2l\}\cup\{t_2t_3\}$.
    Then the ordering
    \[\sigma'=u_1,u_3,u_2,\dots,u_{2l-5},u_{2l-6},u_{2l-3},u_{2l-4},u_{2l-1},u_{2l-2},t_1,t_2,u_{2l},t_3\]
    induces a matching backedge graph of $T$.

    We give an example with $l=3$ in \Cref{fig: even path with two isolated vertices on one end plus an incoming edge}.
\end{lemma}

\begin{proof}
    Identical to the proof of \Cref{lemma: even path with two isolated vertices on one end}, further noting that the position of $t_1$ relative to $\{t_2,t_3\}$ is the same in $\sigma$ and $\sigma'$.
\end{proof}

\begin{figure}[ht]
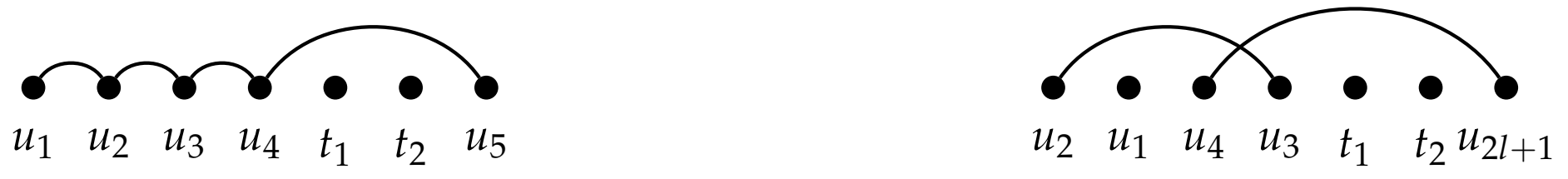

    \centering
    % \begin{subfigure}{0.9\linewidth}\centering
    %     \BackedgeGraphWithEllipsis{1/2, 2/3, 3/4, 4/5, 5/6, 6/EL, ER/7, 7/8, 8/9, 9/10, 10/11, 11/12, 12/15}{15}{$u_1$, $u_2$, $u_3$, $u_4$, $u_5$, $u_6$, $u_{2l-5}$, $u_{2l-4}$, $u_{2l-3}$, $u_{2l-2}$, $u_{2l-1}$, $u_{2l}$, $t_1$, $t_2$, $u_{2l+1}$}{6}
    %     \caption{$B_{\sigma}$ of $T$ in \Cref{lemma: odd path with two isolated vertices on one end}}
    % \end{subfigure}
    % \begin{subfigure}{0.9\linewidth}\centering
    %     \BackedgeGraphWithEllipsis{1/4, 3/6, 5/EL, ER/8, 7/10, 9/12, 11/15}{15}{$u_2$, $u_1$, $u_4$, $u_3$, $u_6$, $u_5$, $u_{2l-4}$, $u_{2l-5}$, $u_{2l-2}$, $u_{2l-3}$, $u_{2l}$, $u_{2l-1}$, $t_1$, $t_2$, $u_{2l+1}$}{6}
    %     \caption{$B_{\sigma'}$ of $T$ in \Cref{lemma: odd path with two isolated vertices on one end}}
    % \end{subfigure}
    \begin{subfigure}{0.49\linewidth}\centering
        \BackedgeGraphLabeled{1/2, 2/3, 3/4, 4/7}{7}{$u_1$, $u_2$, $u_3$, $u_4$, $t_1$, $t_2$, $u_{5}$}
        \caption{$B_{\sigma}$ of $T$ in \Cref{lemma: odd path with two isolated vertices on one end}}
    \end{subfigure}
    \begin{subfigure}{0.49\linewidth}\centering
        \BackedgeGraphLabeled{1/4, 3/7}{7}{$u_2$, $u_1$, $u_4$, $u_3$, $t_1$, $t_2$, $u_{2l+1}$}
        \caption{$B_{\sigma'}$ of $T$ in \Cref{lemma: odd path with two isolated vertices on one end}}
    \end{subfigure}
    \caption{Backedge graphs of the tournaments in \Cref{lemma: odd path with two isolated vertices on one end}}
    \label{fig: odd path with two isolated vertices on one end}
\end{figure}

\begin{lemma}
    \label{lemma: odd path with two isolated vertices on one end}
    Let $T$ be a tournament on the vertex set $u_1,\dots,u_{2l+1},t_1,t_2$ with a backedge graph $B_{\sigma}$ induced by the ordering $\sigma=u_1,\dots,u_{2l},t_1,t_2,u_{2l+1}$ and $E(G)\coloneqq\{u_{i}u_{i+1}\mid 1\leq i<2l+1\}$.
    Then the ordering
    \[\sigma'=u_2,u_1,u_4,u_3,\dots,u_{2l},u_{2l-1},t_1,t_2,u_{2l+1}\]
    induces a matching backedge graph of $T$.

    We give an example with $l=2$ in \Cref{fig: odd path with two isolated vertices on one end}.
\end{lemma}

\begin{proof}
    Identical to the proof of \Cref{lemma: even path with two isolated vertices on one end}.
\end{proof}

\begin{figure}[ht]
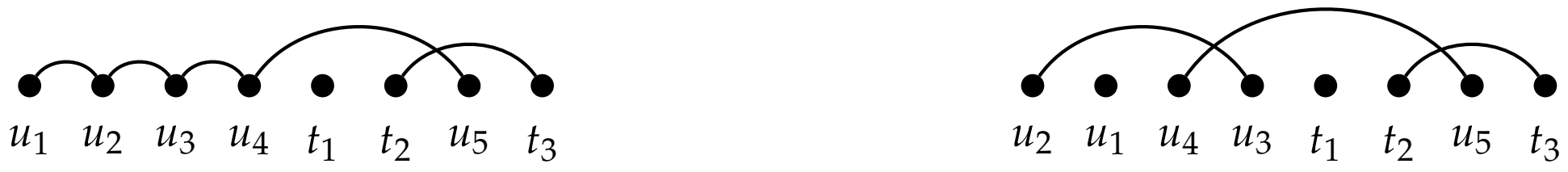

    \centering
    % \begin{subfigure}{0.9\linewidth}\centering
    %     \BackedgeGraphWithEllipsis{1/2, 2/3, 3/4, 4/5, 5/6, 6/EL, ER/7, 7/8, 8/9, 9/10, 10/11, 11/12, 12/15, 14/16}{16}{$u_1$, $u_2$, $u_3$, $u_4$, $u_5$, $u_6$, $u_{2l-5}$, $u_{2l-4}$, $u_{2l-3}$, $u_{2l-2}$, $u_{2l-1}$, $u_{2l}$, $t_1$, $t_2$, $u_{2l+1}$, $t_3$}{6}
    %     \caption{$B_{\sigma}$ of $T$ in \Cref{lemma: odd path with two isolated vertices on one end plus an incoming edge}}
    % \end{subfigure}
    % \begin{subfigure}{0.9\linewidth}\centering
    %     \BackedgeGraphWithEllipsis{1/4, 3/6, 5/EL, ER/8, 7/10, 9/12, 11/15, 14/16}{16}{$u_2$, $u_1$, $u_4$, $u_3$, $u_6$, $u_5$, $u_{2l-4}$, $u_{2l-5}$, $u_{2l-2}$, $u_{2l-3}$, $u_{2l}$, $u_{2l-1}$, $t_1$, $t_2$, $u_{2l+1}$, $t_3$}{6}
    %     \caption{$B_{\sigma'}$ of $T$ in \Cref{lemma: odd path with two isolated vertices on one end plus an incoming edge}}
    % \end{subfigure}
    \begin{subfigure}{0.49\linewidth}\centering
        \BackedgeGraphLabeled{1/2, 2/3, 3/4, 4/7, 6/8}{8}{$u_1$, $u_2$, $u_3$, $u_4$, $t_1$, $t_2$, $u_{5}$, $t_3$}
        \caption{$B_{\sigma}$ of $T$ in \Cref{lemma: odd path with two isolated vertices on one end plus an incoming edge}}
    \end{subfigure}
    \begin{subfigure}{0.49\linewidth}\centering
        \BackedgeGraphLabeled{1/4, 3/7, 6/8}{8}{$u_2$, $u_1$, $u_4$, $u_3$, $t_1$, $t_2$, $u_{5}$, $t_3$}
        \caption{$B_{\sigma'}$ of $T$ in \Cref{lemma: odd path with two isolated vertices on one end plus an incoming edge}}
    \end{subfigure}
    \caption{Backedge graphs of the tournaments in \Cref{lemma: odd path with two isolated vertices on one end plus an incoming edge}}
    \label{fig: odd path with two isolated vertices on one end plus an incoming edge}
\end{figure}

\begin{lemma}
    \label{lemma: odd path with two isolated vertices on one end plus an incoming edge}
    Let $T$ be a tournament on the vertex set $u_1,\dots,u_{2l+1},t_1,t_2,t_3$ with a backedge graph $B_{\sigma}$ induced by the ordering $\sigma=u_1,\dots,u_{2l},t_1,t_2,u_{2l+1},t_3$ and $E(G)\coloneqq\{u_{i}u_{i+1}\mid 1\leq i<2l+1\}\cup\{t_2t_3\}$.
    Then the ordering
    \[\sigma'=u_2,u_1,u_4,u_3,\dots,u_{2l},u_{2l-1},t_1,t_2,u_{2l+1},t_3\]
    induces a matching backedge graph of $T$.

    We give an example with $l=2$ in \Cref{fig: odd path with two isolated vertices on one end plus an incoming edge}.
\end{lemma}

\begin{proof}
    Identical to the proof of \Cref{lemma: even path with two isolated vertices on one end plus an incoming edge}.
\end{proof}

\begin{figure}[ht]
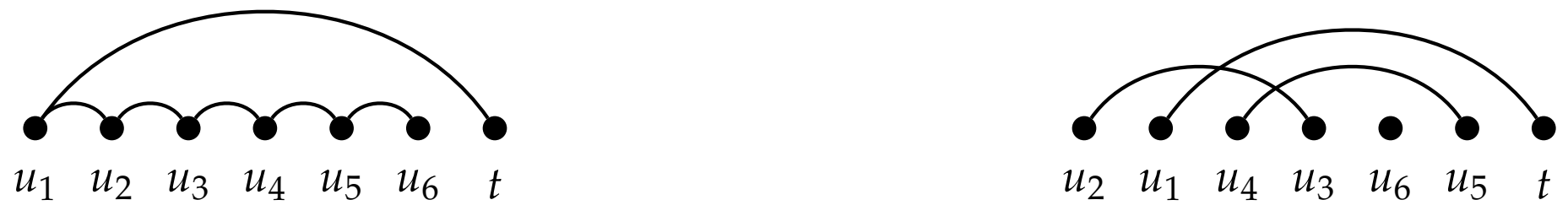

    \centering
    % \begin{subfigure}{0.9\linewidth}\centering
    %     \BackedgeGraphWithEllipsis{1/2, 2/3, 3/4, 4/5, 5/6, 6/EL, ER/7, 7/8, 8/9, 9/10, 10/11, 11/12, 1/13}{13}{$u_1$, $u_2$, $u_3$, $u_4$, $u_5$, $u_6$, $u_{2l-5}$, $u_{2l-4}$, $u_{2l-3}$, $u_{2l-2}$, $u_{2l-1}$, $u_{2l}$, $t$}{6}
    %     \caption{$B_{\sigma}$ of $T$ in \Cref{lemma: even path plus a long forward edge}}
    % \end{subfigure}
    % \begin{subfigure}{0.9\linewidth}\centering
    %     \BackedgeGraphWithEllipsis{1/4, 3/6, 5/EL, ER/8, 7/10, 9/12, 2/13}{13}{$u_2$, $u_1$, $u_4$, $u_3$, $u_6$, $u_5$, $u_{2l-4}$, $u_{2l-5}$, $u_{2l-2}$, $u_{2l-3}$, $u_{2l}$, $u_{2l-1}$, $t$}{6}
    %     \caption{$B_{\sigma'}$ of $T$ in \Cref{lemma: even path plus a long forward edge}}
    % \end{subfigure}
    \begin{subfigure}{0.49\linewidth}\centering
        \BackedgeGraphLabeled{1/2, 2/3, 3/4, 4/5, 5/6, 1/7}{7}{$u_1$, $u_2$, $u_3$, $u_4$, $u_5$, $u_6$, $t$}
        \caption{$B_{\sigma}$ of $T$ in \Cref{lemma: even path plus a long forward edge}}
    \end{subfigure}
    \begin{subfigure}{0.49\linewidth}\centering
        \BackedgeGraphLabeled{1/4, 3/6, 2/7}{7}{$u_2$, $u_1$, $u_4$, $u_3$, $u_6$, $u_5$, $t$}
        \caption{$B_{\sigma'}$ of $T$ in \Cref{lemma: even path plus a long forward edge}}
    \end{subfigure}
    \caption{Backedge graphs of the tournaments in \Cref{lemma: even path plus a long forward edge}}
    \label{fig: even path plus a long forward edge}
\end{figure}

\begin{lemma}
    \label{lemma: even path plus a long forward edge}
    Let $T$ be a tournament on the vertex set $u_1,\dots,u_{2l},t$ with a backedge graph $B_{\sigma}$ induced by the ordering $u_1,\dots,u_{2l-1},u_{2l},t$ and $E(G)\coloneqq\{u_{i}u_{i+1}\mid 1\leq i<2l\}\cup\{u_1t\}$.
    Then the ordering
    \[\sigma'=u_2,u_1,u_4,u_3,\dots,u_{2l},u_{2l-1},t\]
    induces a matching backedge graph of $T$.

    We give an example with $l=3$ in \Cref{fig: even path plus a long forward edge}.
\end{lemma}

\begin{proof}
    Note that the ordering of $u_1,\dots,u_{2l}$ induced by $\sigma'$ induces a matching backedge graph of the subtournament of $T$ induced by this vertex set where $u_1$ has degree zero (due to \Cref{lemma: matching orderings of induced paths}).
    Moreover, the position of $t$ relative to vertices in $V(T)\setminus \{t\}$ is the same in $\sigma$ as in $\sigma'$ (and therefore the only edge incident to $t$ in $B_{\sigma'}$ is $u_1t$).
    Thus, $\sigma'$ is a matching ordering of $T$.
\end{proof}

\begin{figure}[ht]
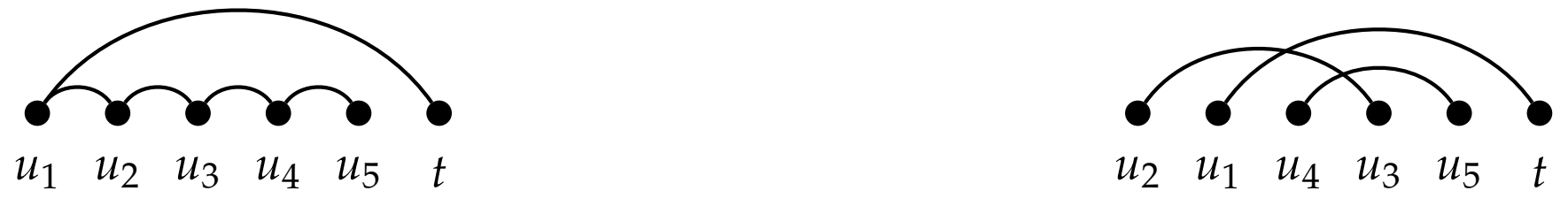

    \centering
    % \begin{subfigure}{0.9\linewidth}\centering
    %     \BackedgeGraphWithEllipsis{1/2, 2/3, 3/4, 4/5, 5/6, 6/EL, ER/7, 7/8, 8/9, 9/10, 10/11, 11/12, 12/13, 1/14}{14}{$u_1$, $u_2$, $u_3$, $u_4$, $u_5$, $u_6$, $u_{2l-5}$, $u_{2l-4}$, $u_{2l-3}$, $u_{2l-2}$, $u_{2l-1}$, $u_{2l}$, $u_{2l+1}$, $t$}{6}
    %     \caption{$B_{\sigma}$ of $T$ in \Cref{lemma: odd path plus a long forward edge}}
    % \end{subfigure}
    % \begin{subfigure}{0.9\linewidth}\centering
    %     \BackedgeGraphWithEllipsis{1/4, 3/6, 5/EL, ER/8, 7/10, 9/12, 11/13, 2/14}{14}{$u_2$, $u_1$, $u_4$, $u_3$, $u_6$, $u_5$, $u_{2l-4}$, $u_{2l-5}$, $u_{2l-2}$, $u_{2l-3}$, $u_{2l}$, $u_{2l-1}$, $u_{2l+1}$, $t$}{6}
    %     \caption{$B_{\sigma'}$ of $T$ in \Cref{lemma: odd path plus a long forward edge}}
    % \end{subfigure}
    \begin{subfigure}{0.49\linewidth}\centering
        \BackedgeGraphLabeled{1/2, 2/3, 3/4, 4/5, 1/6}{6}{$u_1$, $u_2$, $u_3$, $u_4$, $u_5$, $t$}
        \caption{$B_{\sigma}$ of $T$ in \Cref{lemma: odd path plus a long forward edge}}
    \end{subfigure}
    \begin{subfigure}{0.49\linewidth}\centering
        \BackedgeGraphLabeled{1/4, 3/5, 2/6}{6}{$u_2$, $u_1$, $u_4$, $u_3$, $u_5$, $t$}
        \caption{$B_{\sigma'}$ of $T$ in \Cref{lemma: odd path plus a long forward edge}}
    \end{subfigure}
    \caption{Backedge graphs of the tournaments in \Cref{lemma: odd path plus a long forward edge}}
    \label{fig: odd path plus a long forward edge}
\end{figure}

\begin{lemma}
    \label{lemma: odd path plus a long forward edge}
    Let $T$ be a tournament on the vertex set $u_1,\dots,u_{2l+1},t$ with a backedge graph $B_{\sigma}$ induced by the ordering $u_1,\dots,u_{2l},u_{2l+1},t$ and $E(G)\coloneqq\{u_{i}u_{i+1}\mid 1\leq i<2l+1\}\cup\{u_1t\}$.
    Then the ordering
    \[\sigma'=u_2,u_1,u_4,u_3,\dots,u_{2l},u_{2l-1},u_{2l+1},t\]
    induces a matching backedge graph of $T$.

    We give an example with $l=2$ in \Cref{fig: odd path plus a long forward edge}.
\end{lemma}

\begin{proof}
    Identical to the proof of \Cref{lemma: even path plus a long forward edge}.
\end{proof}

Finally, we carry out our main goal for this part of the paper, i.e., show that every proper subtournament of any of $M_{2,2}^{(2k+1)}$, $M_{1,1,is}^{(2k+1)}$, $M_{1,1,adj}^{(2k+1)}$, $M_{1,2}^{(2k+1)}$, $M_{2,1}^{(2k+1)}$ has a matching backedge graph.
To do so, we consider each possible type of tournament and each deletion in each tournament, showing that the corresponding subtournament has a matching ordering.
Note, however, that $M_{2,1}^{(2k+1)}$ may be ignored, since its backedge graph from \Cref{def: infinite classes minimal nonmatching} is the mirror of that of $M_{1,2}^{(2k+1)}$.
\begin{itemize}
    \item [a)]
    Suppose $T$ is $M_{2,2}^{(2k+1)}$.
    % on vertices $v_1,\dots,v_{2k+1}$ and $w_1,\dots,w_4$.
    \begin{enumerate}
        \item $T'\coloneqq T\setminus w$ for some $w\in\{w_1,\dots,w_4\}$.
            Note that, up to isomorphism and reversing all arcs, we may assume that $w=w_1$.
            In this case, the ordering $w_2,v_1,v_2,\dots,v_{2k},w_3,w_4,v_{2k+1}$ shows that \Cref{lemma: even path with two isolated vertices on one end} applies to $T'$ and it therefore has a matching ordering.

        \item $T'\coloneqq T\setminus v$ for $v=v_1$ or $v=v_{2k+1}$.
            Note that, up to reversing all arcs, we may assume $v=v_1$.
            In this case, the ordering $v_2,\dots,v_{2k},w_3,w_4,v_{2k+1}$ shows that \Cref{lemma: odd path with two isolated vertices on one end} applies to $T'\setminus\{w_1,w_2\}$ and therefore this subtournament of $T'$ has a matching ordering $\sigma'$.
            The ordering
            \[w_1,w_2,\sigma'(1),\dots,\sigma'(2k+3)\]
            is a matching ordering of $T'$, since $w_1,w_2$ have degree zero in $B_{\sigma}(T')$ and the relative positions of $w_1,w_2$ and $T'\setminus\{w_1,w_2\}$ are the same as in $\sigma$.

        \item $T'\coloneqq T\setminus v_j$ for $1<j<2k+1$.
            The orderings $v_1,w_1,w_2,v_2,\dots,v_{j-1}$ and $v_{j+1},\dots,v_{2k},w_3,w_4,v_{2k+1}$ show that \Cref{lemma: even path with two isolated vertices on one end} or \Cref{lemma: odd path with two isolated vertices on one end} applies to the corresponding induced subtournaments of $T'$ (or those obtained from them by reversing all arcs).
            % $T'[v_1,w_1,w_2,v_2,\dots,v_{j-1}]$ and $T'[v_{j+1},\dots,v_{2k},w_3,w_4,v_{2k+1}]$ (or $T'[v_{j+1},\dots,v_{2k},w_3,w_4,v_{2k+1}]^*$).
            Therefore, these subtournaments of $T'$ have matching orderings $\sigma'$ and $\sigma''$ respectively.
            The ordering
            \[\sigma'(1),\dots,\sigma'(j+1),\sigma''(1),\dots,\sigma''(2k+2-j)\]
            is a matching ordering of $T'$, since the relative positions of vertices in $\{v_1,w_1,w_2,v_2,\dots,v_{j-1}\}$ and $\{v_{j+1},\dots,v_{2k},w_3,w_4,v_{2k+1}\}$ are the same as in $\sigma$.
    \end{enumerate}
    
    \item [b)]
    Suppose $T$ is $M_{1,1,is}^{(2k+1)}$.
    % on vertices $v_1,\dots,v_{2k+1}$ and $w_1,w_2$.
    \begin{enumerate}
        \item $T'\coloneqq T\setminus w$ for $w\in\{w_1,w_2\}$.
            Note that up to reversing all arcs, we may assume $w=w_2$.
            In this case, the ordering $v_1,v_2,\dots,v_{2k+1},w_2$ shows that \Cref{lemma: odd path plus a long forward edge} applies to $T'$ and it therefore has a matching ordering.
        
        \item $T'\coloneqq T\setminus v$ for $v\in\{v_1,v_{2k+1}\}$.
            Note that up to reversing all arcs, we may assume $v=v_{2k+1}$.
            In this case, the ordering $v_1,v_2,\dots,v_{2k},w_2$ shows that \Cref{lemma: even path plus a long forward edge} applies to $T'\setminus w_1$ and therefore this subtournament of $T'$ has a matching ordering $\sigma'$.
            The ordering
            \[w_1,\sigma'(1),\dots,\sigma'(2k+1)\]
            is a matching ordering of $T'$, since $w_1$ has degree zero in $B_{\sigma}(T')$ and the relative positions of $w_1$ and $T'\setminus w_1$ are the same as in $\sigma$.
        
        \item $T'\coloneqq T\setminus v_j$ for $1<j<2k+1$.
            The orderings $v_1,v_2,v_3,\dots,v_{j-1},w_2$ and $w_1,v_{j+1},\dots,v_{2k+1}$ show that \Cref{lemma: even path plus a long forward edge} or \Cref{lemma: odd path plus a long forward edge} applies to the corresponding subtournaments of $T'$ (or those obtained from them by reversing all arcs).
            % $T'[v_1,v_2,v_3,\dots,v_{j-1},w_2]$ and $T'[w_1,v_{j+1},\dots,v_{2k+1}]$ (or $T'[w_1,v_{j+1},\dots,v_{2k+1}]^*$).
            Therefore, these subtournaments of $T'$ have matching orderings $\sigma'$ and $\sigma''$ where $w_2$ and $w_1$ are last and first respectively.
            The ordering
            \[\sigma''(1)=w_1,\sigma'(1),\dots,\sigma'(j-1),\sigma''(2),\dots,\sigma''(2k+2-j),\sigma'(j)=w_2\]
            is a matching ordering of $T'$, since the relative positions of vertices in $\{v_1,v_2,v_3,\dots,v_{j-1},w_2\}$ and $\{w_1,v_{j+1},\dots,v_{2k+1}\}$ are the same as in $\sigma$.
    \end{enumerate}
    
    \item [c)]
    Suppose $T$ is $M_{1,1,adj}^{(2k+1)}$.
    % a tournament of type c) on vertices $v_1,\dots,v_{2k+1}$ and $w_1,w_2$.
    \begin{enumerate}
        \item $T'\coloneqq T\setminus w$ for some $w\in\{w_1,w_2\}$.
            Note that, up to reversing all arcs, we may assume $w=w_{1}$.
            In this case, the ordering $v_1,\dots,v_{2k+1},w_2$ induces a backedge graph of $T'$ that is an induced monotone path and has a matching ordering by \Cref{lemma: matching orderings of induced paths}.
        
        \item $T'\coloneqq T\setminus v$ for $v\in\{v_1,v_{2k+1}\}$.
            Note that, up to reversing all arcs, we may assume $v=v_{1}$.
            In this case, the ordering $w_1,v_2,v_3,\dots,v_{2k+1},w_2$ shows that \Cref{lemma: even path plus a long forward edge} applies to $T'$ and it therefore has a matching ordering.
        
        \item $T'\coloneqq T\setminus v_j$ for $1<j<2k+1$.
            The orderings $w_1,v_1,v_2,v_3,\dots,v_{j-1}$ and $v_{j+1},\dots,v_{2k+1},w_2$ show that $T'[w_1,v_1,v_2,v_3,\dots,v_{j-1}]$ and $T'[v_{j+1},\dots,v_{2k+1},w_2]$ have backedge graphs that are induced monotone paths.
            In particular, by \Cref{obs: precise matching orderings induced paths odd length at least five}, these subtournaments of $T'$ have matching orderings $\sigma'$ and $\sigma''$ that induce backedge graphs where $w_1$ and $w_2$ have degree zero respectively.
            The ordering
            \[\sigma'(1),\sigma'(2),\dots,\sigma'(j),\sigma''(1),\dots,\sigma''(2k+3-j)\]
            is a matching ordering of $T'$, since the relative positions of vertices in $\{w_1,v_1,v_2,v_3,\dots,v_{j-1}\}$ and $\{v_{j+1},\dots,v_{2k+1},w_2\}$ are the same as in $\sigma$ and the only edge between these two vertex sets in $B_{\sigma}$ is $w_1w_2$.
    \end{enumerate}
    
    \item [d)]
    Suppose $T$ is $M_{1,2}^{(2k+1)}$.
    % a tournament of type d) on vertices $v_1,\dots,v_{2k+1}$ and $w_1,\dots,w_4$.
    \begin{enumerate}
        \item $T'\coloneqq T\setminus w_1$.
            The ordering $v_1,v_2,\dots,v_{2k},w_2,w_3,v_{2k+1}$ shows that \Cref{lemma: odd path with two isolated vertices on one end} applies to the corresponding subtournament of $T'$ and it therefore has a matching ordering $\sigma'$.
            % $T'[v_1,v_2,\dots,v_{2k},w_2,w_3,v_{2k+1}]$ and it therefore has a matching ordering $\sigma'$.
            The ordering
            \[\sigma'(1),\sigma'(2),\dots,\sigma'(2k+3),w_4\]
            is a matching ordering of $T'$, since the relative positions of vertices in $V(T'\setminus w_4)$ and $w_4$ are the same as in $\sigma$.

        \item $T'\coloneqq T\setminus w$ for $w\in\{w_2,w_3\}$.
            Note that we may assume $w=w_2$.
            In this case, the ordering $w_1,v_1,v_2,\dots,v_{2k},v_{2k+1},w_3,w_4$ shows that \Cref{lemma: odd path plus a long forward edge} applies to $T'$ and it therefore has a matching ordering.
        
        \item $T'\coloneqq T\setminus w_4$.
            The ordering $w_1,v_1,v_2,\dots,v_{2k},w_2,w_3,v_{2k+1}$ shows that \Cref{lemma: even path with two isolated vertices on one end} applies to $T'$ and it therefore has a matching ordering.
        
        \item $T'\coloneqq T\setminus v_1$.
            The ordering $v_2,\dots,v_{2k},w_2,w_3,v_{2k+1}$ shows that \Cref{lemma: even path with two isolated vertices on one end} applies to the corresponding subtournament of $T'$ and it therefore has a matching ordering $\sigma'$.
            % $T'[v_2,\dots,v_{2k},w_2,w_3,v_{2k+1}]$ and it therefore has a matching ordering $\sigma'$.
            The ordering
            \[w_1,\sigma'(1),\sigma'(2),\dots,\sigma'(2k+2),w_4\]
            is a matching ordering of $T'$, since the relative positions of vertices in $V(T'\setminus \{w_1,w_4\})$ and $\{w_1,w_4\}$ are the same as in $\sigma$.

        \item $T'\coloneqq T\setminus v_{2k+1}$.
            The ordering $w_1,v_1,v_2,\dots,v_{2k},w_4$ shows that \Cref{lemma: odd path plus a long forward edge} applies to the corresponding subtournament of $T'$ and it therefore has a matching ordering $\sigma'$ where $w_4$ is last.
            % $T'[w_1,v_1,v_2,\dots,v_{2k},w_4]$ and it therefore has a matching ordering $\sigma'$ where $w_4$ is last.
            The ordering
            \[\sigma'(1),\sigma'(2),\dots,\sigma'(2k+1),w_2,w_3,\sigma'(2k+2)=w_4\]
            is a matching ordering of $T'$, since the relative positions of vertices in $V(T'\setminus \{w_2,w_3\})$ and $\{w_2,w_3\}$ are the same as in $\sigma$.
        
        \item $T'\coloneqq T\setminus v_j$ for $1<j<2k+1$.
            The orderings $w_1,v_1,v_2,\dots,v_{j-1},w_4$ and $v_{j+1},\dots,v_{2k},w_2,w_3,v_{2k+1}$ show that \Cref{lemma: even path plus a long forward edge} or \Cref{lemma: odd path plus a long forward edge} and \Cref{lemma: even path with two isolated vertices on one end} or \Cref{lemma: odd path with two isolated vertices on one end} apply to the corresponding subtournaments of $T'$ which we denote by $T'_1$ and $T'_2$ respectively.
            % $T'[w_1,v_1,v_2,\dots,v_{j-1},w_4]$ and $T'[v_{j+1},\dots,v_{2k},w_2,w_3,v_{2k+1}]$ respectively.
            Therefore, $T'_1$ has a matching ordering $\sigma'$ where $w_4$ is last and $T'_2$ has a matching ordering $\sigma''$.
            The ordering
            \[\sigma'(1),\dots,\sigma'(j),\sigma''(1),\dots,\sigma''(2k+3-j),\sigma'(j+1)=w_4\]
            is a matching ordering of $T'$, since the relative positions of vertices in $V(T'_1)$ and $V(T'_2)$ are the same as in $\sigma$.
    \end{enumerate}
\end{itemize}

\subsection{Incomparability}
\begin{lemma}\label{lemma:incomparable}
    If $T$ and $T'$ are tournaments of types a) to e), or $T_5$, $M_7$, $P_7^-$, and $T'$ is a subtournament of $T$, then $T=T'$.
\end{lemma}
\begin{proof}
    Since every tournament in the classes a--e, and $T_5$, $P_7^-$, and $M_7$, is minimal nonmatching, no tournament in this list can contain another as a proper subtournament. It therefore suffices to prove that the listed tournaments are pairwise non-isomorphic.

    The tournament $P_7^-$ has six vertices. However,
    $T_5$ has five vertices, $M_7$ has seven vertices, and the five infinite classes only contain tournaments on an odd number of vertices. Hence, $P_7^-$ is not isomorphic to any other tournament in the list.
    
    We now compare the remaining tournaments. We proceed by computing some simple statistics on the tournaments to differentiate them. Let $T_x$ be a tournament on $2k+3$ vertices of type $x$. We first compare $T_5$ to the other tournaments.
    We are only considering two tournaments on five vertices, namely $T_5$ and, for $k=1$, $T_c$.
    % There are only two tournaments on five vertices in our consideration, namely $T_5$ and, for $k=1$, $T_c$.
    Since $T_5$ has no vertex of in-degree 1 but $T_c$ does, these tournaments are distinct, and thus $T_5$ is not contained in any of the other classes.

    We leave the comparison with $M_7$ to the end, and now compare the tournaments $T_b, T_c, T_d, T_e$. $T_a$ has two vertices of out-degree 1, namely $w_4$ and $v_{2k-1}$, whereas the tournaments $T_b, T_c, T_d$, and $T_e$ do not. Each of $T_b, T_c, T_d$, and $T_e$ has a unique vertex of maximum out-degree and a unique vertex of minimum out-degree, which we denote $v_{\max}$ and $v_{\min}$, respectively. These are shown in the table below. Let $\tau_x(v)$ denote the number of cyclic triangles containing $v$ in $T_x$, for $T_x \in \{T_a, T_b, T_c, T_d\}$. We compute $\tau_x(v_{\max})$ and $\tau_x(v_{\min})$. It suffices to show these statistics are distinct for $T_b, T_c, T_d$, and $T_e$. We show that these statistics are in fact as in \Cref{tab:cyclic-triangle}.

    \begin{table}[ht]
        \centering
        \begin{tabular}{c|c|c|c|c}
            & $v_{\min}$  & $v_{\max}$ & $\tau(v_{\min})$ & $\tau(v_{\max})$\\
            \hline
            $T_b$ & $w_2$& $w_1$ &$2k-1$ &$2k-1$\\
            $T_c$ & $v_{2k+1}$ & $v_1$&$2$ &$2$\\
            $T_d$ & $w_4$ & $v_1$&$2k$ &$2$\\
            $T_e$ & $v_{2k+1}$ & $w_1$&$2$ &$2k$\\
        \end{tabular}
        \caption{Cyclic triangle statistics of $T_b, T_c, T_d$, and $T_e$.}
        \label{tab:cyclic-triangle}
    \end{table}

    Now, a cyclic triangle can be read off a backedge graph. A triple of vertices $u, v, w$ in a backedge graph appearing in that order forms a cyclic triangle if and only if the set of backedges is one of $\{uw\}$, or $\{uv, vw\}$. For $T_b$, we have $E=\{v_iv_{i+1}\mid 1\leq l<2k+1\}\cup \{w_1v_{2k+1},v_1w_2\}$, and the ordering inducing $G$ is $\sigma=w_1,v_1,\dots,v_{2k+1},w_2$. Thus, all cyclic triangles containing $w_2$ are of the form $v_1, v_i, w_2$, where $i\in \{3,\dots,2k+1\}$, so $\tau_b(v_{\min}) = 2k-1$. By symmetry, $\tau_b(v_{\max}) = 2k-1$. 

    For $T_c$, we have $E=\{v_iv_{i+1}\mid 1\leq l<2k+1\}\cup \{w_1w_2\}$, and the ordering inducing $G$ is $\sigma=v_1,w_1,v_2,\dots,$ $v_{2k},w_2,v_{2k+1}$. The only cyclic triangles containing $v_{1}$ are $\{v_1, v_2, v_3\}$ and $\{v_1, w_1, v_2\}$. Hence, $\tau_c(v_{\max}) = 2$, and by symmetry, $\tau_c(v_{\min}) = 2$. 

    For $T_d$, we have $E=\{v_iv_{i+1}\mid 1\leq l<2k-1\}\cup \{w_1w_4\}$, and the ordering inducing $G$ is $\sigma=v_1,w_1,v_2,\dots,$ $v_{2k},w_2,w_3,v_{2k-1},w_4$. Hence, all cyclic triangles containing $w_4$ are of the form $\{w_1, u, w_4\}$ for $u\in \{w_2,w_3\}\cup \{v_2,\dots, v_{2k-1}\}$. Thus, $\tau_d(v_{\min}) = 2k$. The only cyclic triangles containing $v_{1}$ are $\{v_1, v_2, v_3\}$ and $\{v_1, w_1, v_2\}$. Hence, $\tau_d(v_{\max}) = 2$.

    Since $T_e$ is the mirror of $T_d$, we have $\tau_e(v_{\min}) = 2$ and $\tau_e(v_{\max}) = 2k$. Hence, \Cref{tab:cyclic-triangle} is correct. 
    
    Finally, we consider $M_7$. This tournament has only one vertex of out-degree 1, so it is distinct from $T_a$. If $T_x$ is isomorphic to $M_7$, $T_x$ needs to have seven vertices, so $k=2$. Now, $M_7$ also has a unique vertex of maximum and minimum out-degree, and in $M_7$ $\tau(v_{\min}) = \tau(v_{\max}) = 2$. By \Cref{tab:cyclic-triangle}, the only possible tournament equal to $M_7$ is $T_c$. If $M_7$ and $T_c$ are isomorphic, then such an isomorphism maps $w_1$ in $M_7$ to $v_1$ in $T_c$. Let $u$ be the in-neighbour of $v_1$ in $T_c$. Then, $N_{T_c}^-(u)$ consists of two vertices of in-degree two. However, the same analysis on the in-neighbour of $w_1$ yields two vertices, one of in-degree two, the other of in-degree four. This concludes the proof. 
\end{proof}

\subsection{Unavoidability}
% \Mathieu{Maybe more of a high level overview of 1. The technique, and 2. The structure, of this section could help?}
We finally argue that every minimal nonmatching tournament is in the defined set. With this in mind, fix a minimal nonmatching tournament $T$ arbitrarily.
% \Dinis{Yes, that sounds eminently reasonable.
% Is something like this what you had in mind?}
As mentioned earlier, we will start by showing that we succeed when $T$ is non-$1$-in-degenerate, and then assume that $T$ is $1$-in-degenerate (this second case is by far the one requiring the longest analysis; this is not very surprising, since most tournaments in the set of minimal nonmatching tournaments are $1$-in-degenerate).
In the $1$-in-degenerate case, we will quickly argue that $T$ must have at least three vertices $s,a,b$ and we will know the edges between them in a certain special backedge graph of $T$.
Among other things, we will know that $s$ is first in this backedge graph and consider two cases depending on the relative ordering of $a$ and $b$.
In both of these cases, we will quickly find that $T$ must have two other vertices $c,d$ and will again split into cases depending on the relative ordering of $c$ and $d$.
In each case, we will show that the backedge graph of $T$ we are considering is a backedge graph of one of the elements of our infinite set, or that contradictions to the fact that $T$ is minimal non-matching can be derived (either by finding a matching ordering of $T$ or by finding that the backedge we are considering contains an induced copy of a backedge graph of an element of the minimal forbidden set).
% In none of the cases we will consider is success technically difficult.
% The case-analysis is, however, somewhat cumbersome due to the sheer number of cases that are considered.

Begin by noting that every matching tournament has a vertex with in-degree at most one.
% (this must be true of any vertex that is first in a matching ordering of the tournament).
Hence, if $T$ is non-$1$-in-degenerate, then $T$ is minimal with this property.
Moreover, $T_5$ and $P_7^-$ are explicitly included in our list of minimal nonmatching tournaments, and $\Delta(2,2,2),\Delta(1,2,C_3),Q_6,U_5^{(5)}$ all contain $M_{1,1,adj}^{(3)}$ as a subtournament.

From now on, we may therefore assume that $T$ is $1$-in-degenerate.
In particular, $T$ has a vertex $s$ with in-degree at most one.
By minimality, $T\setminus s$ has a matching ordering.
Consider an ordering $\sigma$ of $V(T)$ such that the ordering of $V(T)\setminus s$ is a matching ordering of the subtournament of $T$ induced by this vertex set and $s$ is first along $\sigma$. Since $\sigma$ is not a matching ordering of $T$, it follows that $s$ has in-degree one in $T$. 
Let $a$ be the neighbour of $s$ in $B_{\sigma}$.
Since $\sigma$ is not a matching ordering of $T$, but induces a matching ordering of $T\setminus s$, it follows that $a$ has degree two in $B_{\sigma}$.
Let $b$ be the neighbour of $a$ in $B_{\sigma}$ that is not $s$.
If $a$ and $b$ were consecutive along $\sigma$, then the ordering of $V(T)$ obtained from $\sigma$ by swapping $a$ and $b$ would be a matching ordering of $T$.
Hence, this is never the case.

\begin{figure}[ht]
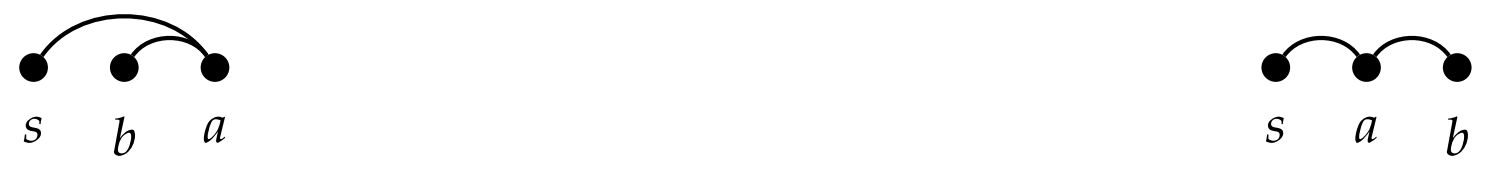

    \centering
    \begin{subfigure}{0.49\linewidth}\centering
        \BackedgeGraphLabeled{2/3, 1/3}{3}{$s$, $b$, $a$}
        \caption{Case 1.}
    \end{subfigure}
    \begin{subfigure}{0.49\linewidth}\centering
        \BackedgeGraphLabeled{1/2, 2/3}{3}{$s$, $a$, $b$}
        \caption{Case 2.}
    \end{subfigure}
    
    \caption{Cases 1 and 2.}
    \label{fig:cases}
\end{figure}
We consider two cases (see Figure \ref{fig:cases}):
\begin{itemize}[align=left,leftmargin=*]
    \item [Case 1.] $b\prec_{\sigma} a$.
        Note that if there exist vertices $c,d$ with $b\prec_{\sigma} c\prec_{\sigma} d\prec_{\sigma} a$, we may assume without loss of generality that they do not neighbour one another in $B_{\sigma}$ (since we could otherwise swap them).
        But then, the ordering $\sigma'=s,c,a,b,d$ induces a backedge graph of $T[s,a,b,c,d]$ that is isomorphic to the backedge graph of $M_{1,1,adj}^{(3)}$ in \Cref{fig:base-tournaments}.
        Hence, if $T$ is minimal nonmatching, we must have $V(T)=\{s,a,b,c,d\}$.
        
        % \begin{figure}[ht]\centering
        %     \BackedgeGraphGeneral{1/5, 2/5}{5}{$s$, $b$, $c$, $d$, $a$}{}
        %     \caption{}\label{subfig:s then b then c then d then a}
        % \end{figure}

        We may therefore assume that there is precisely one vertex, say $c$, between $a$ and $b$ along $\sigma$.
        The ordering $\sigma'$ obtained from $\sigma$ by swapping $c$ and $b$, and then swapping $a$ and $b$
        % , i.e., setting
        % \[\sigma'(\sigma^{-1}(b))=c,\sigma'(\sigma^{-1}(c))=a,\sigma'(\sigma^{-1}(a))=b\]
        % and leaving all other vertices of $T$ in their original places
        is a matching ordering of $T[s,a,b,c]$.
        Since $\sigma'$ cannot be a matching ordering of $T$, $c$ must have a neighbour, say $d$, in $B_{\sigma'}$.
        Since the position of $c$ relative to the vertices in $V(T)\setminus\{s,a,b\}$ is the same in $\sigma$ and $\sigma'$, this vertex must also neighbour $c$ in $B_{\sigma}$.
        We consider two sub-cases (see Figure \ref{fig:subcases1}) depending on its position relative to $c$ in $\sigma$.
        It is worth pointing out that $b,c,a$ are consecutive along $\sigma$,
        % (by assumption)
        and that $c$ has degree at most one in $B_{\sigma}$.
        % (since $B_{\sigma}\setminus s$ has maximum degree one and $c$ does not neighbour $s$ in this graph).
        \begin{figure}[ht]
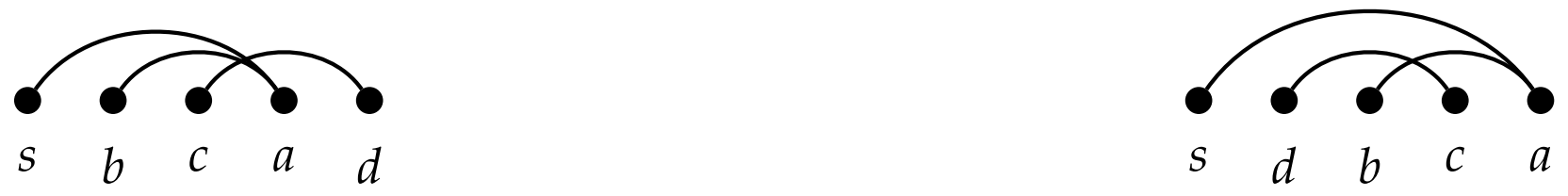

    \centering
    \begin{subfigure}{0.49\linewidth}\centering
        \BackedgeGraphLabeled{1/4, 2/4, 3/5}{5}{$s$, $b$, $c$, $a$, $d$}
        \caption{Sub-case 1.1.}
    \end{subfigure}
    \begin{subfigure}{0.49\linewidth}\centering
        \BackedgeGraphLabeled{1/5, 3/5, 2/4}{5}{$s$, $d$, $b$, $c$, $a$}
        \caption{Sub-case 1.2.}
    \end{subfigure}
    \caption{Sub-cases 1.1 and 1.2.  Note that $b, c, a$ are consecutive in $\sigma$.}
    \label{fig:subcases1}
\end{figure}

        \begin{itemize}[align=left,leftmargin=*]
            \item [Sub-case 1.1.] $c\prec_{\sigma}d$.
                Note that $a\prec_{\sigma}d$.
                Additionally, it can be shown that there is at most one vertex $e$ with $s\prec_{\sigma}e\prec_{\sigma} b$ unless $T$ is isomorphic to $M_{2,2}^{(3)}$.
                Suppose there are at least two such vertices $e_1,e_2$, and assume $e_1e_2\in E(T)$.
                Then the ordering $\sigma''=s,e_1,e_2,a,b,d,c$ induces a backedge graph of the corresponding subtournament of $T$ that is isomorphic to that of $M_{2,2}^{(3)}$ in \Cref{fig:base-tournaments}.
                We may therefore assume that there is at most one vertex $e$ between $s$ and $b$ along $\sigma$.
                \begin{claim} \label{cl:1}
                    If $e$ exists, then either it is isolated in $B_{\sigma}$ or $T$ is a seven-vertex tournament with no matching backedge graph.
                \end{claim}
                \begin{proof}
                  \begin{figure}[ht]
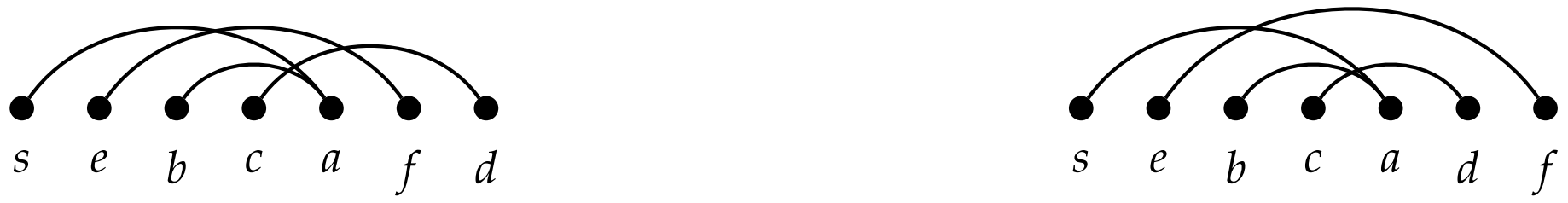

    \centering
    \begin{subfigure}{0.49\linewidth}\centering
        \BackedgeGraphLabeled{1/5, 3/5, 4/7, 2/6}{7}{$s$, $e$, $b$, $c$, $a$, $f$, $d$}
        \caption{$f \prec_\sigma d$.}
    \end{subfigure}
    \begin{subfigure}{0.49\linewidth}\centering
        \BackedgeGraphLabeled{1/5, 3/5, 4/6, 2/7}{7}{$s$, $e$, $b$, $c$, $a$, $d$, $f$}
        \caption{$d \prec_\sigma f$.}
    \end{subfigure}
    \caption{Proof of Claim \ref{cl:1}.}
    \label{fig:claim1}
\end{figure}
                    Assume otherwise, i.e., $e$ exists and has a neighbour $f$ in $B_{\sigma}$
                    Note that $f\notin\{s,a,b,c,d\}$.
                    It follows that $a\prec_{\sigma}f$.
                    We consider two cases (see Figure \ref{fig:claim1}). 
                    \begin{itemize}
                        \item $f\prec_{\sigma} d$.
                            Then the ordering $s,e,c,b,a,f,d$ induces a backedge graph of the corresponding sub-tournament of $T$ that is isomorphic to that of $M_7$ in \Cref{fig:M7}.
                        \item $d\prec_{\sigma} f$.
                            Then the ordering $s,e,a,b,d,c,f$ induces a backedge graph of the corresponding sub-tournament of $T$ that is isomorphic to that of $M_{1,2}^{(3)}$ in \Cref{fig:base-tournaments}.\qedhere
                    \end{itemize}
                \end{proof}
                We may therefore assume that $e$ is isolated in $B_{\sigma}$ (if $e$ exists at all).
                Let $r$ be the number of vertices that are after $a$ along $\sigma$ and let these vertices be $u_1,\dots,u_r$ sorted increasingly along $\sigma$.
                Note that $d=u_i$ for some $1\leq i\leq r$.
                Throughout the rest of the proof, we will repeatedly rely on another ordering $\sigma''$.
                If $e$ exists, then $\sigma''$ is obtained from $\sigma$ by setting its initial five-vertex fragment to $e,a,s,b,c$ and leaving the remaining vertices in their original positions.
                If $e$ does not exist, then $\sigma''$ is obtained from $\sigma$ by setting its initial four-vertex fragment to $a,s,b,c$ and leaving the remaining vertices in their original positions (see Figure \ref{fig:orderings1}). 
                
        \begin{figure}[ht]
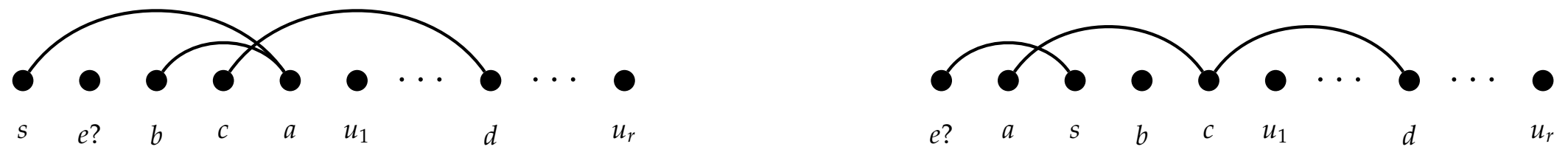

    \centering
    \begin{subfigure}{0.49\linewidth}\centering
        \BackedgeGraphGeneral{1/5, 3/5, 4/7}{8}{$s$, $e?$, $b$, $c$, $a$, $u_1$, $d$,  $u_r$}{7, 6}
        \caption{Ordering $\sigma$.}
    \end{subfigure}
    \begin{subfigure}{0.49\linewidth}\centering
        \BackedgeGraphGeneral{1/3, 2/5, 5/7}{8}{$e?$, $a$, $s$, $b$, $c$, $u_1$, $d$,  $u_r$}{6, 7}
        \caption{Ordering $\sigma''$.}
    \end{subfigure}
    \caption{The orderings $\sigma$ and $\sigma''$.}
    \label{fig:orderings1}
\end{figure}
                \begin{claim}
                    The unique neighbour of $c$ in $B_{\sigma}$ is $u_2$ or $T$ is a seven-vertex tournament with no matching backedge graph.
                \end{claim}
                \begin{proof}
                    Assume, for the sake of contradiction, that the claim is false.
                    Letting $i$ have $u_i=d$, we consider two cases that together encompass all possible values of $i\neq 2$.
                    \begin{itemize}[align=left,leftmargin=*]
                        \item [Option 1:] $i\geq 3$.
                            Consider the ordering $\sigma'''$ of $\{s,a,b,c,u_{i-2},u_{i-1},u_i\}$ obtained from $\sigma''$ by swapping $u_{i-2}$ and $u_{i-1}$ if they are adjacent in $B_{\sigma}$, and that is equal to the ordering of this set induced by $\sigma''$ if they are not.
                            The ordering $\sigma'''$ induces a backedge graph of the corresponding subtournament of $T$ that is isomorphic to that of $M_{2,2}^{(3)}$ in \Cref{fig:base-tournaments}.
                        \item [Option 2:] $i=1$.
                            Note that all neighbours of vertices in $\{u_2,\dots,u_r\}$ in $B_{\sigma}$ are in this set.
                            But then the ordering obtained from $\sigma''$ by swapping $c$ and $u_1$ is a matching ordering of $T$, a contradiction.\qedhere
                    \end{itemize}
                \end{proof}
                We may therefore assume that the only neighbour of $c$ in $B_{\sigma}$ is $u_2$.
                Let $I$ be the largest odd number that is at most $r-4$. For the following claim, see Figure \ref{fig:claim3}. 
             
                \begin{claim}
                    \label{claim: matching proof case 1.1. odd i points at i+3}
                    For each odd $j$ between $1$ and $I-1$, the unique neighbour of $u_j$ in $B_{\sigma}$ is $u_{j+3}$.
                    Moreover, $u_{I}$ is neither isolated nor adjacent to $u_{I+1}$, $u_{I+l}$ with $l\geq 5$, or $u_{I+2}$ in $B_{\sigma}$.
                \end{claim}
                  \begin{figure}[ht]
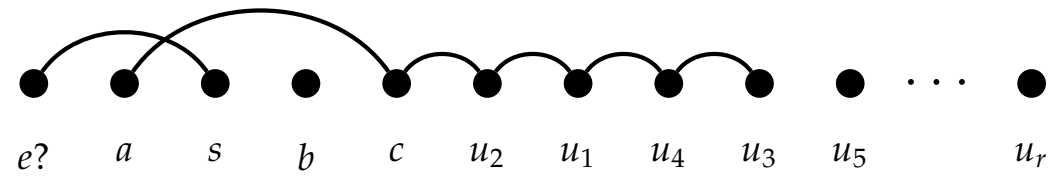

    \centering
        \BackedgeGraphGeneral{1/3, 2/5, 5/6, 6/7, 7/8, 8/9}{11}{$e?$, $a$, $s$, $b$, $c$, $u_2$, $u_1$, $u_4$, $u_3$, $u_5$, $u_r$}{10}
        \caption{Ordering $\sigma_3''$.}

    \label{fig:claim3}
\end{figure}
                \begin{proof}
                    For every odd $k$ between $1$ and $r-1$, we inductively define $\sigma''_k$ to be the ordering obtained from $\sigma''_{k-2}$ (with $\sigma''_{-1}=\sigma''$) by swapping $u_k$ and $u_{k+1}$.
    
                    Assume, for the sake of contradiction, that the claim is false, and let $j$ be the smallest odd number falsifying it.
                    Since $j$ is smallest, the ordering $\sigma''_{j-2}$ induces a backedge graph of $T[\{c,u_1,\dots,u_{j-1}\}]$ that is isomorphic to a path on $j$ vertices.
                    The following also turns out to be useful.
                    \begin{observation}
                        \label{obs: vertex sets closed under nbhd case 1.1}
                        All neighbours of vertices in $S\coloneqq\{s,a,b,c,e,u_1,\dots,u_{j+1}\}\cup N_{B_{\sigma}}(u_j)$ in the backedge graphs $B_{\sigma}$, $B_{\sigma''}$, and $B_{\sigma''_{j-2}}$ are in $S$ and an analogous statement applies to $V(T)\setminus S$.
    
                        In particular, if $\sigma'''$ is an ordering of $T$ that induces a matching ordering of $T[S]$ and a matching ordering of $T\setminus S$, and the relative positions of vertices in $S$ and $V(T)\setminus S$ are the same in $\sigma'''$ as in $\sigma$, then $\sigma'''$ is a matching ordering of $T$.
                    \end{observation}
                    We consider three possibilities for the neighbourhood of $u_j$ in $B_{\sigma}$:
                    \begin{itemize}[align=left,leftmargin=*]
                        \item[Option 1:\label{case 1.1 opt 1}] $u_j$ is isolated in $B_{\sigma}$.
                            Note that \Cref{lemma: even path with two isolated vertices on one end} or \Cref{lemma: even path with two isolated vertices on one end plus an incoming edge} applies to $T[S]^*$ (the precise lemma that applies depends on the existence of $e$, but the mirror of $\sigma''_{j}$ always certifies that the relevant lemma applies).
                            In particular, $T[S]$ has a matching ordering $\sigma'''$.
                            By \Cref{obs: vertex sets closed under nbhd case 1.1}, the ordering obtained by appending the final $r-j-1$-vertex fragment of $\sigma$ to $\sigma'''$ is a matching ordering of $T$, a contradiction.
                            
                        \item[Option 2:] $u_j$ neighbours $u_l$ in $B_{\sigma}$ with $l\geq j+4$ (or $l\geq j+5$ if $j=I$).
                            Then the ordering of $S\cup\{u_{l-1},u_{l-2}\}$ induced by $\sigma''_{j}$ induces a backedge graph of the corresponding subtournament of $T$ that is isomorphic to that of $M_{2,2}^{(j+4)}$ in \Cref{def: infinite classes minimal nonmatching} (up to swapping $u_{l-2}$ and $u_{l-1}$ if they are adjacent in $B_{\sigma}$).
                            Since this is a proper subtournament of $T$, as evidenced by the fact that its vertex set does not contain $u_{l-3}$ (when $l\geq j+5$) or $u_{l+1}$ (which is always applicable when $l=j+4$, because in this case $j < I \leq r-4$) that has no matching backedge graph, we obtain a contradiction to minimality of $T$.
                            
                        \item[Option 3:] $u_j$ neighbours $u_{j+2}$ in $B_{\sigma}$.
                            Once more, $\sigma''_j$ certifies that \Cref{lemma: odd path with two isolated vertices on one end} or \Cref{lemma: odd path with two isolated vertices on one end plus an incoming edge} applies to $T[S]$ and we can find a matching ordering of $T[S]$.
                            Precisely as in Option 1, \Cref{obs: vertex sets closed under nbhd case 1.1} enables us to find a matching ordering of $T$, a contradiction.\qedhere
                    \end{itemize}
                \end{proof}
    
                Finally, we claim that $u_I$ neighbours $u_{I+4}$ in $B_{\sigma}$.
                In particular, this implies that the ordering $\sigma''_{I}$ induces a backedge graph of $T[s,a,b,c,u_1,\dots,u_{I+4}]$ that is isomorphic to that of $M_{2,2}^{(I+4)}$ in \Cref{def: infinite classes minimal nonmatching}.
                By minimality, $T$ must be isomorphic to that tournament.

                \begin{claim}
                    The unique neighbour of $u_I$ in $B_{\sigma}$ is $u_{I+4}$.
                \end{claim}
                \begin{proof}
                \begin{figure}[ht]
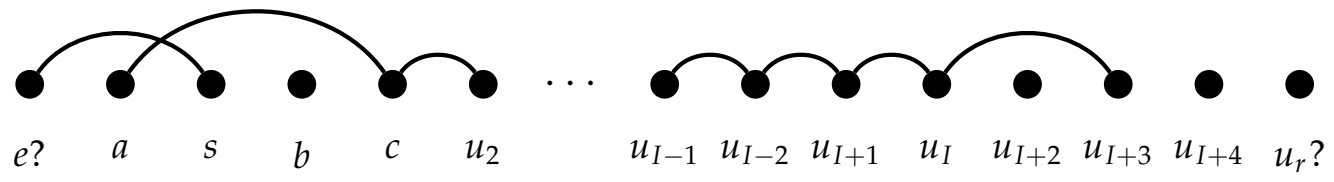

    \footnotesize
    \centering
        \BackedgeGraphGeneral{1/3, 2/5, 5/6, 7/8, 8/9, 9/10, 10/12}{14}{$e?$, $a$, $s$, $b$, $c$, $u_2$, $u_{I-1}$, $u_{I-2}$, $u_{I+1}$, $u_I$, $u_{I+2}$, $u_{I+3}$, $u_{I+4}$, $u_r?$}{6}
        \caption{Ordering $\sigma_I''$, assuming $u_I$ is adjacent to $u_{I+3}$ in $B_\sigma$.}

    \label{fig:claim4}
\end{figure}
                    By \Cref{claim: matching proof case 1.1. odd i points at i+3}, the only other option is for $u_{I}u_{I+3}$ to be an edge in $B_{\sigma}$. See Figure \ref{fig:claim4}. 
                    We consider three different cases depending on the neighbourhood of $u_{I+2}$ in $B_{\sigma}$, noting that it is either empty, contains only $u_{I+4}$, or contains only $u_{I+5}$.
                    \begin{enumerate}
                        \item $u_{I+2}$ is isolated in $B_{\sigma}$.
                        Then \Cref{lemma: even path with two isolated vertices on one end} or \Cref{lemma: even path with two isolated vertices on one end plus an incoming edge} applies to $T[s,a,b,c,e,u_1,\dots,u_{I+3}]$ (as evidenced by $\sigma''_{I+2}$) and, much as in Option 1 of the proof of \Cref{claim: matching proof case 1.1. odd i points at i+3}, the corresponding matching ordering can be combined with $\sigma$ to obtain a matching backedge graph of $T$.
    
                        \item $u_{I+2}$ neighbours $u_{I+4}$.
                        Then \Cref{lemma: odd path with two isolated vertices on one end} or \Cref{lemma: odd path with two isolated vertices on one end plus an incoming edge} applies to $T[s,a,b,c,e,u_1,\dots,u_{I+4}]$ (as evidenced by $\sigma''_{I+2}$).
                        This either gives a matching ordering of $T$ (if $r=I+4$), or an ordering to which $u_{I+5}$ can be appended to obtain a matching ordering of $T$.
    
                        \item $u_{I+2}$ neighbours $u_{I+5}$.
                        Then $\sigma_{I+4}''$ shows that \Cref{lemma: even path with two isolated vertices on one end} or \Cref{lemma: even path with two isolated vertices on one end plus an incoming edge} applies to $T$ (as evidenced by $\sigma''_{I+4}$).
                        This implies the existence of a matching ordering of $T$.
                    \end{enumerate}
    
                    In each of the considered possibilities, we managed to find a matching ordering of $T$.
                    Such an ordering does not exist by assumption.
                    Hence, $u_I$ cannot neighbour $u_{I+3}$.
                \end{proof}
                This concludes our examination of Sub-case 1.1.
            
            \item [Sub-case 1.2.] $d\prec_{\sigma}c$.
                Let $l$ be the number of vertices that are between $s$ and $b$ along $\sigma$ and let these vertices be $u_1,\dots,u_{l}$ sorted decreasingly along $\sigma$ (i.e., $u_l$ is the closest vertex to $s$ and $u_1$ the closest vertex to $b$).
                % Note that $d=u_i$ for some $1\leq i\leq l$.
                Throughout the analysis of this case, we will make repeated use of the ordering $\sigma'$.
                Recall that this is the ordering obtained from $\sigma$ by swapping $b,c$ and then $b,a$. See Figure \ref{fig:subcases12}.
                 \begin{figure}[ht]
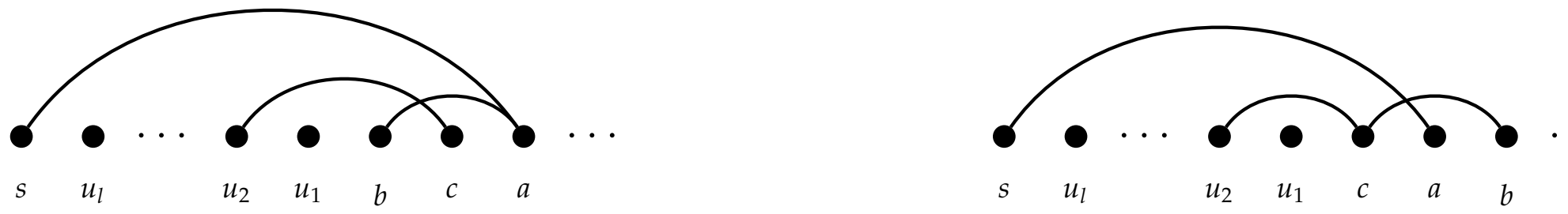

    \centering
    \begin{subfigure}{0.49\linewidth}\centering
        \BackedgeGraphGeneral{1/7, 5/7, 3/6}{7}{$s$, $u_l$, $u_2$, $u_1$, $b$, $c$, $a$}{2, 7}
        \caption{Sub-case 1.2, ordering $\sigma$.}
    \end{subfigure}
    \begin{subfigure}{0.49\linewidth}\centering
       \BackedgeGraphGeneral{1/6, 5/7, 3/5}{7}{$s$, $u_l$, $u_2$, $u_1$,  $c$, $a$, $b$}{2, 7}
        \caption{Sub-case 1.2, ordering $\sigma'$.}
    \end{subfigure}
    \caption{Sub-case 1.2, assuming Claim \ref{cl:5}. Specifically, in Sub-case 1.2, we know that $c$ is adjacent to one of $u_l, \dots, u_1$.}
    \label{fig:subcases12}
\end{figure}
                \begin{claim}
                    The unique neighbour of $c$ in $B_{\sigma}$ is $u_2$ or $T$ is a seven-vertex tournament with no matching ordering.\label{cl:5}
                \end{claim}
                \begin{proof}
                    Assume, for the sake of contradiction, that the claim is false.
                    Letting $i$ have $u_i=d$, we consider two cases that together encompass all possible values of $i\neq 2$.
                    \begin{itemize}[align=left,leftmargin=*]
                        \item [Option 1:] $i\geq 3$.
                            Consider the ordering $\sigma''$ of $\{s,a,b,c,u_{i-2},u_{i-1},u_i\}$ that can be obtained from $\sigma'$ by swapping $u_{i-1}$ and $u_{i-2}$ if they are adjacent in $B_{\sigma}$, and that is equal to the ordering of this set induced by $\sigma'$ if they are not.
                            The ordering $\sigma''$ induces a backedge graph of the corresponding subtournament of $T$ that is isomorphic to that of $M_{2,1}^{(3)}$ in \Cref{fig:base-tournaments}.
                            
                        \item [Option 2:] $i=1$.
                            Note that the ordering of $\{s,a,b,c,u_1\}$ that can be obtained from $\sigma'$ by swapping $a$ and $b$ induces a backedge graph of $T[\{s,a,b,c,u_1\}]$ that is the mirror of that in \Cref{lemma: even path plus a long forward edge}.
                            By \Cref{lemma: even path plus a long forward edge}, $T[\{s,a,b,c,u_1\}]$ has a matching ordering $\sigma''$ where $s$ is first.
                            It follows that the ordering
                            \[\sigma''(1)=s,\sigma(2),\dots,\sigma(r-1),\sigma''(2),\sigma''(3),\sigma''(4),\sigma''(5),\sigma(r+4),\dots,\sigma(n(T))\]
                            is a matching ordering of $T$ (since the relative positions of vertices in $\{s,a,b,c,u_1\}$ are the same in $\sigma$ and this ordering, $\sigma$ is a matching ordering of $V(T)\setminus\{s,a,b,c,u_1\}$, and $\sigma''$ is a matching ordering of $\{s,a,b,c,u_1\}$).
                            This is a contradiction to the fact that $T$ is nonmatching.\qedhere
                    \end{itemize}
                \end{proof}
                We may therefore assume that the only neighbour of $c$ in $B_{\sigma}$ is $u_2$.
                Let $I$ be the largest odd number that is at most $l-4$.
                
                \begin{claim}
                    \label{claim: matching proof case 1.2. odd i points at i+3}
                    For each odd $j$ between $1$ and $I-1$ the unique neighbour of $u_j$ in $B_{\sigma}$ is $u_{j+3}$.
                    Moreover, $u_{I}$ is neither isolated nor adjacent to $u_{I+1}$, $u_{I+l}$ with $l\geq 5$, $u_{I+2}$, or $e$ in $B_{\sigma}$. \label{cl:6}
                \end{claim}
                 \begin{figure}[ht]
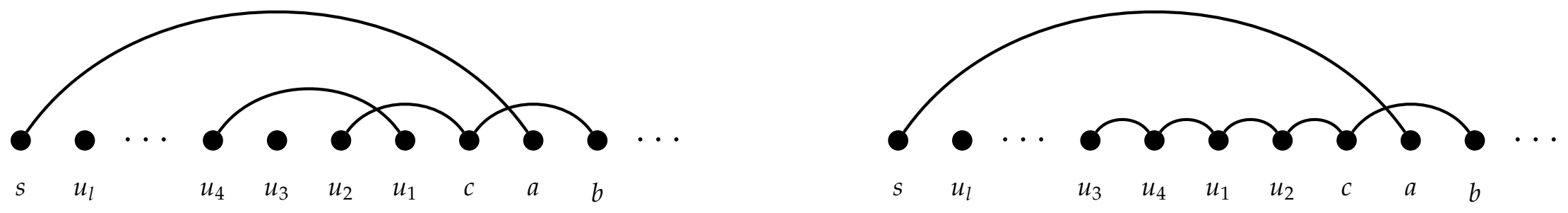

    \centering
    \begin{subfigure}{0.49\linewidth}\centering
        \BackedgeGraphGeneral{1/8, 7/9, 5/7, 3/6}{9}{$s$, $u_l$, $u_4$, $u_3$, $u_2$, $u_1$,  $c$, $a$, $b$}{2, 9}
        \caption{Claim \ref{cl:6}, ordering $\sigma'$.}
    \end{subfigure}
    \begin{subfigure}{0.49\linewidth}\centering
        \BackedgeGraphGeneral{1/8, 7/9, 6/7, 4/5, 3/4, 5/6}{9}{$s$, $u_l$, $u_3$, $u_4$, $u_1$, $u_2$,  $c$, $a$, $b$}{2, 9}
        \caption{Claim \ref{cl:6}, ordering $\sigma'_3$.}
    \end{subfigure}
    \caption{The proof of Claim \ref{cl:6}.}
    \label{fig:claim6}
\end{figure}
                \begin{proof}
                    For every odd $k$ between $1$ and $l$, we inductively define $\sigma'_{k}$ to be the ordering obtained from $\sigma'_{k-2}$ (with $\sigma'_{-1}=\sigma'$) by swapping $u_{k}$ and $u_{k+1}$. See Figure \ref{fig:claim6}.

                    Assume, for the sake of contradiction, that the claim is false, and let $j$ be the smallest odd number falsifying it.
                    Since $j$ is smallest, the ordering $\sigma'_{j-2}$ induces a backedge graph of $T[c,u_1,\dots,u_{j-1}]$ that is isomorphic to a path on $j$ vertices.
                    The following also turns out to be useful.
                    \begin{observation}
                        \label{obs: vertex sets closed under nbhd case 1.2}
                        All neighbours of vertices in $S\coloneqq\{s,a,b,c,e,u_1,\dots,u_{j+1}\}\cup N_{B_{\sigma}}(u_j)$ in the backedge graphs $B_{\sigma}$, $B_{\sigma''}$, and $B_{\sigma''_{j-2}}$ are in $S$ and an analogous statement applies to $V(T)\setminus S$.
    
                        In particular, if $\sigma'''$ is an ordering of $T$ that induces a matching ordering of $T[S]$ and a matching ordering of $T\setminus S$, and the relative positions of vertices in $S$ and $V(T)\setminus S$ are the same in $\sigma'''$ as in $\sigma$, then $\sigma'''$ is a matching ordering of $T$.
                    \end{observation}
                    We consider four possibilities for the neighbourhood of $u_j$ in $B_{\sigma}$:
                    \begin{itemize}[align=left,leftmargin=*]
                        \item [Option 1:] $u_j$ is isolated in $B_{\sigma}$.
                            Note that \Cref{lemma: odd path plus a long forward edge} applies to the mirror of $T[S]$ (the ordering obtained from $\sigma'_{j}$ by swapping $a$ and $b$ shows that this is the case).
                            In particular, this subtournament of $T$ has a matching ordering $\sigma''$ where $s$ is first.
                            By \Cref{obs: vertex sets closed under nbhd case 1.2} the ordering of $T$ that induces the same ordering of $S$ as $\sigma''$ and where the vertices in $V(T)\setminus S$ occupy the same place as in $\sigma$ is a matching ordering, a contradiction.
                            
                        \item [Option 2:] $u_j$ neighbours $u_r$ in $B_{\sigma}$ with $r\geq j+4$ (or $r\geq j+5$ if $j=I$).
                            Then the ordering of $S\cup\{u_{r-1},u_{r-2}\}$ induced by $\sigma'_{j}$ induces a backedge graph of $T[S\cup\{u_{r-1},u_{r-2}\}]$ that is isomorphic to that of $M_{2,1}^{(j+4)}$ in \Cref{def: infinite classes minimal nonmatching} (up to swapping $u_{r-2}$ and $u_{r-1}$ if they are adjacent in $B_{\sigma}$).
                            Since this is a proper subtournament of $T$, as evidenced by the fact that its vertex set does not contain $u_{l-3}$ (when $l\geq j+5$) or $u_{l+1}$ (when $l=j+4$), that has no matching backedge graph, we obtain a contradiction to minimality of $T$.
                            
                        \item [Option 3:] $u_j$ neighbours $u_{j+2}$ in $B_{\sigma}$.
                            By \Cref{lemma: even path plus a long forward edge}, we can find a matching ordering of $T[S]$ (this is certified by the ordering obtained from $\sigma'_{j}$ by swapping $a$ and $b$).
                            Precisely as in Option 1, this enables us to find a matching ordering of $T$, a contradiction.
                        \item [Option 4:] $u_j$ neighbours $e$ with $a\prec_{\sigma} e$ in $B_{\sigma}$.
                            Then the ordering obtained from $\sigma'_{j}$ by swapping $a$ and $b$ induces a backedge graph of $T[S]$ that is isomorphic to that of $M_{1,1,is}^{(j+4)}$ in \Cref{def: infinite classes minimal nonmatching}.
                            Since this is a proper subtournament of $T$, as evidenced by the fact that its vertex set does not contain $u_{j+2}$, that has no matching backedge graph, we obtain a contradiction to minimality of $T$.
                    \end{itemize}
                    Since all possibilities lead to contradictions, no such $j$ can exist.
                \end{proof}

                Finally, we investigate the neighbourhood of vertices $u_j$ for odd $j\geq I$.
                \begin{claim}
                    If $l$ is odd, then the unique neighbour of $u_{I}$ in $B_{\sigma}$ is $u_{I+4}$.
                    If $l$ is even, then $u_{I}$ neighbours $u_{I+3}$, $u_{I+2}$ neighbours $u_{I+5}$, and $u_{I+4}$ neighbours $e$ with $a\prec_{\sigma} e$ in $B_{\sigma}$. \label{cl:7}
                \end{claim}
 \begin{figure}[ht]
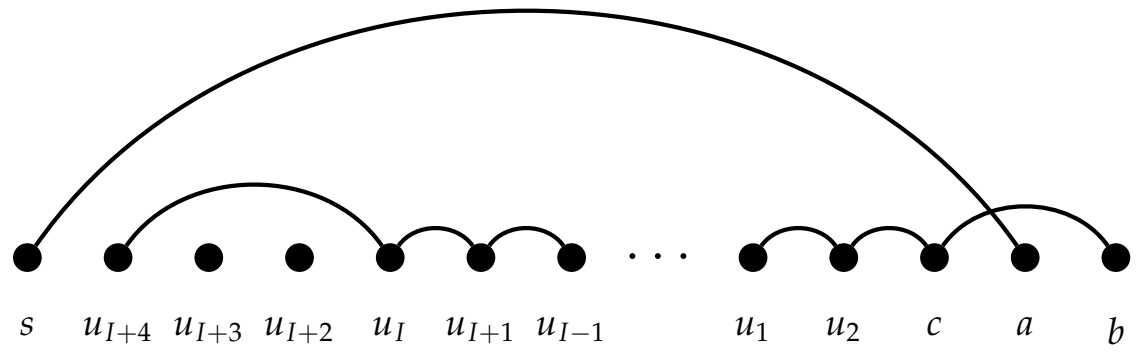
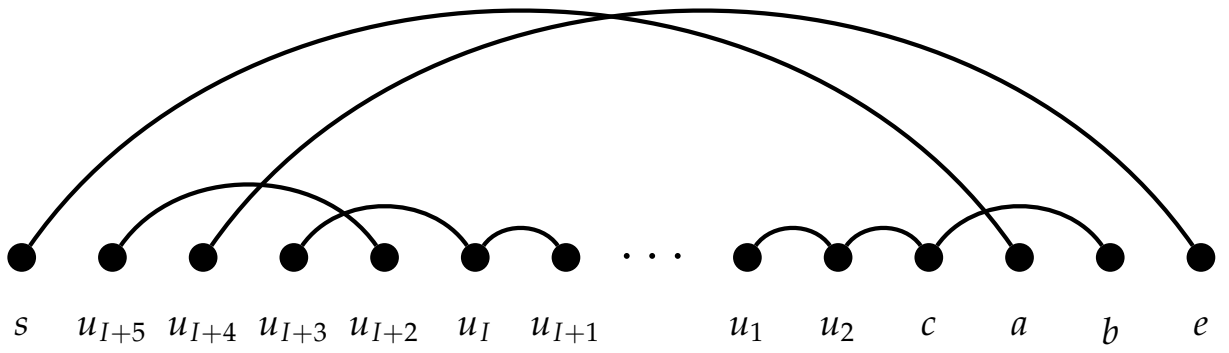

    \centering
    \begin{subfigure}{0.5\linewidth}\centering
        \BackedgeGraphGeneral{1/11, 10/12, 8/9, 9/10, 5/6, 6/7, 2/5}{12}{$s$, $u_{I+4}$, $u_{I+3}$, $u_{I+2}$, $u_{I}$, $u_{I+1}$, $u_{I-1}$,  $u_1$, $u_2$,  $c$, $a$, $b$}{7}
        \caption{Claim \ref{cl:7} when $l$ is odd.}
    \end{subfigure}
    \begin{subfigure}{0.5\linewidth}\centering
        \BackedgeGraphGeneral{1/11, 10/12, 8/9, 9/10, 6/7, 4/6, 2/5, 3/13}{13}{$s$, $u_{I+5}$, $u_{I+4}$, $u_{I+3}$, $u_{I+2}$, $u_{I}$, $u_{I+1}$,  $u_1$, $u_2$,  $c$, $a$, $b$, $e$}{7}
        \caption{Claim \ref{cl:7} when $l$ is even.}
    \end{subfigure}
    \caption{The outcomes of Claim \ref{cl:7}, using ordering $\sigma_I'$.}
    \label{fig:claim7}
\end{figure}
                
                Note that this implies that $T$ is $M_{2,1}^{(I+4)}$ when $l$ is odd (as certified by $\sigma'_{I}$, or the ordering obtained from it by swapping $u_{I+2}$ and $u_{I+3}$ if they are adjacent in $B_{\sigma}$).
                Similarly, when $l$ is even, it follows that $T$ is $M_{1,1,is}^{I+8}$ (as certified by the ordering obtained from $\sigma_{I+4}'$ by swapping $a$ and $b$). See Figure \ref{fig:claim7}. 
                \begin{proof}
                    We handle each parity of $l$ separately:
                    \begin{itemize}[align=left,leftmargin=*]
                        \item [Option 1:] $l$ is odd.
                            Note that this implies that $l=I+4$.
                            The only other possible neighbour of $u_{I}$ is $u_{I+3}$.
                            By \Cref{lemma: odd path plus a long forward edge}, we can find a matching ordering of $T[\{s,a,b,c,u_1,\dots,u_{I+3}\}]$ (this is certified by the ordering obtained from $\sigma'_{I+2}$ by swapping $a$ and $b$).
                            As in Option 1 of the proof of \Cref{claim: matching proof case 1.2. odd i points at i+3}, this enables us to find a matching ordering of $T$, a contradiction.
                            
                        \item [Option 2:] $l$ is even.
                            Note that this implies that $l=I+5$.
                            
                            Begin by noting that the proof of \Cref{claim: matching proof case 1.2. odd i points at i+3} still goes through for $u_I$ (since all it really requires is for $u_{I+5}$ to be a vertex of $T$).
                            This implies that $u_{I}$ neighbours $u_{I+3}$ $B_{\sigma}$.

                            Similarly, the proof of \Cref{claim: matching proof case 1.2. odd i points at i+3} still goes through for $u_{I+2}$ (as it is impossible for $u_{I+2}$ to neighbour $u_{I+2+4}$ in $B_{\sigma}$, i.e., the only potentially problematic case, Option 2 in that proof, never occurs).
                            
                            Finally, we show that $u_{I+4}$ neighbours $e$ with $a\prec_{\sigma}e$ in $B_{\sigma}$.
                            The only other option (since $B_{\sigma}\setminus s$ has maximum degree one) is for $u_{I+4}$ to be isolated.
                            But then, much like in Option 3 of the proof of \Cref{claim: matching proof case 1.2. odd i points at i+3}, \Cref{lemma: even path plus a long forward edge} applies to the mirror of $T[\{s,a,b,c,u_1,\dots,u_{l}\}]$ (the ordering obtained from $\sigma'_{I+4}$ by swapping $a$ and $b$ shows that this is the case).
                            In particular, this subtournament of $T$ has a matching ordering $\sigma''$.
                            Prepending $\sigma''$ to the the ordering of $V(T)\setminus\{s,a,b,c,u_1,\dots,u_{l}\}$ induced by $\sigma$ yields a matching ordering of $T$, a contradiction.\qedhere
                    \end{itemize}
                \end{proof}
        \end{itemize}
        This concludes our examination of Sub-case 1.2., and therefore of Case 1.

\begin{figure}[ht]
    \centering
    
    \begin{subfigure}{0.35\linewidth}\centering
        \BackedgeGraphGeneral{1/3, 3/6}{6}{$s$, $c$, $a$, $d_1$, $d_2$, $b$}{5,6}
        \caption{Sub-case 2.1.}
    \end{subfigure}
    \begin{subfigure}{0.35\linewidth}\centering
        \BackedgeGraphGeneral{1/4, 4/6}{6}{$s$, $c_1$, $c_2$, $a$, $d$, $b$}{3,6}
        \caption{Sub-case 2.2.}
    \end{subfigure}
    \begin{subfigure}{0.25\linewidth}\centering
        \BackedgeGraphGeneral{1/3, 3/5}{5}{$s$, $c$, $a$, $d$, $b$}{5}
        \caption{Sub-case 2.3.}
    \end{subfigure}
    
        \caption{Case 2.}
    \label{fig:case2}
\end{figure}

    \item [Case 2.] $a\prec_{\sigma} b$.
        Note that neither $s,a$, nor $a,b$ are consecutive along $\sigma$, since the ordering obtained from $\sigma$ by swapping the vertices in one consecutive pair would otherwise yield a matching ordering of $T$.
        Hence, there exist $c$ with $s\prec_{\sigma}c\prec_{\sigma}a$ and $d$ with $a\prec_{\sigma}d\prec_{\sigma} b$.
        Moreover, we may assume that no such $c,d$ are adjacent in $B_{\sigma}$, since the backedge graph of $T[s,a,b,c,d]$ is otherwise isomorphic to that of $M_{1,1,adj}^{(3)}$ in \Cref{fig:base-tournaments}.
        Additionally, if there are two vertices $c_1,c_2$ between $s$ and $a$ along $\sigma$, and two vertices $d_1,d_2$ between $a$ and $b$ along $\sigma$, we assume without loss of generality that $c_1\Rightarrow_{T}c_2$ and $d_1\Rightarrow_{T}d_2$.
        But then the ordering $s,c_1,c_2,a,d_1,d_2,b$ induces a backedge graph of $T[s,a,b,c_1,c_2,d_1,d_2]$ that is isomorphic to that of $M_{2,2}^{(3)}$ in \Cref{fig:base-tournaments}.
        By minimality, $T$ would have to be isomorphic to that tournament.
        Hence, we may assume that $c$ is the only vertex between $s$ and $a$, or $d$ is the only vertex between $a$ and $b$ along $\sigma$ (or both).
        We now consider three different sub-cases depending on the number of vertices between $s,a$ and $a,b$ along $\sigma$. See Figure \ref{fig:case2}. 
        
        \begin{itemize}[align=left,leftmargin=*]
            \item [Sub-case 2.1.] There is precisely one vertex $c$ between $s$ and $a$ along $\sigma$, and there are two vertices $d_1,d_2$ with $d_1d_2\in E(T)$ between $a$ and $b$ along $\sigma$.
                Note that $c$ must have a neighbour in $B_{\sigma}$, since the ordering obtained from $\sigma$ by swapping $s$ and $c$, and then swapping $s$ and $a$ is otherwise a matching ordering of $T$.
                Moreover, this vertex, $e$, must come after $b$ in $\sigma$.
                But then, the ordering $s,c,a,d_1,d_2,b,e$ induces a backedge graph of $T[s,a,b,c,d_1,d_2,e]$ that is isomorphic to that of $M_{1,2}^{(3)}$ in \Cref{fig:base-tournaments}.
                Minimality of $T$ now implies that it is isomorphic to this tournament and we are done.
            \item [Sub-case 2.2.] There are at least two vertices $c_1,c_2$ with $c_1\Rightarrow c_2$ between $s$ and $a$ along $\sigma$, and there is precisely one vertex $d$ between $a$ and $b$ along $\sigma$.
                Note that $d$ must have a neighbour in $B_{\sigma}$, since the ordering obtained from $\sigma$ by swapping $d$ and $b$, and then $a$ and $b$ is otherwise a matching ordering of $T$.
                Moreover, this vertex must come after $b$ in $\sigma$.
                Throughout the analysis of Sub-case 2.2, we will use the ordering $\sigma'$ which is obtained from $\sigma$ by swapping $d$ and $b$ if $c_1\prec_{\sigma} c_2$, and obtained from $\sigma$ by swapping $d$ and $b$, and swapping $c_1$ and $c_2$ otherwise.
                Let $r$ be the number of vertices that are after $b$ along $\sigma$ and let these vertices be $u_1,\dots,u_r$ sorted increasingly along $\sigma$. See Figure \ref{fig:claim8}.
                \begin{figure}[ht]
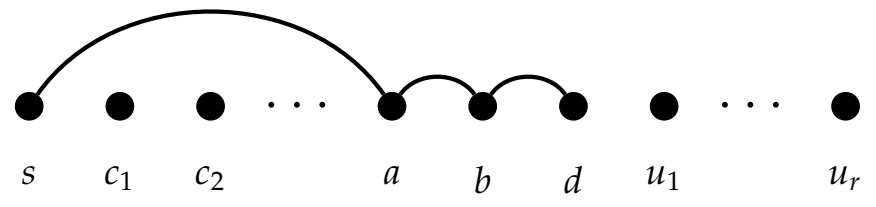
 \centering
                   \BackedgeGraphGeneral{1/4, 4/5, 5/6}{8}{$s$, $c_1$, $c_2$, $a$, $b$, $d$, $u_1$, $u_r$}{3,7}
                   \caption{The ordering $\sigma'$; neighbour of $d$ not shown.}
                   \label{fig:claim8}
                   \end{figure}
                \begin{claim}
                    The unique neighbour of $d$ in $B_{\sigma}$ is $u_2$ or $T$ is a nine-vertex tournament in the forbidden set.
                \end{claim}
                \begin{proof}
                    By the argumentation from above, $d$ neighbours $u_i$ for some $1\leq i\leq r$.
                    We consider two cases that together encompass all possible values of $i\neq 2$.
                    \begin{itemize}[align=left,leftmargin=*]
                        \item [Option 1:] $i\geq 3$.
                            Consider the ordering $\sigma''$ of $\{s,a,b,c_1,c_2,d,u_{i-2},u_{i-1},u_i\}$ obtained from $\sigma'$ by swapping $u_{i-2}$ and $u_{i-1}$ if they are adjacent in $B_{\sigma}$, and that is equal to $\sigma'$ if they are not.
                            The ordering $\sigma''$ induces a backedge graph of $T[\{s,a,b,c_1,c_2,d,u_{i-2},u_{i-1},u_{i}\}]$ isomorphic to that of $M_{2,2}^{(5)}$ in \Cref{def: infinite classes minimal nonmatching}.
                        \item [Option 2:] $i=1$.
                            Then the ordering $\sigma''$ obtained from $\sigma'$ by swapping $a$ and $b$,  and swapping $d$ and $u_1$, is a matching ordering of $T$, a contradiction.\qedhere
                    \end{itemize}
                \end{proof}
                We may therefore assume that the only neighbour of $d$ in $B_{\sigma}$ is $u_2$.
                Let $I$ be the largest odd number that is at most $r-4$.
                \begin{claim}
                    \label{claim: matching proof case 2.2. odd i points at i+3}
                    For each odd $j$ between $1$ and $I-1$, the unique neighbour of $u_j$ in $B_{\sigma}$ is $u_{j+3}$.
                    Moreover, $u_{I}$ is neither isolated nor adjacent to $u_{I+1}$, $u_{I+l}$ with $l\geq 5$, $u_{I+2}$, or $c$ in $B_{\sigma}$.
                \end{claim}

                \begin{proof}
                    For every odd $k$ between $1$ and $r$, we inductively define $\sigma'_k$ to be the ordering obtained from $\sigma'_{k-2}$ (with $\sigma'_{-1}=\sigma'$) by swapping $u_k$ and $u_{k+1}$.

                    Assume, for the sake of contradiction, that the claim is false, and let $j$ be the smallest odd number falsifying it.
                    Since $j$ is smallest, the ordering $\sigma'_{j-2}$ induces a backedge graph of $T[a,b,d,u_1,\dots,u_{j-1}]$ that is isomorphic to a path on $j+2$ vertices.
                    The following also turns out to be useful.
                    \begin{observation}
                        \label{obs: vertex sets closed under nbhd case 2.2}
                        All neighbours of vertices in $S=\{s, a, b, d,u_1,\dots,u_{j+1}\}\cup N_{B_{\sigma}}(u_j)$ in the backedge graphs $B_{\sigma}$, $B_{\sigma'}$, and $B_{\sigma_{j-2}'}$ are in $S$ and an analogous statement applies to $V(T)\setminus S$.
                    \end{observation}
                    Moreover, all neighbours of vertices in $\{d,u_1,\dots,u_{j-3},u_{j-1}\}$ in $B_{\sigma}$ are in $\{d,u_1,\dots,u_{j-1}\}$.
                    We consider four possibilities for the neighbourhood of $u_j$ in $B_{\sigma}$:
                    \begin{itemize}[align=left,leftmargin=*]
                        \item [Option 1:] $u_j$ is isolated in $B_{\sigma}$.
                            Note that $T[S \setminus \{s\}]$ has a backedge graph, $B_{\sigma'_j}$, that is an induced monotone path on $j+4$ vertices.
                            By \Cref{lemma: matching orderings of induced paths}, there is an ordering $\sigma''$ of $S \setminus \{s\}$ that induces a matching backedge graph of this subtournament of $T$ where $a$ has degree zero.
                            Now the ordering of $T$ that induces the same ordering of $S \setminus \{s\}$ as $\sigma''$ and where the vertices in $V(T)\setminus(S \setminus \{s\})$ occupy the same place as in $\sigma$ is a matching ordering of $T$, a contradiction.

                        \item [Option 2:] $u_j$ neighbours $u_l$ with $l\geq j+4$ in $B_{\sigma}$ (or $l\geq j+5$ if $j=I$).
                            Then the ordering of $S\cup\{c_1, c_2,u_{l-2},u_{l-1}\}$ induced by $\sigma'_{j}$ induces a backedge graph of $T[S\cup\{c_1, c_2, u_{l-2},u_{l-1}\}]$ that is isomorphic to that of $M_{2,2}^{(j+6)}$ in \Cref{def: infinite classes minimal nonmatching} (up to swapping $u_{l-2}$ and $u_{l-1}$ if they are adjacent in $B_{\sigma}$).
                            Since this is a proper subtournament of $T$, as evidenced by the fact that its vertex set does not contain $u_{l-3}$ (if $l\geq j+5$) or $u_{l+1}$ (if $l=j+4$), that has no matching backedge graph, we obtain a contradiction to minimality of $T$.

                        \item [Option 3:] $u_j$ neighbours $u_{j+2}$ in $B_{\sigma}$.
                            Note that $T[\{a,b,d,u_1,\dots,u_{j+2}\}]$ has a backedge graph, $B_{\sigma'_j}$, that is an induced monotone path on $j+5$ vertices.
                            By \Cref{lemma: matching orderings of induced paths}, there is an ordering $\sigma''$ of $S \setminus \{s\}$ that induces a matching backedge graph of this subtournament of $T$ where $a$ has degree zero.
                            By \Cref{obs: vertex sets closed under nbhd case 2.2} the ordering of $T$ that induces the same ordering of $S \setminus \{s\}$ as $\sigma''$ and where the vertices in $V(T)\setminus(S\setminus\{s\})$ occupy the same place as in $\sigma$ is a matching ordering of $T$, a contradiction.
                            
                        \item [Option 4:] $u_j$ neighbours $c$ with $s\prec_{\sigma} c\prec_{\sigma} a$ in $B_{\sigma}$.
                            Then the ordering of $S\cup\{c\}$ induced by $\sigma'_{j-2}$ induces a backedge graph of $T[S\cup\{c\}]$ that is isomorphic to that of $M_{1,1,adj}^{(j+4)}$ in \Cref{def: infinite classes minimal nonmatching}.
                            Since this is a proper subtournament of $T$, as evidenced by the fact that its vertex set does not contain at least one of $c_1,c_2$, that has no matching backedge graph, we obtain a contradiction to minimality of $T$. See Figure \ref{fig:claim8-4}.\qedhere
                    \end{itemize}
                \end{proof}

                \begin{figure}[ht]
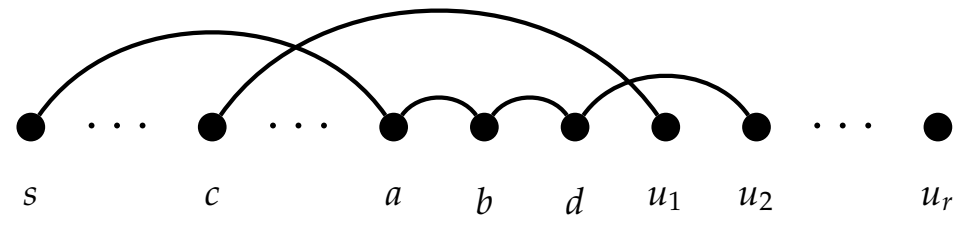
 \centering
                   \BackedgeGraphGeneral{1/3, 4/5, 3/4, 5/7,  2/6}{8}{$s$, $c$, $a$, $b$, $d$, $u_1$, $u_2$, $u_r$}{1,2,7}
                   \caption{Option 4 in the proof of Claim \ref{claim: matching proof case 2.2. odd i points at i+3}, ordering $\sigma_{-1}'$.}
                   \label{fig:claim8-4}
                   \end{figure}
                        
                Finally, we claim that $u_I$ neighbours $u_{I+4}$ in $B_{\sigma}$.
                If this is true, a similar argument as in Option 4 above rules out backedges between $u_{I+2}, u_{I+3}$ and $c_1, c_2$, and thus the ordering $\sigma'_{I}$ induces a backedge graph of $T[s,a,b,c_1,c_2,d,u_1,\dots,u_{I+4}]$ that is isomorphic to that of $M_{2,2}^{(I+6)}$ in \Cref{def: infinite classes minimal nonmatching},
                By minimality of $T$, it must be isomorphic to that tournament.
                \begin{claim}
                    The unique neighbour of $u_I$ in $B_{\sigma}$ is $u_{I+4}$. \label{cl:10}
                \end{claim}
                
                \begin{proof}
                    By \Cref{claim: matching proof case 2.2. odd i points at i+3}, the only remaining possibility is for $u_{I}u_{I+3}$ to be an edge in $B_{\sigma}$. See Figure \ref{fig:claim10}.
                    \begin{figure}[ht]
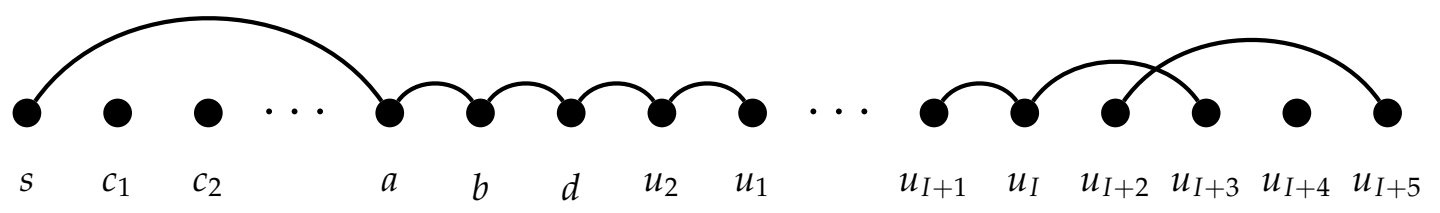
 \centering
                   \BackedgeGraphGeneral{1/4, 4/5, 9/10, 5/6,  6/7, 7/8, 10/12, 11/14}{14}{$s$, $c_1$, $c_2$, $a$, $b$, $d$, $u_2$, $u_1$, $u_{I+1}$, $u_I$,  $u_{I+2}$, $u_{I+3}$, $u_{I+4}$, $u_{I+5}$}{3,8}
                   \caption{Proof of Claim \ref{cl:10}, ordering $\sigma_{I}'$.}
                   \label{fig:claim10}
                   \end{figure}
                    Moreover, the proof of \Cref{claim: matching proof case 2.2. odd i points at i+3} still goes through for $u_{I+2}$ (as it is impossible for $u_{I+2}$ to neighbour $u_{I+2+4}$ in $B_{\sigma}$, i.e., the only potentially problematic case, Option 2 in that proof, never occurs).
                    Hence, $u_{I+2}$ must neighbour $u_{I+5}$ in $B_{\sigma}$.
                    We consider the neighbourhood of $u_{I+4}$ and note that $u_{I+4}$ is either isolated, or exclusively contains a $c$ with $s\prec_{\sigma}c\prec_{\sigma}a$ in $B_{\sigma}$.
                    In either case, the considerations in Option 1 or 4 of the proof of \Cref{claim: matching proof case 2.2. odd i points at i+3} respectively apply to $u_{I+4}$, leading to contradictions.
                    It follows that $u_{I}u_{I+3}$ cannot be an edge of $B_{\sigma}$ and therefore $u_I$ must neighbour $u_{I+4}$ in $B_{\sigma}$.
                \end{proof}
                This concludes our examination of Sub-case 2.2.
                
                % \Dinis{The above tells us that we find the tournaments of type a) \emph{twice} in our proof (here and in sub-case 1.1.).
                % Do the cases really not reduce to one another (I don't see how they can, we have the whole weirdness surrounding $s$ which seems to prevent it, but it is a shame that the proof ends up being this long and redundant)}
            
            \item [Sub-case 2.3.] There is precisely one vertex $c$ between $s$ and $a$ along $\sigma$, and precisely one vertex $d$ between $a$ and $b$ along $\sigma$.
                If $c$ or $d$ have no neighbours in $B_{\sigma}$, then we can find matching orderings of $T$ (by swapping $s$ and $c$, and then $s$ and $a$, or $a$ and $d$, and then $a$ and $b$).
                Moreover, their neighbours must be distinct and come after $b$ along $\sigma$.
                In particular, at least one vertex of $T$ comes after $b$ in $\sigma$.
                Throughout the analysis of this case, we will use the ordering $\sigma'$ which is obtained from $\sigma$ by swapping $d$ and $b$.
                Let $r$ be the number of vertices that are after $b$ in $\sigma$ and let these vertices be $u_1,\dots,u_r$ sorted increasingly along $\sigma$. See Figure \ref{fig:case23}. 
\begin{figure}[ht]
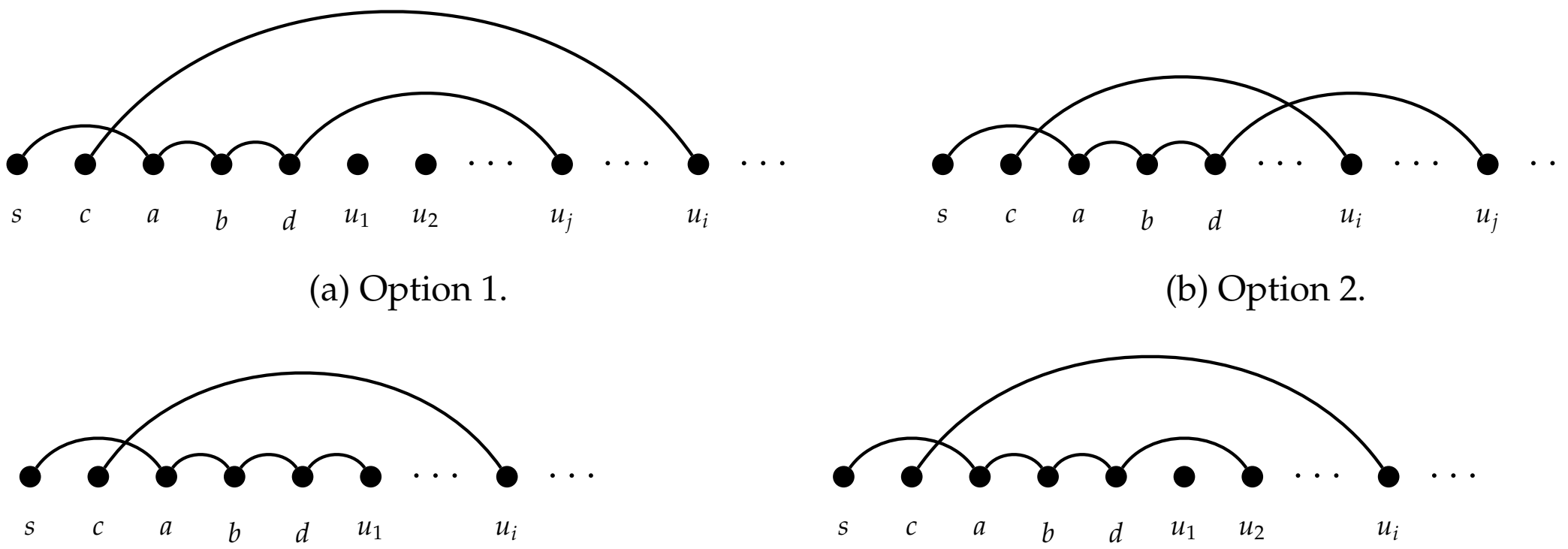

    \centering
    
    \begin{subfigure}{0.5\linewidth}\centering
        \BackedgeGraphGeneral{1/3, 3/4, 4/5, 2/9, 5/8}{9}{$s$, $c$, $a$, $b$, $d$, $u_1$, $u_2$, $u_j$, $u_i$}{7, 8, 9}
        \caption{Option 1.}
    \end{subfigure}
    \begin{subfigure}{0.4\linewidth}\centering
          \BackedgeGraphGeneral{1/3, 3/4, 4/5, 5/7, 2/6}{7}{$s$, $c$, $a$, $b$, $d$,  $u_i$, $u_j$}{5, 6, 7}
        \caption{Option 2.}
    \end{subfigure}

    \begin{subfigure}{0.4\linewidth}\centering
        \BackedgeGraphGeneral{1/3, 3/4, 4/5, 2/7, 5/6}{7}{$s$, $c$, $a$, $b$, $d$, $u_1$, $u_i$}{6, 7}
        \caption{Option 3.}
    \end{subfigure}\begin{subfigure}{0.5\linewidth}\centering
        \BackedgeGraphGeneral{1/3, 3/4, 4/5, 2/8, 5/7}{8}{$s$, $c$, $a$, $b$, $d$, $u_1$, $u_2$, $u_i$}{7, 8}
        \caption{The outcome asserted by Claim \ref{cl:11}.}
    \end{subfigure}
    
        \caption{Sub-case 2.3, ordering $\sigma'$.}
    \label{fig:case23}
\end{figure}
                
                \begin{claim} \label{cl:11}
                    Let $1\leq i,j\leq r$ be such that $u_i$ is the unique neighbour of $c$ in $B_{\sigma}$ and $u_j$ is the unique neighbour of $d$ in $B_{\sigma}$.
                    If $T$ is not a tournament with no matching ordering and at most nine vertices, then $j=2$ and $i>j$.
                \end{claim}
                \begin{proof}
                    By the argumentation from above, there exist $i$ and $j$ such that $c$ neighbours $u_i$ and $d$ neighbours $u_j$ in $B_{\sigma}$.
                    Hence, it suffices to prove that for such $i$ and $j$ we must have $j=2$ and $i>j$ (unless $T$ is a small tournament with no matching ordering).
                    To do so, we consider different cases depending on the values of $i$ and $j$.
                    \begin{itemize}[align=left,leftmargin=*]
                        \item [Option 1:] $j\geq 3$ and $i>j$.
                            Then the ordering of $\{s,a,b,c,d,u_{j-2},u_{j-1},u_j,u_i\}$ induced by $\sigma'$ induces a backedge graph of the corresponding subtournament of $T$ that is isomorphic to that of $M_{1,2}^{(5)}$ in \Cref{def: infinite classes minimal nonmatching} (up to swapping $u_{j-2}$ and $u_{j-1}$ if they are adjacent in $B_{\sigma}$).
                            Minimality implies that $T$ is this tournament.
                        \item [Option 2:] $j\geq 2$ and $i<j$.
                            Then the ordering of $\{s,a,b,c,d,u_{i},u_{j}\}$ induced by $\sigma'$ induces a backedge graph of the corresponding subtournament of $T$ that is isomorphic to that of $M_{1,1,adj}^{(5)}$ in \Cref{def: infinite classes minimal nonmatching}.
                            Minimality implies that $T$ is this tournament.
                        \item [Option 3:] $j=1$.
                            Then the ordering obtained from $\sigma'$ by swapping $a$ with $b$, and swapping $d$ with $u_1$, is a matching ordering of $T$.\qedhere
                    \end{itemize}
                \end{proof}
                We may therefore assume that the only neighbour of $d$ in $B_{\sigma}$ is $u_2$ and the only neighbour of $c$ in $B_{\sigma}$ is $u_i$ with $i>2$.
                Let $I$ be the largest odd number that is at most $r-5$. 
                \begin{claim}
                    \label{claim: matching proof case 2.3. odd i points at i+3}
                    For each odd $j$ between $1$ and $I-1$, the unique neighbour of $u_j$ in $B_{\sigma}$ is $u_{k+3}$.
                    Moreover, the unique neighbour of $c$ in $B_{\sigma}$, $u_i$, has $i>I$, and $u_I$ is neither isolated, nor does it neighbour $u_{I+1}$, $u_{I+l}$ with $l\geq 5$, $u_{I+2}$ in $B_{\sigma}$.
                \end{claim}
                \begin{figure}
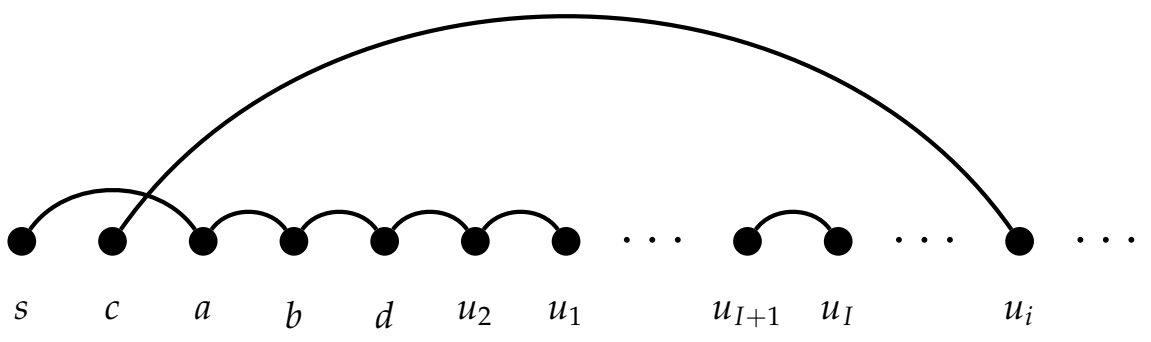

                \centering\BackedgeGraphGeneral{1/3, 3/4, 4/5, 2/10, 5/6, 6/7, 8/9}{10}{$s$, $c$, $a$, $b$, $d$, $u_2$, $u_1$, $u_{I+1}$, $u_I$, $u_i$}{7, 9, 10}
        \caption{The outcome asserted by Claim \ref{claim: matching proof case 2.3. odd i points at i+3} with ordering $\sigma_I'$.}
        \label{fig:27}
            \end{figure}
                \begin{proof}
                    For every odd $k$ between $1$ and $r$, we inductively define $\sigma'_k$ to be the ordering obtained from $\sigma'_{k-2}$ (with $\sigma'_{-1}=\sigma'$) by swapping $u_k$ and $u_{k+1}$. See Figure \ref{fig:27}. 

                    Assume, for the sake of contradiction, that the claim is false, and let $j$ be the smallest odd number falsifying it.
                    Since $j$ is smallest, the ordering $\sigma'_{j-2}$ induces a backedge graph of $T[s,a,b,d,u_1,\dots,u_{j-1}]$ that is isomorphic to a path on $j+2$ vertices.
                    The following also turns out to be useful.
                    \begin{observation}
                        \label{obs: vertex sets closed under nbhd case 2.3}
                        All neighbours  of vertices in $S=\{s, a, b, d,u_1,\dots,u_{j+1}\}\cup N_{B_{\sigma}}(u_j)$ in the backedge graphs $B_{\sigma}$, $B_{\sigma'}$, and $B_{\sigma'_{j-2}}$ are in $S$ and an analogous statement applies to $V(T)\setminus S$.
                    \end{observation}
                    Also note that $i\geq j$, since $i$ cannot be an odd number smaller than $j$, and all $u_l$ with even $l<j$ have a neighbour in $\{d,u_1,\dots,u_{j-3},u_{j-1}\}$.
                    We consider four possibilities for the neighbourhood of $u_j$ in $B_{\sigma}$:
                    \begin{itemize}[align=left,leftmargin=*]
                        \item [Option 1:] $u_j$ is isolated in $B_{\sigma}$.
                            Note that $T[S]$ has a backedge graph, $B_{\sigma'_j}$, that is an induced monotone path on $j+5$ vertices.
                            By \Cref{lemma: matching orderings of induced paths}, there is an ordering $\sigma''$ of $S$ that induces a matching backedge graph of this subtournament and whose first vertex is $s$. 
                            The ordering of $T$ that induces the same ordering of $S$ as $\sigma''$ and where the remaining vertices occupy the same place as in $\sigma$ is a matching ordering, a contradiction.
                        
                        \item [Option 2:] $u_j$ neighbours $u_l$ with $l\geq j+4$ in $B_{\sigma}$ (or $u_{l}$ with $l\geq j+5$ if $j=I$).
                            If the neighbour of $c$, $u_i$, is before $u_l$ in $\sigma$, then the ordering of $S\cup\{c,u_i\}$ induced by $\sigma'_{j}$ induces a backedge graph of $T[S\cup\{c,u_i\}]$ that is isomorphic to that of $M_{1,1,adj}^{(j+6)}$ in \Cref{def: infinite classes minimal nonmatching}.
                            Otherwise, the ordering of $S\cup\{c,u_i,u_{l-2},u_{l-1}\}$ induced by $\sigma'_{j}$ induces a backedge graph of $T[S\cup\{a,b,c,s,u_i,u_{l-2},u_{l-1}\}]$ that is isomorphic to that of $M_{2,2}^{(j+6)}$ in \Cref{def: infinite classes minimal nonmatching} (up to swapping $u_{l-2}$ and $u_{l-1}$ if they are adjacent in $B_{\sigma}$).
                            In either case, we have found a proper subtournament of $T$, as evidenced by the fact that its vertex set does not contain $u_{l-3}$ (when $l\geq j+5$) or $u_{l+1}$ (when $l=j+4$) that has no matching backedge graph, contradicting minimality of $T$.
                        
                        \item [Option 3:] $u_j$ neighbours $u_{j+2}$ in $B_{\sigma}$.
                            Note that $T[\{a,b,d,u_1,\dots,u_{j+2}\}]$ has a backedge graph, $B_{\sigma'_j}$, that is an induced monotone path on $j+5$ vertices.
                            By \Cref{lemma: matching orderings of induced paths}, there is an ordering $\sigma'''$ of $\{a,b,d,u_1,\dots,u_{j+1}\}$ that induces a matching backedge graph of this subtournament of $T$ where $a$ has degree zero. The ordering of $T$ that induces the same ordering of $\{a,b,d,u_1,\dots,u_{j+2}\}$ as $\sigma''$ and where the remaining vertices in occupy the same place as in $\sigma$ is a matching ordering, a contradiction.
                            
                        \item [Option 4:] $u_j$ neighbours $c$ in $B_{\sigma}$.
                            Then the ordering of $S$ induced by $\sigma'_{j-2}$ induces a backedge graph of $T[S]$ that is isomorphic to that of $M_{1,1,adj}^{(j+4)}$ in \Cref{def: infinite classes minimal nonmatching}.
                            Since this is a proper subtournament of $T$ that has no matching backedge graph, , as evidenced by the fact that its vertex set does not contain $u_{j+2}$, we obtain a contradiction to minimality of $T$.\qedhere
                    \end{itemize}
                \end{proof}
                Now we pay closer attention to the neighbourhood of $u_I$ in $B_{\sigma}$.
                \begin{claim}
                    \label{claim: matching proof case 2.3 large i}
                    One of the following must be true:
                    \begin{itemize}
                        \item $u_{I}$ neighbours $u_{I+4}$ and $u_i$ has $i=I+5$.
                        \item $u_{I}$ neighbours $u_{I+3}$, $u_{I+2}$ neighbours $u_{I+5}$, and the unique neighbour of $c$, $u_{i}$, has $i=I+4$.
                    \end{itemize}
                    \label{cl:13}
                \end{claim}
                \begin{proof}
                    Note that by \Cref{claim: matching proof case 2.3. odd i points at i+3}, $u_I$ neighbours $u_{I+3}$ or $u_{I+4}$ in $B_{\sigma}$.
                    We consider each of these possibilities separately. See Figure \ref{fig:28}. 

        \begin{figure}
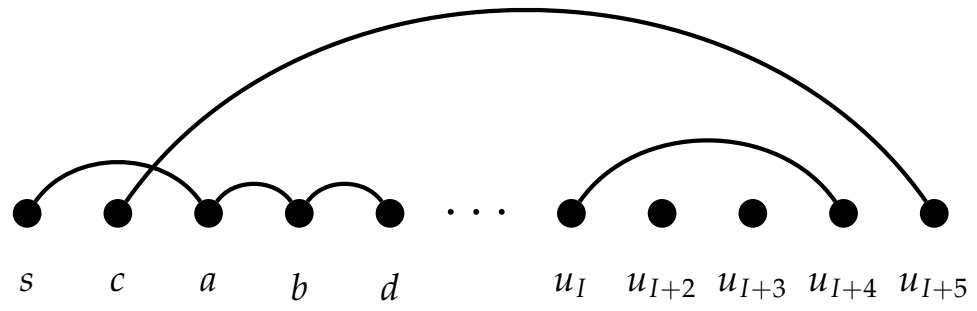
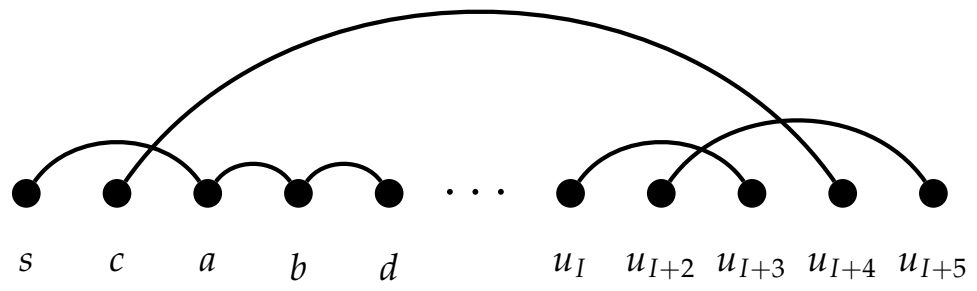

            \centering
            \begin{subfigure}{0.48\linewidth}
                \centering\BackedgeGraphGeneral{1/3, 3/4, 4/5, 6/9, 2/10}{10}{$s$, $c$, $a$, $b$, $d$,  $u_I$, $u_{I+2}$, $u_{I+3}$, $u_{I+4}$, $u_{I+5}$}{5}
        \caption{Option 1 in the proof of Claim \ref{cl:13}.}
            \end{subfigure}
            \begin{subfigure}{0.48\linewidth}
                \centering\BackedgeGraphGeneral{1/3, 3/4, 4/5, 6/8, 7/10, 2/9}{10}{$s$, $c$, $a$, $b$, $d$,  $u_I$, $u_{I+2}$, $u_{I+3}$, $u_{I+4}$, $u_{I+5}$}{5}
        \caption{Option 2 in the proof of Claim \ref{cl:13}.}
            \end{subfigure}
                \caption{Proof of Claim \ref{cl:13}.}
        \label{fig:28}
            \end{figure}
                    \begin{itemize}[align=left,leftmargin=*]
                        \item [Option 1:] $u_I$ neighbours $u_{I+4}$ in $B_{\sigma}$.
                            If $i\neq I+5$, then we will now show we can always derive a contradiction.
                            We consider two different types of value $i$ can take:
                            \begin{itemize}[align=left,leftmargin=*]
                                \item [Sub-option 1.1:] $I+2\leq i\leq I+3$.
                                    Then the ordering of $\{s,a,b,c,d,u_1,\dots,u_{I+1},u_i,u_{I+4}\}$ induced by $\sigma'_{I}$ induces a backedge graph of $T[\{s,a,b,c,d,u_1,\dots,u_{I+1},u_i,u_{I+4}\}]$ that is isomorphic to that of $M_{1,1,adj}^{(I+6)}$ in \Cref{def: infinite classes minimal nonmatching}.
                                    Since this is a proper subtournament of $T$, as evidenced by the fact that its vertex set does not contain $u_{I+5}$, that has no matching backedge graph, we obtain a contradiction to minimality of $T$.
                                    
                                \item [Sub-option 1.2:] $i=I+6$.
                                    Then the ordering of $\{s,a,b,c,d,u_1,\dots,u_{I+4},u_{i}\}$ induced by $\sigma'_{I}$ (or $\sigma'_{I+2}$ if $u_{I+2}u_{I+3}\in E(B_{\sigma})$) induces a backedge graph of $T[\{s,a,b,c,d,u_1,\dots,u_{I+4},u_{i}\}]$ that is isomorphic to that of $M_{1,2}^{(I+6)}$ in \Cref{def: infinite classes minimal nonmatching}.
                                    Since this is a proper subtournament of $T$, as evidenced by the fact that its vertex set does not contain $u_{I+5}$, that has no matching backedge graph, we obtain a contradiction to minimality of $T$.
                            \end{itemize}
                            Hence, $i=I+5$ as we set out to prove.

                        \item [Option 2:] $u_{I}$ neighbours $u_{I+3}$ in $B_{\sigma}$.
                            Note that the considerations in Options 1, 3, and 4 of the proof of \Cref{claim: matching proof case 2.3. odd i points at i+3} still go through for $u_{I+2}$, enabling us to conclude that $u_{I+2}$ neighbours $u_{I+k}$ with $j\geq 4$ in $B_{\sigma}$.
                            We prove that $k=5$ by considering each other option separately.
                            \begin{itemize}[align=left,leftmargin=*]
                                \item [Sub-option 2.1:] $u_{I+2}$ neighbours $u_{I+4}$ in $B_{\sigma}$.
                                    Note that, much as in Option 1 of the proof of \Cref{claim: matching proof case 2.3. odd i points at i+3}, $T[\{a,b,d,u_1,\dots,u_{I+4}\}]$ has a backedge graph, $B_{\sigma'_{I+2}}$, that is an induced monotone path on $j+7$ vertices.
                                    Precisely as in that proof, this implies the existence of a matching ordering of $T$, a contradiction.
    
                                \item [Sub-option 2.2:] $u_{I+2}$ neighbours $u_{I+6}$.
                                    Then the ordering of $\{s,a,b,c,d,u_1,\dots,u_{I+3},u_{i},u_{I+6}\}$ induced by $\sigma'_{I+2}$ induces a backedge graph of $T[\{s,a,b,c,d,u_1,\dots,u_{I+3},u_{i},u_{I+6}\}]$ that is isomorphic to that of $M_{1,1,adj}^{(I+8)}$ in \Cref{def: infinite classes minimal nonmatching}.
                                    Since this is a proper subtournament of $T$, as evidenced by the fact that its vertex set does not contain $u_{I+4}$ or $u_{I+5}$ (only one of these can be $u_i$), that has no matching backedge graph, we obtain a contradiction to minimality of $T$.
                            \end{itemize}
                            Hence, $u_{I+2}$ neighbours $u_{I+5}$.
                            All that remains to be shown is that $c$ neighbours $u_{I+4}$ in $B_{\sigma}$ (equiv. $i=I+4$).
                            Aiming for a contradiction, assume otherwise.
                            Then $i=I+6$.
                            This implies that $T[\{a,b,d,u_1,\dots,u_{I+5}\}]$ has a backedge graph, $B_{\sigma'_{I+4}}$, that is an induced monotone path on $j+7$ vertices.
                            Much as in Sub-option 2.2, this implies the existence of a matching ordering of $T$, a contradiction.
                            % \Dinis{Maybe it is no longer clear, but this only works because the only edge going into the vertex set of this tournament is between $s$ and $a$ (and that fails for $i=I+4$).
                            % This should probably be pointed out.}
                    \end{itemize}
                \end{proof}
                In the cases that were not excluded in \Cref{claim: matching proof case 2.3 large i}, it is not difficult to see that $T[s,a,b,c,d,u_1,\dots,u_{I+5}]$ is a subtournament of $T$ in the forbidden set ($M_{1,2}^{(6)}$ if $u_Iu_{I+4}\in E(B_{\sigma})$, and $M_{1,1,adj}^{(I+8)}$ if $u_Iu_{I+3}\in E(B_{\sigma})$).
                Minimality of $T$ implies that $T$ is, in fact, equal to this subtournament, concluding our examination of Sub-case 2.3.
        \end{itemize}
\end{itemize}
Since in each case $T$ was shown to be a tournament in the forbidden set (and $T$ must fall into one of the categories described by the analyzed cases), $T$ is always in the forbidden set.

\section{Non-crossing matching tournaments}\label{section:NCM}
In this section, we prove the characterization of NCM tournament announced earlier.
\ncm

We will use the following theorem due to Latka \cite{https://doi.org/10.1002/jgt.10081}. 
\begin{theorem}\label{thm:latka}
    A prime tournament on at least five vertices contains $\Delta(1,2,2), T_5, U_5$ or $W_5$.
\end{theorem}
Observe the similarities between the above two theorems. 

\begin{lemma}\label{lem:NCM-suff}
    The tournaments $\Delta(1,2,2), T_5, U_5,$ and $W_5$ are not NCM.
\end{lemma}
\begin{proof}
    Any tournament having a matching backedge graph has a vertex of in-degree at most 1. $T_5$ is 2-regular, and thus, not NCM. The verification that the remaining three tournaments can be done by brute force, and is not very instructive. We omit the details.
\end{proof}

A set $H$ in a tournament $T$ is \emph{homogeneous} if $H\neq V(T)$ and for all $v\in V(T)\setminus H$ either $v\Rightarrow H$ or $H\Rightarrow v$.
After \cite{boussairi2004c3}, we call a homogeneous set $H$ in a tournament $T$ \emph{strong} if for any homogeneous set $H'$ in $T$,  $H\cap H'\neq\emptyset$ implies one of $H$ and $H'$ is contained in the other. The strong homogeneous sets of a tournament form a laminar family. Moreover, since every vertex by itself is a strong homogeneous set, the set $\mathcal{P}(T)$ of maximal strong homogeneous sets of $T$ is a partition of $V(T)$. Let $Q(T)$ be the tournament on $\mathcal{P}(T)$ where $X\Rightarrow Y$ if for $x\in X$ and $y\in Y$, $x\Rightarrow y$ in $T$. 

\begin{theorem}[Theorem 4 in \cite{boussairi2004c3}]\label{thm:quotient}
    Let $T$ be a tournament.
    \begin{enumerate}
        \item $T$ is not strongly connected if and only if $Q(T)$ is a transitive tournament.
        \item $T$ is strongly connected if and only if $Q(T)$ has at least three vertices and is prime.
    \end{enumerate}
\end{theorem}

\begin{proof}[Proof of \Cref{thm:NCM}]
    By \Cref{lem:NCM-suff} it suffices to prove that if $T$ is not NCM, then $T$ contains one of $\Delta(1,2,2), T_5, U_5,$ and $W_5$.    
    Suppose this is false, and let $T$ be a minimal counterexample, meaning $T$ is not NCM, does not contain $\Delta(1,2,2),$ $T_5,$ $U_5,$ or $W_5$, and every proper subtournament of $T$ is NCM. Every tournament on at most four vertices is NCM, so $T$ has at least five vertices. By \Cref{thm:latka}, $T$ is not prime. Consider $Q = Q(T)$. By \Cref{thm:quotient}, $Q$ is prime or transitive. If $Q$ is transitive, we immediately obtain a contradiction by minimality: take an NCM ordering of each homogeneous set in $\mathcal{P}(T)$, and combine them in the natural way. Thus $Q$ is prime and not transitive. If $Q$ has at least five vertices, we obtain by \Cref{thm:latka} one of the four desired tournaments as a subtournament of $Q$ and thus of $T$. Thus $Q$ has at most four vertices. There are four tournaments on four vertices up to isomorphism, none of which are prime, and all tournaments on at most two vertices are transitive. Thus $Q$ is isomorphic to $C_3$. Then $\mathcal{P}(T)$ contains exactly one non-trivial homogeneous set, say $S$, as otherwise $T$ contains $\Delta(1,2,2)$. So $T$ is of the form $\Delta(1,1,S)$. By minimality of $T$, we conclude that $S$ is NCM. But then it is easy to obtain an NCM backedge graph for $T$ from that of $S$, a contradiction.
\end{proof}

\section{Future Work}
\label{sec: conclusion}
In the preceding text we have precisely characterized the sets of minimal forbidden subtournaments for membership in three classes of tournaments defined in terms hereditary properties of their backedge graphs: matching tournaments, paving tournaments, and non-crossing matching tournaments.
These characterizations give rise to some open questions that we collect below.

Our characterizations of paving and non-crossing matching tournaments both imply the existence of a polynomial-time algorithm for recognition of tournaments in these classes (respectively with runtime $O(n^6)$ and $O(n^5)$).
It is natural to wonder whether substantially faster algorithms exist. \cite{keeney2024degreewidth} shows that matching tournaments can be recognized in $O(n^2)$ time. Generalizing \cite{keeney2024degreewidth}, the work \cite{evangelista2026backedge} gives an $O(n^2)$ algorithm for recognizing $n$-vertex paving tournaments. 
Can one check whether a given tournament is NCM in time $O(n^2)$?

It also seems like each of the classes we examined here contain \enquote{very few} (non-isomorphic) tournaments (where by very few, we mean a slowly growing function of their number of vertices).
How many matching/paving/NCM tournaments are there (asymptotically)?

We conclude with what seems to be the most interesting question suggested by our results.
We have precisely characterized the set of minimal forbidden tournaments for membership in $\mathcal{C}$-FAS for three classes hereditary classes $\mathcal{C}$.
All three classes $\mathcal{C}$ that we considered can be defined in terms of finitely many forbidden subgraphs.
However, for two of these classes, we have found finitely many forbidden subtournaments, whereas for the class of matching tournaments the set of minimal forbidden subtournaments is infinite.
% I vaguely recall something about "always a finite set" for having a $k$-degenerate backedge graph being mentioned in some meeting? If true, it might be worth bringing up here
Precisely determining which hereditary classes $\mathcal{C}$ result in $\mathcal{C}$-FAS being characterizable by finitely many forbidden tournaments seems widely out of reach.
It is not even known which hereditary classes admit characterizations in terms of finitely many forbidden induced subgraphs, although this question was raised as early as \cite{damaschke1990}.
In a similar spirit, Latka asks which classes of tournaments are well-quasi-ordered by induced subtournaments \cite{latka_antichains_bounded_size}, and here too the situation seems to not be very well-understood.
In this context, it is worth pointing out that any class of tournaments that is closed under subtournaments and admits a characterization in terms of finitely many forbidden subtournaments also admits a characterization in terms of finitely many forbidden backedge graphs, as forbidding a tournament amounts to the same as forbidding all its (finitely many) backedge graphs.
We conjecture that the following holds:
\begin{conjecture}
    \label{conj: forbidding only one backedge graph gives infinitely many tournaments}
    For any backedge graph $B$ on at least three vertices, the class of minimal tournaments $T$ such that all backedge graphs of $T$ contain $B$ as an induced (ordered) subgraph is infinite.
\end{conjecture}
Our main reason for believing in \Cref{conj: forbidding only one backedge graph gives infinitely many tournaments} is that it is true whenever $B$ is a clique (unless P$=$NP) due to the NP-hardness results mentioned in the introduction.
For similar reasons, we conjecture that the same is true when we forbid $B$ as an unordered graph.
As an initial step towards resolving \Cref{conj: forbidding only one backedge graph gives infinitely many tournaments}, it may be useful to first identify infinite antichains in the class of minimal forbidden subtournaments for $K_3$-free-FAS.
There is some hope that it may be possible to adapt such constructions and obtain other antichains of tournaments that are minimal such that all their backedge graphs contain other graphs (on three vertices).

\section*{Comment on the use of AI:}
AI was not used at any point during the preparation of this document.
All the research, mathematical ideas and arguments, writing, and figures are due to the authors.

\bibliographystyle{abbrv}
\bibliography{bibliography.bib}
\end{document}